\documentclass[preprint,twocolumn,5p]{elsarticle}

\usepackage[utf8]{inputenc}
\usepackage{amsmath}
\usepackage{amsfonts}
\usepackage{amssymb}
\usepackage{amsthm}
\usepackage{booktabs}
\usepackage{graphicx}
\usepackage{subcaption}
\usepackage{hyperref}
\usepackage{mathtools}
\usepackage{microtype}
\usepackage{multirow}
\usepackage{pdfpages}
\usepackage{placeins}
\usepackage{stix2}
\usepackage{tikz}
\usepackage{transparent}
\usepackage{xcolor}
\usepackage[normalem]{ulem} %

\usepackage{cleveref}

\hypersetup{
    colorlinks=true,
    pageanchor=false,
    linkcolor=[rgb]{0.10,0.20,0.55},
    citecolor=[rgb]{0.10,0.20,0.55},
    urlcolor=[rgb]{0.10,0.20,0.55},
}

\newtheorem{definition}{Definition}

\newcommand{\e}[1]{\times 10^{#1}}

\title{Approaching the optimal closure: equivariance, inductive bias, and Reynolds-number generalization in data-driven LES}

\date{\today}

\begin{document}

\begin{frontmatter}

\author[cwi,eindhoven]{Syver Døving Agdestein\corref{cor1}}
\ead{sda@cwi.nl}
\author[cwi,eindhoven]{Benjamin Sanderse}
\ead{b.sanderse@cwi.nl}

\cortext[cor1]{Corresponding author}

\affiliation[cwi]{
    organization={Scientific Computing Group, Centrum Wiskunde \& Informatica},
    addressline={Science Park 123}, 
    city={Amsterdam},
    postcode={1098 XG}, 
    country={The Netherlands}
}
\affiliation[eindhoven]{
    organization={Centre for Analysis, Scientific Computing and Applications, Eindhoven University of Technology},
    addressline={PO Box 513}, 
    city={Eindhoven},
    postcode={5600 MB}, 
    country={The Netherlands}
}

\begin{abstract}
Data-driven closures for large-eddy simulation (LES) are commonly built to
respect the symmetries of the Navier--Stokes equations, on the premise that
this improves accuracy and generalization. We test this premise in a
controlled comparison of three data-driven LES closures that share a
pointwise, Galilean-invariant velocity-gradient construction but span
non-equivariant, octahedral-equivariant, and tensor-basis designs: an
unconstrained multi-layer perceptron (MLP), a group-convolutional network
whose exactly equivariant weights we synthesize in closed form, and a
tensor-basis neural network (TBNN). The designs follow from an analysis of
which symmetries survive discretization on a uniform grid, where the
continuous orthogonal group reduces to the 48-element octahedral group. Across a range of network sizes the three closures
saturate to the same a priori and a posteriori accuracy, and a direct
conditional-mean estimate identifies the a priori floor as the one-point
optimal closure of Langford and Moser. The equivariant and tensor-basis
models reach this floor with $25$ times fewer parameters than the MLP: the
inductive bias buys parameter efficiency rather than a lower error floor.
Finally, we train the closures across several viscosities and supply the
global filter-scale Reynolds number
$\operatorname{Re}_\Delta = \Delta^2 \| \nabla \bar{u} \| / \nu$ as an input,
a scaling-invariant feature dictated by the same symmetry
analysis. The closures then generalize across Reynolds number: they hold
their dissipation calibration at held-out viscosities and filter ratios where
Reynolds-blind closures mis-dissipate, and partially correct it on an
out-of-distribution Taylor--Green flow. Reynolds-number generalization is
thus largely a calibration that the right input feature supplies.

\end{abstract}

\begin{graphicalabstract}
\begin{center}
\includegraphics[width=\textwidth]{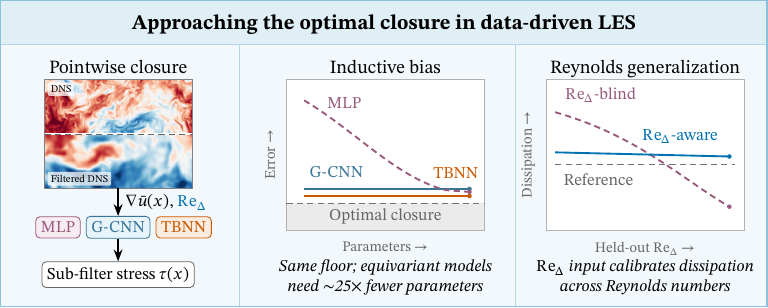}
\end{center}
\end{graphicalabstract}

\begin{highlights}
\item We determine which Navier--Stokes symmetries survive discretization in forced LES.
\item We build exactly equivariant octahedral group-convolution weights in closed form.
\item Data-driven closures saturate at the directly estimated one-point optimal closure.
\item Equivariant and tensor-basis closures reach the floor at 25 times fewer parameters.
\item The filter-scale Reynolds number as an input restores the dissipation calibration.
\end{highlights}


\end{frontmatter}

\section{Introduction} \label{sec:introduction}

Direct numerical simulation (DNS) of turbulent flows resolves all scales of motion,
but its computational cost often exceeds available resources.
Large-eddy simulation (LES) offers a tractable alternative by resolving only
the large-scale features and modeling the effect of unresolved
scales through a closure model~\cite{popeTenQuestionsConcerning2004}.

Symmetries are fundamental properties of both
turbulence~\cite{frischTurbulenceLegacyKolmogorov1995} and differential
equations in
general~\cite{blumanSymmetriesDifferentialEquations1989,olverApplicationsLieGroups1986}.
The incompressible Navier--Stokes equations
exhibit symmetries including Galilean invariance and rotation
invariance.
The LES equations, obtained by filtering the Navier--Stokes equations,
may inherit these symmetries depending on the filter,
but the complete LES system preserves them only if the closure model does as well.
Closure models that violate these symmetries can introduce spurious forces
and unphysical behavior~\cite{oberlackInvariantModelingLargeeddy1997,
razafindralandyAnalysisDevelopmentSubgrid2007,
silvisPhysicalConsistencySubgridscale2017}.
Classical functional closures such as eddy-viscosity models enforce
Galilean invariance by construction, since they depend only on invariants of the
velocity-gradient tensor.
Structural closures, including Bardina's scale-similarity
model~\cite{bardinaImprovedSubgridscaleModels1980}
and
Clark's gradient model~\cite{clarkEvaluationSubgridscaleModels1979},
also preserve most symmetries of the original
equations~\cite{silvisPhysicalConsistencySubgridscale2017}.
Beyond these classical models, a dedicated class of subgrid closures has been
constructed to explicitly preserve the Lie symmetries of the Navier--Stokes
equations~\cite{razafindralandyAnalysisDevelopmentSubgrid2007,
akhalPrioriAnalysisSubgrid2026,
cosseratVorticitydependentSymmetrypreservingModels2026}.

Machine learning (ML) has emerged as a powerful tool for turbulence closure
modeling~\cite{duraisamyTurbulenceModelingAge2019,
sanderseScientificMachineLearning2025}.
Early data-driven LES closures learned the subgrid term directly from
high-fidelity data~\cite{maulikSubgridModellingTwodimensional2019,
beckDeepNeuralNetworks2019},
and ML has since been applied to both
functional~\cite{kurzHarnessingEquivarianceModeling2025}
and structural~\cite{prakashInvariantDatadrivenSubgrid2022,
prakashInvariantDatadrivenSubgrid2024}
LES closure modeling, achieving high a priori accuracy.
However, models that perform well on training data
can still exhibit instabilities when deployed in
actual simulations~\cite{kurzInvestigatingModelDataInconsistency2021},
a problem that affects both ML-based and classical
closures~\cite{torresRationalBooleanStabilization2023}.
Symmetry preservation addresses one source of such instabilities:
closure models that respect the symmetries of the Navier--Stokes equations
cannot introduce spurious frame-dependent forces or violate rotational isotropy,
improving physical consistency; symmetry-constrained closures have also been
reported to improve numerical
stability~\cite{guanLearningPhysicsconstrainedSubgridscale2023}.
Ensuring symmetry preservation in ML closures is thus an important step toward
their reliable deployment~\cite{mcconkeyTurbulenceTeachesEquivariance2026a}.

Several approaches have been developed to embed symmetries in ML-based closure models.
One family of methods employs group-equivariant neural
networks~\cite{cohenGroupEquivariantConvolutional2016,
lafargeRototranslationEquivariantConvolutional2021}.
Group convolutional neural networks equivariant to the discrete rotation
groups of the cube were developed for volumetric vision
tasks~\cite{worrallCubeNetEquivariance3D2018,
winkelsPulmonaryNoduleDetection2019}---the weight-sharing pattern our
octahedral G-CNN instantiates for stress-tensor-valued fields---and have
been applied to LES~\cite{guanLearningPhysicsconstrainedSubgridscale2023}
and to mesoscale-eddy closures in ocean
models~\cite{perezhoginGeneralizableNeuralNetworkParameterization2025},
but their feature channels scale with the number of group elements,
making them computationally intractable for
large symmetry groups such as the continuous rotation group.
Steerable CNNs address this scaling limitation by representing
equivariance-preserving weights via a truncated Fourier
series~\cite{weilerGeneralE2equivariantSteerable2019,
weilerEquivariantCoordinateIndependent2023}.
Connolly et al.\ applied rotation-equivariant steerable CNNs to atmospheric
LES~\cite{connollyDeepLearningTurbulence2025},
though only for rotations around the vertical axis.
Jalaali and Okabayashi~\cite{jalaaliPrioriAssessmentRotational2026} achieved
rotational invariance in a CNN-based SGS model by incorporating spatial transformer networks.
Physical invariances can also be embedded directly into CNN closures through
specialized equivariant layers or input--output
transformations~\cite{frezatPhysicalInvarianceNeural2021,
pawarFrameInvariantNeural2023}.
Graph neural networks (GNNs) provide yet another equivariant
framework~\cite{kipfSemiSupervisedClassificationGraph2017}:
Shankar et al.\ used GNNs to preserve symmetries in
2D LES~\cite{shankarImportanceEquivariantInvariant2023},
Kurz et al.\ applied them to learn eddy viscosity coefficients in
3D~\cite{kurzHarnessingEquivarianceModeling2025},
and Lino et al.~\cite{linoMultiscaleRotationequivariantGraph2022}
and List et al.~\cite{listRotationalEquivariantGraph2025}
developed rotation-equivariant GNN architectures for fluid dynamics.

Rather than enforcing equivariance through the network architecture,
an alternative approach encodes symmetries through the choice of input features.
This approach employs Pope's tensor
basis~\cite{popeMoreGeneralEffectiveviscosity1975}, rooted in the theory of
Spencer and Rivlin~\cite{spencerTheoryMatrixPolynomials1958,
spencerFurtherResultsTheory1959}.
Expressing the subgrid-scale stress tensor as a linear
combination of basis tensors that are equivariant to the symmetry groups of the
Navier--Stokes equations automatically enforces the required symmetries,
provided the expansion coefficients depend only on invariants of the velocity-gradient tensor (VGT).
Tensor basis closures for LES remain an active area of
research~\cite{lundParameterizationSubgridscaleStress1993,
silvisPhysicalConsistencySubgridscale2017}.
Ling et al.~\cite{lingReynoldsAveragedTurbulence2016} pioneered the use of
tensor basis neural networks (TBNNs), in which a neural network predicts the invariant
coefficients, initially for Reynolds-averaged
Navier--Stokes (RANS) closures, with subsequent developments
in~\cite{milaniGeneralityTensorBasis2021,
berroneInvariancespreservingVectorBasis2022,
cinnellaDatadrivenTurbulenceModeling2024,
brenerHighlyAccurateStrategy2024,
millerSymmetryAwareReynolds2025}.
TBNNs have also been applied to
LES~\cite{boseInvarianceEmbeddedPhysicsinfused2024,
shankarDifferentiableTurbulenceClosure2025,
wuTwoNeuralNetwork2025}.
Stallcup et al.~\cite{stallcupAdaptiveScaleSimilarClosure2022,
stallcupAdaptiveScaleSimilarClosure2022a}
proposed an alternative tensor basis derived from Smith's
analysis~\cite{smithIsotropicFunctionsSymmetric1971},
later employed in a TBNN by Wu and Lele~\cite{wuTwoNeuralNetwork2025}.
A key advantage of TBNNs is that they enforce symmetries regardless of the
functional form of the invariant coefficients, imposing no architectural
constraints on the neural network. This contrasts with group convolutional
neural networks, steerable CNNs,
and GNNs, where the architecture itself must be designed to preserve
equivariance. Kaszuba et al.~\cite{kaszubaImplicitModelingEquivariant2025}
proposed Euclidean neural networks as an alternative that implicitly encodes
the tensor basis structure.
The same architectural freedom can also be obtained by canonicalization rather
than basis expansion: Prakash et
al.~\cite{prakashInvariantDatadrivenSubgrid2022} represent the input and output
tensors in the filtered strain-rate eigenframe, reducing an invariant one-point
closure to four nondimensional inputs, and later embedded filter anisotropy in
the same framework~\cite{prakashInvariantDatadrivenSubgrid2024}.
The two routes parametrize the same class of equivariant one-point closures;
the eigenframe form trades the tensor basis for frame-fixing conventions
(eigenvalue ordering and a vorticity-based eigenvector orientation) that become
discontinuous at degenerate strain configurations.
Notably, Prakash et al.\ report that a single-hidden-layer network trained on
one snapshot at a single filter width suffices, with further data bringing no
improvement---an observation our capacity study revisits as saturation at the
one-point optimal closure (\cref{sec:results-capacity,sec:discussion}).

In this article, we compare three data-driven LES closures that share the same
pointwise velocity-gradient construction, which is Galilean-invariant and,
through the output prefactor, dimensionally consistent under scaling.
Beyond this shared construction, they represent three distinct strategies for
rotational and reflectional equivariance: tensor-basis neural networks (TBNNs)
and group-convolutional neural networks enforce it, while an unconstrained
pointwise multi-layer perceptron (MLP) does not.
The three are representative non-equivariant, tensor-basis-equivariant, and
octahedral-equivariant closures: they hold the shared construction fixed,
together with the Galilean and scaling symmetries it carries, but they
necessarily differ in their internal representation, output parametrization,
and inductive bias as well as in equivariance.
Rather than asking only whether equivariance matters, we use the symmetry
analysis twice: once to constrain the network architecture, and once to choose
its input features.
The first determines how much inductive bias each design carries; the second,
through the scaling-invariant filter-scale Reynolds number, determines how the
closures respond when the Reynolds number changes.
We therefore study the three designs along two axes that single-configuration
comparisons leave open: network capacity, to separate inductive bias from the
attainable accuracy floor, and Reynolds number, to test generalization beyond
the training regime.
Data-driven closures are still typically trained and evaluated at a single
Reynolds number; where transfer to other Reynolds numbers has been pursued, it
has mostly relied on retraining, such as transfer
learning~\cite{subelDatadrivenSubgridscaleModeling2021,
guanStablePosteriori2D2022}.
A notable exception is the reinforcement-learning closure of Novati et
al.~\cite{novatiAutomatingTurbulenceModelling2021}, which generalizes across
grid sizes and Reynolds numbers, with nondimensionalized flow states doing
the carrying---an early indication that such generalization hinges on the
inputs rather than on retraining. Here, likewise, the Reynolds-number
dependence is carried by an input feature of a single trained closure.
Beyond the comparison itself, this work makes four contributions.
First, we analyze which continuous Navier--Stokes symmetries survive in the
discrete, forced LES equations the closures are embedded in.
Second, we give an explicit, closed-form instance of group-convolutional
weight sharing~\cite{cohenGroupEquivariantConvolutional2016,
worrallCubeNetEquivariance3D2018} for the
octahedral group and its stress-tensor-valued input and output layers, with
weights that are exactly equivariant in floating-point arithmetic and require
neither assembling nor eigendecomposing the group-averaging projector.
Third, a capacity study shows that the equivariant and tensor-basis closures
reach the same accuracy floor as the unconstrained network with a factor of
$25$ fewer parameters; all three saturate to a common floor that a direct,
training-free conditional-mean estimate---the optimal estimator of a priori
analysis~\cite{moreauOptimalEstimationLargeeddy2006}---identifies as the
one-point optimal closure in the sense of Langford and
Moser~\cite{langfordOptimalFormulationsIsotropic1999}.
Fourth, a Reynolds-number generalization study shows that the filter-scale
Reynolds number, supplied as a network input and trained across several
viscosities, separates the Reynolds-dependent dissipation calibration from
the Reynolds-independent tensor structure; we test this on held-out
viscosities and filter ratios and on a decaying Taylor--Green vortex.
The remainder of this article is organized as follows.
\Cref{sec:symmetries} reviews the symmetries of the
incompressible Navier--Stokes equations and their implications for LES.
\Cref{sec:models} introduces the closure models: the classical baselines and
the three data-driven closures---an unconstrained MLP baseline, the
group-convolutional G-CNN, and the tensor-basis TBNN---together with the
filter-scale Reynolds-number input feature.
\Cref{sec:results} compares the closures on forced isotropic turbulence,
studies how their accuracy depends on network capacity, and tests their
generalization across Reynolds number, on held-out viscosities and filter
ratios and on a decaying Taylor--Green vortex.
\Cref{sec:conclusion} summarizes our findings and provides directions for future research.

\section{Symmetries in large-eddy simulation} \label{sec:symmetries}

The incompressible Navier--Stokes equations admit a set of symmetries, but not
all of them survive discretization and closure modeling. This section provides
the conceptual basis for the rest of the paper. We first motivate why a closure
model should respect these symmetries (\cref{sec:why-symmetries}), then review
the formal notions of transform, invariance, and equivariance
(\cref{sec:transforms}) and the symmetry group of the continuous equations
(\cref{sec:symmetries-navier-stokes}). We then filter the equations
(\cref{sec:symmetries-filtered}) and ask which of these symmetries are inherited
by the \emph{discrete} LES equations that our closures actually live in
(\cref{sec:symmetries-discrete}). This leads to a question that frames the
comparison carried out in the rest of the paper: whether a closure should
enforce a symmetry that the discrete equations do not actually possess
(\cref{sec:symmetries-enforce}).

\subsection{Why symmetries?} \label{sec:why-symmetries}

Symmetries are fundamental properties of both turbulence and the differential
equations that govern it~\cite{frischTurbulenceLegacyKolmogorov1995,
oberlackInvariantModelingLargeeddy1997}. The guiding principle of this work is
that a closure model should not break a symmetry possessed by the equations it
closes. A closure that violates rotational symmetry, for instance, predicts a
sub-filter stress that depends on the orientation of the coordinate axes, so
that merely rotating the grid changes the modeled dynamics; a closure that
violates Galilean invariance makes the stress depend on the velocity of the
observer's frame. The exact sub-filter stress does neither, so such a model
introduces spurious, frame-dependent forces that have no physical
counterpart~\cite{oberlackInvariantModelingLargeeddy1997,
razafindralandyAnalysisDevelopmentSubgrid2007,
silvisPhysicalConsistencySubgridscale2017}. For homogeneous isotropic
turbulence the argument is also statistical: the exact sub-filter stress is, by
Kolmogorov's hypothesis of local isotropy, invariant in distribution under
rotations and reflections, and a model that respects this invariance cannot
manufacture a preferred direction. Beyond physical consistency, enforcing a
symmetry is a parameter-free inductive bias: it restricts the hypothesis space
to functions compatible with the known invariances, which improves data
efficiency and has been reported to improve numerical
stability~\cite{guanLearningPhysicsconstrainedSubgridscale2023}. The remainder of this
section makes these notions precise and determines exactly which symmetries are
available to enforce in our discrete setting.

\subsection{Transforms and symmetries} \label{sec:transforms}

We first make precise what it means for the equations to possess a symmetry, and
then what it means for a closure to respect one.

Let $\mathcal{N}(u) = 0$ denote an equation for a field $u(x, t)$, and let
$\tilde{x}(x, t)$, $\tilde{t}(x, t)$, and $\tilde{u}(\tilde{x}, \tilde{t})$ denote
transformed coordinates and fields, with $\tilde{\mathcal{N}}(\tilde{u})$ the
correspondingly transformed equation operator. A transform $\widetilde{(\cdot)}$ is a
\emph{symmetry} of the equation $\mathcal{N} = 0$ if
\begin{equation}
    \mathcal{N}(u) = 0 \implies \tilde{\mathcal{N}}(\tilde{u}) = 0,
\end{equation}
i.e.\ solutions are mapped to solutions under the transform. This is the general
notion: the symmetries of the incompressible Navier--Stokes equations collected in
\cref{sec:symmetries-navier-stokes}---translations, Galilean boosts, scalings, and
rotations/reflections---are all transforms of this kind.

The spatial symmetries at the center of this paper, rotations and reflections, act
\emph{linearly} on physical space $\mathbb{R}^3$. They form a subgroup
$G \subseteq O(3)$, and each element $g \in G$ is represented by an orthogonal matrix
$R_g \in \mathbb{R}^{3 \times 3}$, so that $R_g^{-1} = R_g^T$. For such a transform,
the induced action on a scalar-, vector-, or tensor-valued field $p$, $u$, or
$\sigma$ depends only on its tensor character:
\begin{equation} \label{eq:group-action}
\begin{split}
    g p(x) & = p(R_g^{-1} x), \\
    g u(x) & = R_g u(R_g^{-1} x), \\
    g \sigma(x) & = R_g \sigma(R_g^{-1} x) R_g^{-1}.
\end{split}
\end{equation}
The remaining Navier--Stokes symmetries (translations, Galilean boosts, and
scalings) do not act through a single orthogonal $R_g$ on the field values; they are
captured by the general transform notion above rather than by \cref{eq:group-action}.

Finally, we make precise what it means for a closure to respect such a symmetry. An
operator $\mathcal{M}$, such as a closure model, is \emph{invariant} under $G$ if, for
all $g \in G$,
\begin{equation}
    \mathcal{M}(g \varphi) = \mathcal{M}(\varphi),
\end{equation}
and \emph{equivariant} if, for all $g \in G$,
\begin{equation}
    \mathcal{M}(g \varphi) = g \mathcal{M}(\varphi),
\end{equation}
where $\varphi$ is a scalar-, vector-, or tensor-valued field and the action $g$ is
that of \cref{eq:group-action}. Equivariance is the property we will ask of a closure
under rotations and reflections: rotating the input velocity rotates the predicted
stress in the same way.

\subsection{Symmetries of the incompressible Navier--Stokes equations}
\label{sec:symmetries-navier-stokes}

The incompressible Navier--Stokes equations on a domain $\Omega$ are given by
\begin{equation} \label{eq:navier-stokes}
    \partial_j u_j = 0, \quad
    \partial_t u_i
    + \partial_j \left( \sigma_{i j}(u) + p \delta_{i j} \right)
    = f_i,
\end{equation}
where
$u(x, t) \in \mathbb{R}^3$ is the velocity field,
$p(x, t) \in \mathbb{R}$ is the pressure field,
$(i, j) \in \{1, 2, 3\}$ are spatial indices,
$\delta_{i j}$ is the Kronecker symbol,
$f(x, t) \in \mathbb{R}^3$ is a body force,
$\sigma(u)$ contains the nonlinear and viscous stresses defined as
\begin{equation}
    \sigma_{i j}(u)
    \coloneq u_i u_j
    - \nu \left( \partial_{j} u_i + \partial_{i} u_j \right),
\end{equation}
and $\nu > 0$ is the viscosity.
Unless otherwise stated, the Einstein summation convention applies for repeated
indices.

We also define the velocity gradient tensor (VGT)
$A_{i j} \coloneq \partial_j u_i$ and the strain-rate and rotation-rate tensors
$S_{i j} \coloneq (A_{i j} + A_{j i}) / 2$ and $W_{i j} \coloneq (A_{i j} - A_{j i}) / 2$.

The complete symmetries of the
incompressible Navier--Stokes equations~\eqref{eq:navier-stokes} are listed below
\cite{oberlackInvariantModelingLargeeddy1997}.
\begin{itemize}
    \item Time invariance:
        $\tilde{t} = t + a$,
        $\tilde{x} = x$,
        $\tilde{u} = u$,
        $\tilde{p} = p$, where $a \in \mathbb{R}$ is a constant.
    \item Rotation invariance:
        $\tilde{x}_i = R_{i j} x_j$,
        $\tilde{t} = t$,
        $\tilde{u}_i = R_{i j} u_j$,
        $\tilde{p} = p$,
        where $R \in \mathbb{R}^{3 \times 3}$ is an orthogonal matrix.
    \item Reflection invariance in the direction $x_i$
        (with $j \neq i$):
        $\tilde{x}_i = -x_i$,
        $\tilde{x}_j = x_j$,
        $\tilde{t} = t$,
        $\tilde{u}_i = -u_i$,
        $\tilde{u}_j = u_j$,
        $\tilde{p} = p$.
    \item Generalized Galilean invariance:
        $\tilde{x}_i = x_i + X_i(t)$,
        $\tilde{t} = t$,
        $\tilde{u}_i = u_i + \dot{X}_i(t)$,
        $\tilde{p} = p - x_j \ddot{X}_j(t)$,
        where $X : \mathbb{R} \to \mathbb{R}^3$
        is a time-dependent twice differentiable frame translation.
    \item Scaling invariance:
        $\tilde{x}_i = a x_i$,
        $\tilde{t} = b t$,
        $\tilde{u}_i = \frac{a}{b} u_i$,
        $\tilde{p} = \left(\frac{a}{b}\right)^2 p$,
        $\tilde{\nu} = \frac{a^2}{b} \nu$,
        where $a > 0$ and $b > 0$ are constants (negative $a$ composes a
        scaling with the point reflection already covered above).
        Note that $\nu$ must be scaled to preserve the symmetry
        (unless $a^2 = b$).
    \item Pressure invariance:
        $\tilde{x} = x$,
        $\tilde{t} = t$,
        $\tilde{u} = u$,
        $\tilde{p} = p + a(t)$,
        where $a : \mathbb{R} \to \mathbb{R}$
        is constant in space.
\end{itemize}

We note that the symmetries above hold for the unforced equations.
A forcing term $f$ that is fixed in a specific frame
(such as the spectrally banded forcing used in our experiments,
see \cref{sec:datagen})
breaks the time-dependent Galilean transform with non-zero $\ddot{X}$,
so that only the constant-translation subgroup
$X(t) = X_0 + V t$ remains a symmetry of the forced system.
A closure model that depends only on $\bar{A}$ is, however, invariant under
the full generalized Galilean group, regardless of the forcing.

\subsection{Symmetries of the filtered equations} \label{sec:symmetries-filtered}

Consider a convolutional filter $\overline{(\cdot)}$ defined by
\begin{equation}
    \bar{u}(x) \coloneq \int_\Omega H(x - y) u(y) \, \mathrm{d} y
\end{equation}
for some kernel $H$. The filter commutes with differentiation:
$\partial_j \bar{u} = \overline{\partial_j u}$.
The filtered Navier--Stokes equations are given by
\begin{equation}
    \partial_j \bar{u}_j = 0, \quad
    \partial_t \bar{u}_i +
    \partial_j \left( \sigma_{i j}(\bar{u}) + \tau_{i j}(u) + \bar{p} \delta_{i j} \right)
    = \bar{f}_i,
\end{equation}
where
\begin{equation}
    \tau_{i j}(u) \coloneq \overline{u_i u_j} - \bar{u}_i \bar{u}_j
\end{equation}
is the sub-filter stress tensor (SFS).

Introducing a closure model $m(\bar{u})$ for $\tau(u)$ yields the LES equations
\begin{equation}
    \partial_j v_j = 0, \quad
    \partial_t v_i +
    \partial_j \left( \sigma_{i j}(v) + m_{i j}(v) + p_v \delta_{i j} \right)
    = \bar{f}_i,
\end{equation}
where $p_v$ is the pressure that enforces the divergence-free constraint on the LES solution $v$,
which in general differs from the filtered pressure $\bar{p}$.
Functional closures are designed to match the dissipation rate of the true SFS, such that
$\bar{S}_{i j} m_{i j}(\bar{u}) \approx \bar{S}_{i j} \tau_{i j}(u)$.
They are often expressed in eddy-viscosity
form~\cite{smagorinskyGeneralCirculationExperiments1963,
lillyRepresentationSmallscaleTurbulence1966}:
\begin{equation} \label{eq:eddy-viscosity}
    m^\text{EV}_{i j}(\bar{u}) \coloneq - \nu^\Delta \bar{S}_{i j}
\end{equation}
for some eddy viscosity $\nu^\Delta$.
Structural closures are designed to approximate the SFS tensor itself:
$m_{i j}(\bar{u}) \approx \tau_{i j}(u)$. An example is Clark's gradient model
\begin{equation} \label{eq:clark}
    m_{i j}^\text{Clark}(\bar{u}) = \frac{\Delta^2}{12} \bar{A}_{i k} \bar{A}_{j k}
\end{equation}
which is obtained by truncating a Taylor-series expansion of $\tau(u)$.

Since replacing $m_{i j}$ with $m_{i j} + n \delta_{i j}$ for any scalar function $n$
does not change the dynamics---the isotropic part is
absorbed by the pressure $p_v$---it is convenient to choose $m$ to be deviatoric (trace-free) and to
fit it to the deviatoric part
$\tau^\text{dev}_{i j} \coloneq \tau_{i j} - \tau_{k k} \delta_{i j} / 3$.

Depending on the choice of closure $m$, the symmetry properties of the original
equations may be violated. To enforce Galilean invariance, it is standard practice
to express $m$ as a function
of $\bar{A}$ rather than $\bar{u}$~\cite{lundParameterizationSubgridscaleStress1993},
since $\bar{u} + a(t)$ and $\bar{u}$ have identical velocity gradients for any scalar $a(t)$.
Eddy viscosity models and Clark's gradient model depend solely on $\bar{A}$, and thus
preserve Galilean invariance.
More precisely, since $\bar{A}$ is invariant under
the generalized Galilean transform of \cref{sec:symmetries-navier-stokes},
any closure depending solely on $\bar{A}$ preserves
generalized Galilean invariance as well.
All closures considered in this paper are of this type, and we henceforth
write them as functions $m(\bar{A})$; where a closure is evaluated on a
resolved velocity field, as in the LES equations above or the error metrics
of \cref{sec:results}, $m(\bar{u})$ abbreviates the composition of $m$ with
the gradient $\bar{u} \mapsto \bar{A}$.

\subsection{Symmetries inherited by the discrete LES equations}
\label{sec:symmetries-discrete}

The discussion so far is continuous. The closures in this paper, however, are
embedded in a discrete LES solver: a pseudo-spectral discretization on a uniform
Cartesian grid (\ref{sec:pseudo-spectral}), closed by a model $m(\bar{A})$ that
is applied independently at each grid point. Discretization and grid sampling
break some of the continuous symmetries of \cref{sec:symmetries-navier-stokes},
and a closure can only meaningfully preserve those that the discrete equations
actually retain. We go through them in turn; \cref{tab:symmetry-inheritance}
summarizes the outcome.

\paragraph{Rotation and reflection}
This is the most consequential case. A uniform Cartesian grid is not invariant
under arbitrary rotations: only rotations by multiples of $\pi / 2$ and
reflections aligned with the grid axes map the grid onto itself. These transforms
form the \emph{octahedral group} $G \subset O(3)$, with $|G| = 48$ elements
(\ref{sec:group-parametrization}). The continuous orthogonal group $O(3)$ is
therefore broken down to $G$ the moment the velocity field is sampled on the
grid. The breaking originates in the grid and the dealiasing, not in the filter:
the $2/3$-rule truncation retains modes with $| k |_\infty \le N / 3$, a
\emph{cube} in wavenumber space that is invariant only under $G$, whereas the
Gaussian data-filter of \cref{eq:filter} depends on $| k |$ and is fully
isotropic. All the discrete operators commute with the signed axis permutations that
represent $G$: the discrete Fourier transform, the cubic truncation, the
Leray projector, the nonlinear product, the shell-energy rescaling forcing
(which depends only on the rotation-invariant shell energies;
\ref{sec:datagen}), and even the adaptive time step, which depends only on
the grid-invariant maximum velocity magnitude and grid spacing.
Hence $G$ is an \emph{exact} symmetry of the discrete equations,
and it is the largest rotation/reflection symmetry that any grid-sampled field
can realize exactly.

\paragraph{Scaling}
Scaling symmetry survives only in a weaker, family-wise sense. Strictly, no
nontrivial scaling maps the discrete system to itself: the transform
$\tilde{x} = a x$ carries the fixed box, grid, filter width, and forcing
shells into rescaled copies, so scaling relates \emph{different} members of
the family of discrete systems parametrized by the domain size, $\Delta$, and
$\nu$, rather than fixing any single member. What a closure can meaningfully
preserve is therefore \emph{covariance} along this family. The output
prefactor $\Delta^2 |\bar{A}|^2$ shared by all our data-driven models provides
exactly that: the predicted stress transforms as a stress, the
scale-invariance long advocated as a design principle for subgrid
models~\cite{meneveauScaleInvarianceTurbulenceModels2000}. A Reynolds-blind
closure built from $\Delta$ and $\bar{A}$ alone, however, never sees the
viscosity $\nu$, so its prediction can be consistent with the rescaled
dynamics only along the $a^2 = b$ subgroup, which leaves $\nu$ (and hence the
filter-scale Reynolds number) unchanged; we defer the unit analysis to
\cref{sec:construction}. The $+\mathrm{Re}$ variants of
\cref{sec:redelta-input} restore precisely the missing viscosity dependence,
and with it covariance under the full two-parameter group.

\paragraph{Galilean invariance}
Constant Galilean shifts are inherited, because the closure depends on $\bar{A}$
rather than $\bar{u}$ (\cref{sec:symmetries-filtered}). The frame-fixed forcing
breaks the time-dependent part of the generalized Galilean group, as already
noted in \cref{sec:symmetries-navier-stokes}.

\paragraph{Translation, time, and pressure}
The remaining symmetries survive. Discrete grid translations are exact, and
because the closure is applied identically and pointwise at every grid point it
is automatically equivariant under them; the pseudo-spectral representation
extends this to continuous translation up to the resolved band. Time-translation
invariance holds for the autonomous scheme with its deterministic per-step
forcing, and the pressure gauge is preserved by taking the closure deviatoric
(\cref{sec:symmetries-filtered}).

\begin{table*}
    \centering
    \caption{%
        Fate of the continuous Navier--Stokes symmetries
        (\cref{sec:symmetries-navier-stokes}) in the discrete LES equations closed
        by a pointwise model $m(\bar{A})$. All of them hold for the continuous
        equations; the grid, the dealiasing, and the frame-fixed forcing break or
        restrict several of them. The $+\mathrm{Re}$ closures of
        \cref{sec:redelta-input} additionally carry the viscosity dependence
        needed for covariance under the full two-parameter scaling group.
    }
    \label{tab:symmetry-inheritance}
    \begin{tabular}{lll}
        \toprule
        Symmetry & Continuous Navier--Stokes & Pseudo-spectral forced LES $+$ our closure \\
        \midrule
        Time translation              & preserved & preserved \\
        Space translation             & preserved & preserved (closure is pointwise) \\
        Constant Galilean             & preserved & preserved (closure uses $\bar{A}$) \\
        Generalized Galilean          & preserved & broken by the frame-fixed forcing \\
        Rotation/reflection $O(3)$    & preserved & reduced to the octahedral group $G$ ($|G| = 48$) \\
        Scaling (two-parameter)       & preserved & reduced to covariance along $a^2 = b$ (via $\Delta^2 |\bar{A}|^2$) \\
        Pressure gauge                & preserved & preserved (deviatoric closure) \\
        \bottomrule
    \end{tabular}
\end{table*}

\subsection{Should a closure enforce a symmetry the discrete equations lack?}
\label{sec:symmetries-enforce}

\Cref{tab:symmetry-inheritance} raises the first of the questions this paper
addresses: should a closure enforce a symmetry that the discrete equations it
is embedded in do not actually possess? It is useful to distinguish two carriers of
symmetry: the discrete \emph{dynamical operator}, and the \emph{constitutive
relation} $m(\bar{A})$ that the closure represents.

For a symmetry that the discrete operator possesses \emph{exactly} (the
octahedral group $G$ and the constant Galilean shift), enforcing it in the
closure is unambiguously correct and comes at no cost. The exact discrete SFS
$\tau^{N \to M}$ of \cref{eq:discrete-sfs} is itself exactly $G$-equivariant,
so constraining the model to be $G$-equivariant removes only degrees of
freedom that the target also lacks; likewise, the target transforms as a
stress across rescaled systems, so the $a^2 = b$ scaling covariance built
into the prefactor discards nothing. This is the regime of the
group-convolutional network, which enforces $G$ exactly, and of the
Galilean/scaling prefactor shared by all three of our models.

Enforcing continuous $O(3)$ equivariance is a different matter. The discrete
equations are not $O(3)$-invariant, and the discrete target $\tau^{N \to M}$,
assembled from cubic truncations, is only octahedral-equivariant; imposing the
full orthogonal group therefore enforces a symmetry the discrete system does not
have. It does not follow that doing so is wrong. A closure is a model of a
\emph{continuum} quantity: the sub-filter stress is defined by the continuous
filtered equations (\cref{sec:symmetries-filtered}), and for homogeneous
isotropic turbulence its statistics are $O(3)$-invariant by Kolmogorov's
hypothesis of local isotropy~\cite{frischTurbulenceLegacyKolmogorov1995}. Building
$O(3)$ equivariance into the model is then a physically motivated inductive bias
on the continuum object being approximated, and the mismatch with the discrete
target is of the order of the discretization anisotropy---small when the filter
is well resolved on the grid. The distinction is in fact partly moot in practice:
a model built to be $O(3)$-equivariant in the continuum, such as the tensor-basis
network of \cref{sec:tbnn}, collapses to exact octahedral equivariance once its
input gradient is sampled on the grid, recovering the same symmetry that the
group-convolutional network enforces directly.

The one genuine caveat is that enforcing an absent symmetry can
\emph{over-constrain} the closure when the discrete or physical target truly
violates it. In anisotropic or wall-bounded flows, or with coarse filters whose
cutoff reaches into the anisotropic energy-containing range, the sub-filter stress
need not be $O(3)$-equivariant even in the continuum, and a model that insists on
full rotational equivariance may be too rigid. In the isotropic, well-resolved
setting studied here the constraint is benign, and our experiments bear this out:
at sufficient capacity the equivariant and unconstrained networks attain
practically equivalent accuracy, with remaining differences at the level of
a percent of relative error, resolvable across training seeds but immaterial
a posteriori (\cref{sec:results}). We return to the anisotropic case, and to relaxed-equivariance
constructions~\cite{wangDiscoveringSymmetryBreaking2025} that interpolate between the two
regimes, in \cref{sec:conclusion}.

\section{Closure models} \label{sec:models}

This section defines the closure models compared in this work: three classical
baselines (\cref{sec:classical}) and three data-driven closures that share a
pointwise velocity-gradient construction (\cref{sec:construction}) but differ
in how they treat rotational and reflectional symmetry
(\cref{sec:mlp,sec:groupconv,sec:tbnn}), each data-driven closure also coming
in a $+\mathrm{Re}$ variant that receives the filter-scale Reynolds number as
an extra input (\cref{sec:redelta-input}).

\subsection{Classical baselines} \label{sec:classical}

We compare the data-driven closures against three classical models.
The first is the trivial closure $m^\text{no-model}(\bar{A}) \coloneq 0$
(No-model), i.e.\ a coarse simulation without any closure.
The second is the dynamic Smagorinsky model, a functional closure of
eddy-viscosity form \eqref{eq:eddy-viscosity},
\begin{equation}
    m^\text{dyn.~smag}(\bar{A}) \coloneq -2 c \Delta^2 |\bar{S}| \bar{S},
\end{equation}
where $|\bar{S}| = \sqrt{\bar{S}_{i j} \bar{S}_{i j}}$ and
the coefficient $c$ is a single, volume-averaged scalar per
time step, determined from the resolved field via the Germano--Lilly
dynamic procedure~\cite{germanoDynamicSubgridscaleEddy1991,
lillyProposedModificationGermano1992}; details are given in
\ref{sec:dynamic-smagorinsky}.
The third is Clark's gradient model \eqref{eq:clark}, the structural closure
obtained by truncating the Taylor expansion of the filter kernel
(\cref{sec:symmetries-filtered}).
The two nontrivial baselines illustrate the contrasting closure philosophies of
\cref{sec:symmetries-filtered}: dynamic Smagorinsky is purely dissipative and
does not account for backscatter of energy from small to large scales, whereas
Clark captures the tensor structure of the SFS but is insufficiently
dissipative when used on its
own~\cite{vremanLargeeddySimulationTemporal1996}.

\subsection{A shared construction for the data-driven closures}
\label{sec:construction}

The three data-driven closures all predict the sub-filter stress directly from
the local filtered
velocity-gradient tensor (VGT) $\bar{A}$, and all three share the same input
normalization $\bar{A} / |\bar{A}|$ and output prefactor $\Delta^2 |\bar{A}|^2$.
Here $|\cdot|$ denotes the Frobenius norm, $|T| \coloneq \sqrt{T_{i j} T_{i j}}$ for a
tensor $T$.
Collected into one form, each data-driven closure reads
\begin{equation} \label{eq:shared-form}
    m(\bar{A}) \coloneq \Delta^2 |\bar{A}|^2 \, m^*(\bar{A}^*),
    \qquad
    \bar{A}^* \coloneq \bar{A} / |\bar{A}|,
\end{equation}
where the dimensionless \emph{closure core} $m^*$ maps the normalized VGT to
a normalized stress. The core is the only part in which the three closures
differ (\cref{sec:mlp,sec:groupconv,sec:tbnn}).
They are thus all \emph{one-point} closures: the predicted stress at a grid
point depends only on $\bar{A}$ at that same point.
This restriction is a premise of the paper rather than a mere simplification.
It fixes a common hypothesis space---functions of the local $\bar{A}$---so
that the three architectures differ only in how they parametrize this space,
and the capacity and generalization studies of \cref{sec:results} compare like
with like; and it makes the \emph{best attainable} closure in this class
directly computable as a conditional mean, giving those comparisons an
absolute reference point (\cref{sec:results-capacity}).
The classical baselines below are one-point closures as well (the
dynamic-Smagorinsky coefficient is determined globally, but the resulting
stress is still a local function of $\bar{A}$).
This shared construction enforces two symmetries and dimensional consistency
at once.
Using the VGT rather than the velocity enforces Galilean invariance, since
$\bar{u}$ and $\bar{u} + a(t)$ have identical gradients
(\cref{sec:symmetries-filtered}).
The prefactor provides the correct units of $\text{length}^2 / \text{time}^2$:
the network maps the dimensionless input $\bar{A} / |\bar{A}|$ to a
dimensionless output, and $\Delta^2 |\bar{A}|^2$ is the unique combination of
$\Delta$ and $|\bar{A}|$ with the units of $\tau$.
Nondimensionalizing a learned closure's inputs and outputs by such
filter-scale factors is an established device for improving generalization to
unseen grid resolutions~\cite{perezhoginGeneralizableNeuralNetworkParameterization2025}.
This matches the scaling of Clark's gradient model, which is second order in
the VGT; incorporating higher-order contributions would require the
viscosity $\nu$.
The same construction makes the closures \emph{covariant} under the scaling
transform $\tilde{x} = a x$, $\tilde{t} = b t$, $\tilde{u} = (a / b) u$ from
\cref{sec:symmetries-navier-stokes}. Under this transform the velocity gradient
scales as $\bar{A} \mapsto \bar{A} / b$ (it has units of inverse time), so the
normalized input $\bar{A} / |\bar{A}|$ is invariant, while the prefactor scales
as $\Delta^2 |\bar{A}|^2 \mapsto (a / b)^2 \, \Delta^2 |\bar{A}|^2$---exactly as
a stress does ($\tilde{\tau} = (a / b)^2 \tau$).
The predicted stress therefore transforms identically to the true sub-filter
stress for the full two-parameter group: this is dimensional covariance, not
invariance.
What singles out the $a^2 = b$ subgroup is the viscosity. Along it the viscosity
is unchanged ($\tilde{\nu} = (a^2 / b) \nu = \nu$), so the filter-scale Reynolds
number is preserved and a closure built from $\Delta$ and $\bar{A}$ alone is
fully consistent with the dynamics. Off this subgroup the same transform rescales
the viscosity, $\tilde{\nu} = (a^2 / b) \nu \neq \nu$, and hence the filter-scale
Reynolds number; a closure with no dependence on $\nu$ stays dimensionally
covariant but cannot reproduce the accompanying change in the sub-filter stress.
The closures are thus dimensionally consistent at any single Reynolds number but
carry no explicit $\operatorname{Re}_\Delta$ dependence, so they are not
guaranteed to generalize across Reynolds numbers; \cref{sec:redelta-input}
introduces the input feature that restores this dependence.

The three models differ in how, if at all, they additionally enforce
equivariance to rotations and reflections.
The unconstrained multi-layer perceptron (MLP) of \cref{sec:mlp} is a baseline
that enforces no such equivariance.
The group-convolutional network (G-CNN) of \cref{sec:groupconv} enforces it
through the network architecture, by constraining the weights, whereas the
tensor-basis neural network (TBNN) of \cref{sec:tbnn} enforces it through the
input and output features, by expanding the stress in an equivariant tensor
basis.

\subsection{Feed-forward neural network (MLP)} \label{sec:mlp}

We first define a standard feed-forward neural network closure that is not
equivariant to rotation or reflection, serving as an unconstrained baseline.
A feed-forward layer $\ell: \xi \mapsto \zeta$, mapping an input channel
vector $\xi$ to an output channel vector $\zeta$, is defined by
\begin{equation}
    \zeta_i = \varphi\Bigl(\textstyle\sum_j w_{i j} \xi_j + b_i\Bigr),
\end{equation}
where $i$ and $j$ index the output and input channels,
$w_{i j} \in \mathbb{R}$ are weights, $b_i \in \mathbb{R}$ are biases,
and $\varphi$ is a non-linear activation function.
The input and output shapes are given by the size of the weight matrix $w$.

We define a neural structural closure model as the closure
\eqref{eq:shared-form} whose core is the multi-layer perceptron (MLP)
\begin{equation} \label{eq:m-conv}
    m^{*}_\text{MLP} \coloneq
    \ell_n \circ \ell_{n - 1} \circ \dots \circ \ell_1,
\end{equation}
where $\ell_1, \dots, \ell_n$ are different layers with independent parameters.
This is a structural closure: the core directly predicts
the six independent components of the normalized symmetric stress tensor,
from which the trace is removed to obtain the deviatoric closure.
Following the one-point premise of \cref{sec:construction}, the closure is
applied at each grid point independently: each layer is implemented as a
convolutional layer with a point kernel of physical
size $1 \times 1 \times 1$.
This architecture is not equivariant to rotation or reflection.
To separate the \emph{learned function} from the \emph{architectural
constraint}, we additionally evaluate a \emph{symmetrized} MLP,
\begin{equation} \label{eq:mlp-sym}
    m^\text{MLP-sym}(\bar{A})
    \coloneq \frac{1}{|G|} \sum_{g \in G} R_g^T
    m^\text{MLP}(R_g \bar{A} R_g^T) R_g,
\end{equation}
the trained MLP with its prediction averaged over the $48$ octahedral group
elements at inference time.
This wrapper is exactly equivariant while reusing the MLP's parameters
unchanged, at $48$ forward passes per closure evaluation: if the symmetrized
MLP matches the raw MLP, the function the MLP learned is already
near-equivariant.
The next two subsections enforce these additional symmetries in two different
ways: the G-CNN (\cref{sec:groupconv}) constrains the network weights, and the
TBNN (\cref{sec:tbnn}) constrains the input and output features.

\subsection{Group-convolutional network (G-CNN)} \label{sec:groupconv}

Since Galilean and scaling invariance are already enforced through the
input normalization and output prefactor,
we focus on equivariance to rotations and reflections.
The continuous symmetry group is the \emph{orthogonal group} $O(3)$, but
infinitesimal rotations cannot be represented exactly on a Cartesian grid:
only rotations by multiples of $\pi / 2$ and reflections aligned with the grid
axes map the grid onto itself.
These form the \emph{octahedral group} $G$, with $|G| = 48$ elements.
The G-CNN is $G$-equivariant by construction; once the velocity gradient is
sampled on the grid, the octahedral group is the largest symmetry any
grid-based closure can realize exactly.
Our parametrization of $G$ is described in \ref{sec:group-parametrization}.
Like the MLP, the G-CNN acts pointwise in physical space (a $1 \times 1 \times
1$ kernel); the ``convolution'' in group-convolution is over the group $G$, not
over space.

To make a layer $\ell$ group-equivariant, we require that
$g \ell(\xi) = \ell(g \xi)$ for all group elements $g \in G$.
A simple way to achieve this is to use group-convolutional
layers~\cite{cohenGroupEquivariantConvolutional2016}.

In physical space $\mathbb{R}^3$,
a group element $g \in G$ is represented by a
roto-reflection matrix $R_g \in \mathbb{R}^{3 \times 3}$,
the orthogonal physical-space representation introduced in \cref{sec:transforms}.
These matrices can have negative entries, causing the
signs of transformed physical quantities to change depending on the group element.
For sign-sensitive activation functions, such as
$\operatorname{ReLU}(x) \coloneq \max(0, x)$ or the smooth
$\operatorname{GELU}$ used in our networks, this representation is
problematic.
Although Shang et al.\ proposed concatenated activation functions compatible with
signed permutation matrices~\cite{shangUnderstandingImprovingConvolutional2016},
biases in layers operating on physical $3 \times 3$ matrices break rotational
equivariance, since biases do not transform under group actions.
We therefore employ the \emph{regular representation} space
$\mathbb{R}^{|G|} = \mathbb{R}^{48}$ for the hidden layers of the neural network.
In this representation, group elements are
represented by permutation matrices $P_g \in \mathbb{R}^{48 \times 48}$
with entries that are either $0$ or
$1$~\cite{cohenGroupEquivariantConvolutional2016,cohenSteerableCNNs2016}.
Intuitively, each of the $48$ channels corresponds to one group element,
and the group action simply permutes these channels.

We index the $48$ entries of a regular-representation vector by the group
elements themselves: $\xi_j(g) \in \mathbb{R}$ denotes the component of input
channel $j$ associated with $g \in G$.
In this indexing the permutation matrix $P_g$ acts by relabeling the
components according to the group action,
\begin{equation} \label{eq:regular-action}
    (P_g \xi)(h) = \xi(g^{-1} h), \qquad P_g(h, k) = \delta_{h, \, g k},
\end{equation}
where $\delta_{h, \, g k}$ equals $1$ if $h = g k$ and $0$ otherwise.
A layer maps $n_\text{in}$ regular-representation vectors $\xi_j \in \mathbb{R}^{48}$
to $n_\text{out}$ regular-representation vectors $\zeta_i \in \mathbb{R}^{48}$ by a
\emph{group convolution}: the spatial shifts of an ordinary convolution are
replaced by group multiplication, and a single kernel is shared across the group,
\begin{equation} \label{eq:gconv-layer}
    \zeta_i(g) = \varphi\Biggl(
        \sum_{j = 1}^{n_\text{in}} \sum_{h \in G} k_{i j}(h^{-1} g) \, \xi_j(h) + b_i
    \Biggr),
\end{equation}
where $k_{i j} : G \to \mathbb{R}$ is the learnable kernel for the channel pair
$(i, j)$, i.e.\ $|G| = 48$ numbers.
The inner sum over $h \in G$ contracts the group dimension, while the outer
sum over $j$ runs over the input channels.
The bias $b_i \in \mathbb{R}$ is a single scalar added to every group component
$g$; being constant across the $48$ components, it is left invariant by the
permutation $P_g$.

This layer is equivariant by construction: a change of variables in the group
sum shows that its linear part commutes with every permutation $P_a$, while the
bias (constant across the group components) and the element-wise activation
$\varphi$ commute with permutations trivially~\cite{cohenSteerableCNNs2016}, so
$P_a \, \ell(\xi) = \ell(P_a \xi)$ for all $a \in G$; the explicit computation is
given in \ref{sec:weight-sharing}.
Conversely, every $G$-equivariant linear map between regular representations is of
this group-convolutional form~\cite{cohenGroupEquivariantConvolutional2016,
cohenSteerableCNNs2016}, so the shared kernel $k_{i j}(h^{-1} g)$ captures the
entire equivariant subspace without loss of expressivity---the weight sharing
we make explicit below.

For the first neural network layer, we need a similar construction.
The input is the (filtered) velocity gradient tensor $\bar{A}$
(\cref{sec:symmetries}), which transforms as
\begin{equation}
    g \bar{A} = R_g \bar{A} R^T_g.
\end{equation}
If we ``flatten'' $\bar{A}$ into a vector $a \in \mathbb{R}^{9}$ with
components $a(\mu)$, $\mu \in \{1, \dots, 9\}$,
the left and right multiplications by $R_g$ combine into the single
tensor-representation matrix $Q_g \coloneq R_g \otimes R_g \in \mathbb{R}^{9 \times 9}$,
so that $a$ transforms as $(g a)(\mu) = \sum_\nu Q_g(\mu, \nu) a(\nu)$.
There is now a single tensor-valued input channel, so in place of the shared
kernel $k_{i j}$ we tie the weights through the tensor representation: the initial
layer $\ell : a \mapsto \zeta$ is the \emph{orbit} under $Q$ of a single $9$-vector
$c_i \in \mathbb{R}^9$ per output channel,
\begin{equation}
\begin{split}
    \zeta_i(g) & = \varphi\Biggl( \sum_{\mu = 1}^{9} w_{i}(g, \mu) \, a(\mu) + b_i \Biggr), \\
    w_{i}(g, \mu) & = \sum_{\nu = 1}^{9} Q_g(\mu, \nu) \, c_i(\nu),
\end{split}
\end{equation}
where the sum over $\mu$ contracts the $9$ flattened-tensor components.
Equivariance follows as before: using the orthogonality $Q_g^T = Q_{g^{-1}}$, the
orbit intertwines the input action $Q_g$ with the output permutation $P_g$, so
that $g \ell(a) = \ell(g a)$ for all $g \in G$ (\ref{sec:weight-sharing}).
Finally, in the last layer we predict the sub-filter stress tensor,
represented by its $9$ components as a flattened vector $m \in \mathbb{R}^9$, as
the orbit of a single $9$-vector $d_j \in \mathbb{R}^9$ per input channel,
\begin{equation}
\begin{split}
    m(\mu) & = \sum_{j = 1}^{n_\text{in}} \sum_{h \in G} w_{j}(\mu, h) \, \xi_j(h), \\
    w_{j}(\mu, h) & = \sum_{\nu = 1}^{9} Q_h(\mu, \nu) \, d_j(\nu).
\end{split}
\end{equation}
Each boundary block therefore carries only $9$ free parameters: an intertwiner
into or out of the regular representation is fixed by its value at the identity
(Schur's lemma / Frobenius reciprocity), so these $9$ numbers already span the
full intertwiner space.
We use no activation function and no bias in the final layer,
since $Q^{-1}_g$ does not commute with
elementwise activation functions $\varphi$
(unlike $P_g$, which is a permutation matrix with exactly one nonzero entry per row).
We reshape $m$ into a $3 \times 3$ tensor.
Since the SFS tensor is symmetric ($\tau = \tau^T$),
we discard the antisymmetric part by symmetrization $m \leftarrow (m + m^T) / 2$,
which retains six independent components.
We additionally remove the trace, $m \leftarrow \operatorname{dev}(m) = m - \operatorname{tr}(m) \, I / 3$,
so that the prediction is deviatoric;
the isotropic part is absorbed into the LES pressure as discussed in
\cref{sec:symmetries}.
Both operations commute with the action of $g \in G$ on $3 \times 3$ tensors,
$\sigma \mapsto R_g \sigma R_g^T$, so equivariance is preserved.

The initial, inner, and final layers are thus all equivariant by construction,
parametrized respectively by the boundary orbits $c_i$, the shared kernels
$k_{i j}$, and the boundary orbits $d_j$.
The group-convolutional closure model $m^\text{G-CNN}$ therefore has the same
architecture as the MLP core in \eqref{eq:m-conv}, except that every
feed-forward layer is replaced by the corresponding group-convolutional layer
$\ell^G : \xi \mapsto \zeta$ of \eqref{eq:gconv-layer} (with the tied boundary layers
above).
The final model is the closure \eqref{eq:shared-form} with core
\begin{equation}
    m^{*}_\text{G-CNN} \coloneq
    \ell^G_n \circ
    \ell^G_{n - 1} \circ
    \dots \circ
    \ell^G_1.
\end{equation}
\begin{figure*}
    \centering
    \def\svgwidth{0.9\textwidth}
    \includegraphics[width=\textwidth]{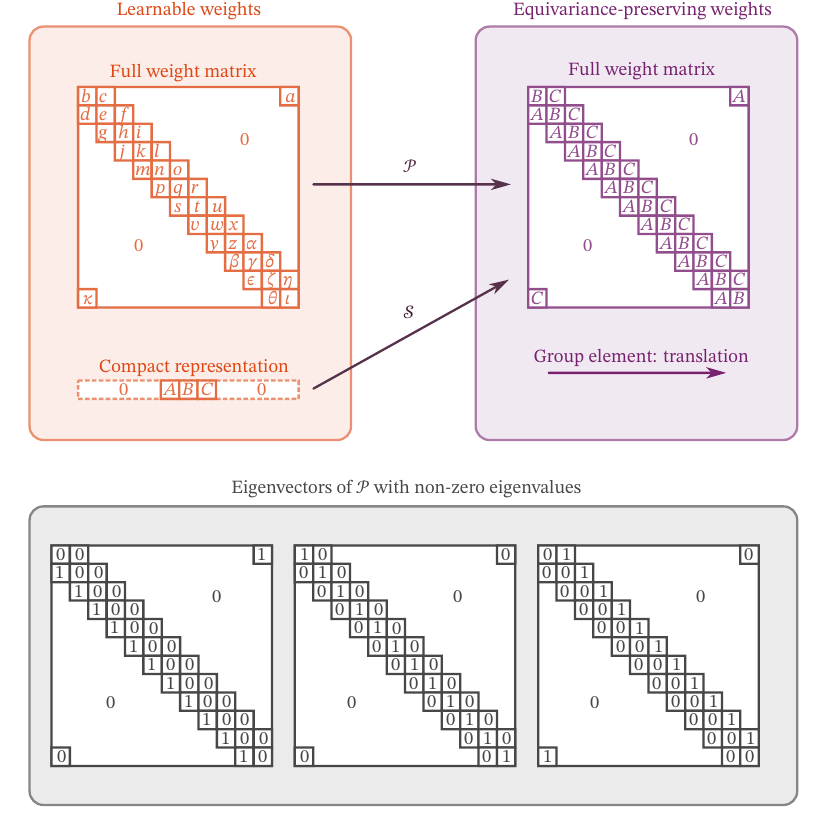}
    \caption{
        Two routes to equivariant weights, illustrated for
        translation-equivariant convolutions in 1D (classical CNNs): the
        projection $\mathcal{P}$ group-averages an arbitrary weight matrix,
        while the synthesis $\mathcal{S}$ builds the circulant matrix directly
        from the convolution kernel (the compact representation). Classical CNNs
        combine weight sharing (each row is a shift of the previous one) with
        locality (most weights forced to zero). Our octahedral
        group-convolutions use the same weight sharing but no
        locality (roto-reflections have no natural notion of nearness), so the
        synthesized weight blocks are dense.
    }
    \label{fig:group-conv}
\end{figure*}

The equivariant weight blocks can always be reached by \emph{projection}:
group-averaging an arbitrary block $\tilde{w}_{i j}$ over $G$,
$w_{i j} = |G|^{-1} \sum_{g \in G} P^{-1}_g \tilde{w}_{i j} P_g$, which makes
$w_{i j}$ commute with every $P_g$~\cite{cohenGroupEquivariantConvolutional2016}.
Rather than materialize this sum, or assemble and eigendecompose the dense
projector to extract its eigenvalue-one subspace when no closed form is
available, we \emph{synthesize} each block directly.
Commuting with the regular representation forces the group-circulant form
$w_{i j}(g, h) = k_{i j}(h^{-1} g)$ established above, so the synthesis
$\mathcal{S}$ is a pure \emph{gather}: it copies the learnable value
$k_{i j}(h^{-1} g)$ into entry $(g, h)$ of the block, and analogously copies the
boundary $9$-vectors $c_i, d_j$ along the orbits of $Q_g$.
Because $P_g$ and $Q_g = R_g \otimes R_g$ have integer entries in
$\{0, \pm 1\}$, $\mathcal{S}$ introduces no irrational numbers and the
synthesized blocks are \emph{exactly} equivariant in floating-point arithmetic.
This weight-sharing scheme is that of Cohen and
Welling~\cite{cohenGroupEquivariantConvolutional2016}. What is specific to our
setting is the explicit octahedral instance, including the
tensor-representation boundary layers that carry the stress-tensor input and
output; it keeps the construction closed-form and floating-point-exact, and
is illustrated for the 1D-translation group in \cref{fig:group-conv}.

During training, the synthesis $\mathcal{S}$ is part of the
computational graph: at every gradient-descent step it maps the reduced,
learnable coordinates to the full equivariant weight blocks $w$, through which
gradients are backpropagated.
The synthesis is precomputed once and reused.
After training, the resulting weights $w$ are precomputed and stored, so that at
inference time $m^\text{G-CNN}$ operates as a standard neural network with fixed
weights $w$.

\subsection{Tensor-basis neural network (TBNN)} \label{sec:tbnn}

The final data-driven closure enforces equivariance not through the network
architecture, as the G-CNN does, but through the choice of input and output
features. In contrast to the G-CNN, which is equivariant only to the octahedral
subgroup realizable on the grid, the tensor-basis construction is equivariant to
the full continuous orthogonal group $O(3)$ as a pointwise map; on grid-sampled
fields, however, only the octahedral subgroup of this equivariance can be
exercised, since only those rotations map one sampled field to another
(\cref{sec:symmetries-discrete}).

Consider a closure model that depends on the VGT $\bar{A}$.
For notational simplicity, we drop the overbar on $A$, $S$, $W$, etc.
Following tensor representation
theory~\cite{spencerTheoryMatrixPolynomials1958,spencerFurtherResultsTheory1959},
Pope~\cite{popeMoreGeneralEffectiveviscosity1975} proposed expressing such a model as
\begin{equation}
    m^\text{TB} \coloneq \Delta^2 \sum_{k} \alpha_k(\lambda) T_k,
\end{equation}
where $\Delta$ is the filter width,
$(T_k)_k$ is a tensor basis (TB)
of predetermined size,
and the coefficients $\alpha_k$ are
scalar functions of the invariants $\lambda$ of the VGT.
These invariants are
\begin{equation} \label{eq:invariants}
    \begin{alignedat}{2}
        \lambda_1 & \coloneq \operatorname{tr}(S^2), & \quad \lambda_4 & \coloneq \operatorname{tr}(S W^2), \\
        \lambda_2 & \coloneq \operatorname{tr}(W^2), & \quad \lambda_5 & \coloneq \operatorname{tr}(S^2 W^2). \\
        \lambda_3 & \coloneq \operatorname{tr}(S^3), & &
    \end{alignedat}
\end{equation}
Pope's original basis, consisting of $10$
tensors
(or $11$ in the compressible case)~\cite{popeMoreGeneralEffectiveviscosity1975},
was used by Ling et al.~\cite{lingReynoldsAveragedTurbulence2016}
in the context of RANS closures.

Based on the analysis of Smith~\cite{smithIsotropicFunctionsSymmetric1971},
Stallcup et al.~\cite{stallcupAdaptiveScaleSimilarClosure2022} employ a
slightly different basis that is both minimal and complete.
This basis comprises eight tensors $T_0, \dots, T_7$, including the isotropic
identity $T_0 = I$; in the incompressible case the isotropic part is carried by
the pressure, leaving the seven deviatoric tensors $T_1, \dots, T_7$ that are
actually used in the closure.
We adopt this basis, defined as
\begin{equation}
    \begin{alignedat}{2}
        T_0 & \coloneq I,   & \quad T_4 & \coloneq S W - W S, \\
        T_1 & \coloneq S,   & \quad T_5 & \coloneq W S W, \\
        T_2 & \coloneq S^2, & \quad T_6 & \coloneq S^2 W - W S^2, \\
        T_3 & \coloneq W^2, & \quad T_7 & \coloneq W S W^2 - W^2 S W.
    \end{alignedat}
\end{equation}
The fifth tensor $T_5 \coloneq W S W$ was not included in Pope's basis, which
instead included fifth-order terms in $S$ and $W$.

The key property of this formulation is that
the invariants $\lambda$ and basis tensors $T$ transform appropriately under
the rotation and reflection symmetry groups (i.e.\ the orthogonal group $O(3)$)
of the Navier--Stokes equations, so the
full model $m^\text{TB}$ is equivariant regardless of the functional form of the
coefficients $\alpha_k$.
For example, let $R \in \mathbb{R}^{3 \times 3}$ be a rotation matrix.
The invariants, being scalar fields, satisfy
$\lambda_1(A) =
\operatorname{tr}(S^2) =
\operatorname{tr}(R S^2 R^T) =
\operatorname{tr}((R S R^T)^2) =
\lambda_1(R A R^T)$, where we use the invariance of the trace under orthogonal transformations:
$\operatorname{tr}(\sigma) = \operatorname{tr}(R \sigma R^T)$ for all orthogonal $R$.
Similarly, the basis tensors, being tensor-valued fields,
satisfy $R T_2(A) R^T = T_2(R A R^T)$
(and analogously for the other basis tensors).

This motivates the tensor basis neural network (TBNN)
approach~\cite{lingReynoldsAveragedTurbulence2016}:
the coefficients $\alpha_k$ are predicted by a neural network
$\operatorname{NN}$, i.e., $\alpha = (\alpha_1, \dots, \alpha_{K}) =
\operatorname{NN}(\lambda_1, \dots, \lambda_5)$; the RANS and LES lineage of
the approach is reviewed in \cref{sec:introduction}.

\begin{figure*}
    \centering

    \begingroup

    \def\svgwidth{\textwidth}

    \def\ubar{$\bar{u}$}
    \def\nabubar{$\nabla \bar{u}$}
    \def\A{$A \coloneq \nabla \bar{u} / | \nabla \bar{u} |$}
    \def\a{$a \coloneq | \nabla \bar{u} |$}
    \def\S{$S \coloneq (A + A^T) / 2$}
    \def\R{$W \coloneq (A - A^T) / 2$}

    \def\lone{$\lambda_1 \coloneq \operatorname{tr} S^2$}
    \def\ltwo{$\lambda_2 \coloneq \operatorname{tr} W^2$}
    \def\lthr{$\lambda_3 \coloneq \operatorname{tr} S^3$}
    \def\lfou{$\lambda_4 \coloneq \operatorname{tr} S W^2$}
    \def\lfiv{$\lambda_5 \coloneq \operatorname{tr} S^2 W^2$}

    \def\Tzer{$T_{0} \coloneq I$}
    \def\Tone{$T_{1} \coloneq S$}
    \def\Ttwo{$T_{2} \coloneq S^2$}
    \def\Tthr{$T_{3} \coloneq W^2$}
    \def\Tfou{$T_{4} \coloneq S W - W S$}
    \def\Tfiv{$T_{5} \coloneq W S W$}
    \def\Tsix{$T_{6} \coloneq S^2 W - W S^2$}
    \def\Tsev{$T_{7} \coloneq W S W^2 - W^2 S W$}

    \def\azer{$\alpha_{0}$}
    \def\aone{$\alpha_{1}$}
    \def\atwo{$\alpha_{2}$}
    \def\athr{$\alpha_{3}$}
    \def\afou{$\alpha_{4}$}
    \def\afiv{$\alpha_{5}$}
    \def\asix{$\alpha_{6}$}
    \def\asev{$\alpha_{7}$}

    \def\tau{$m(\bar{u}) \coloneq \Delta^2 a^2 \sum_{i = 1}^7 \alpha_i T_i$}

    \def\tonabubar{Galilean \\ invariance}
    \def\toA{Normalize}
    \def\toSW{Decompose}
    \def\Deltatom{Scaling \\ invariance}
    \def\toT{$O(3)$ equivariance}
    \def\tolambda{$O(3)$ invariance}
    \def\toalpha{Neural network}
    \def\totau{Assemble closure stress}

    \def\endresult{$\forall g \in G, \ g m(\bar{u}) = m(g \bar{u})$}

    \input{figures/ink/TBNN/tbnn.tex}

    \endgroup

    \caption{%
        Tensor-basis neural network. Each arrow stage preserves one symmetry:
        the velocity gradient enforces Galilean invariance, the normalization
        and the prefactor $\Delta^2 a^2$ enforce the $a^2 = b$ scaling
        subgroup, and the invariant inputs $\lambda$ and equivariant basis
        tensors $T_k$ make the assembled stress equivariant to rotations and
        reflections regardless of the network predicting the coefficients
        $\alpha_k$.
    }
    \label{fig:tbnn}
\end{figure*}

In our experiments, we employ the TBNN with the reduced and complete basis of
Stallcup et al.~\cite{stallcupAdaptiveScaleSimilarClosure2022}.
This is illustrated in \cref{fig:tbnn}.
For incompressible flow, we retain only the deviatoric part of the basis.
To improve training stability,
we normalize the VGT by its Frobenius norm as
$A^* \coloneq A / | A |$, the normalized input of \eqref{eq:shared-form},
following Prakash et al.~\cite{prakashInvariantDatadrivenSubgrid2022}.
The model is then the closure \eqref{eq:shared-form} with core
\begin{equation} \label{eq:tbnn}
    m^{*}_\text{TBNN}(A^*) \coloneq \sum_{k = 1}^{7} \alpha_k(\lambda^*) T^{*, \text{dev}}_k,
\end{equation}
where the invariants $\lambda^*$ and basis tensors $T^*$ are computed from
the normalized VGT $A^*$.
The identity tensor $T_0$ is omitted since its deviatoric part vanishes.
The coefficients $(\alpha_k)_k$ are predicted by a
standard feedforward neural network with
$5$ inputs $\lambda^*_1, \dots, \lambda^*_5$ and
$7$ outputs $\alpha_1, \dots, \alpha_7$.
The model is applied locally in physical space,
independently at each grid point,
so $m^\text{TBNN}(\bar{u})(x)$ depends only on $A(x)$.
We implement this as a convolutional neural network with kernel size
$1 \times 1 \times 1$,
meaning the convolutions do not access neighboring points.
The convolutional framework is used for consistency with the G-CNN model
described above, and because it allows efficient batched evaluation
over the entire spatial grid using GPU-accelerated libraries.

\subsection{The filter-scale Reynolds number as an input feature}
\label{sec:redelta-input}

The shared construction above carries no dependence on the viscosity $\nu$:
the three closures are \emph{Reynolds-blind}.
To let them adapt across Reynolds numbers, we define the \emph{global}
filter-scale Reynolds number
\begin{equation} \label{eq:redelta}
    \operatorname{Re}_\Delta \coloneq
    \frac{\Delta^2 \| \bar{A} \|}{\nu},
\end{equation}
where $\| \bar{A} \|$ denotes the volume average of the Frobenius norm over the domain, so
that $\| \bar{A} \|$ is the root-mean-square magnitude of
the resolved velocity gradient and $\operatorname{Re}_\Delta$ is a single
scalar per filtered field.
The global average is deliberate: a pointwise
$\operatorname{Re}_\Delta \propto |\bar{A}(x)|$ would re-inject the local
gradient magnitude that the input normalization $\bar{A} / |\bar{A}|$ removes,
confounding local intermittency with the flow regime.
As defined, the local structure stays in the normalized gradient and the regime
information in the single scalar.

A viscosity input of this kind has been used before: Prakash et
al.~\cite{prakashInvariantDatadrivenSubgrid2022} add the nondimensional
viscosity $\hat{\nu} = \nu / (| \bar{A} | \Delta^2)$ to repair the over-dissipation of
their model at dissipation-range filter widths;
$\hat{\nu}$ is a \emph{pointwise} inverse filter-scale Reynolds number.
The definition \eqref{eq:redelta} differs in the two respects that matter here:
the average is global, for the reason above, and its variation during training
comes from an actual viscosity sweep.
Since the training data of Prakash et al.\ span a single viscosity,
$\hat{\nu}$ varies only through the filter width and the local gradient
magnitude, so the learned $\nu$-dependence is never identified against a
change in viscosity itself---the collinearity that the multi-viscosity
training pool of \cref{sec:results-reynolds} is designed to break.

Each closure is additionally trained in a $+\mathrm{Re}$ variant that receives
the standardized logarithm
$z \coloneq (\log \operatorname{Re}_\Delta - \mu) / \varsigma$ as one extra
input, where the mean $\mu$ and standard deviation $\varsigma$ of
$\log \operatorname{Re}_\Delta$ over the training pool are computed once and
stored with the trained parameters; the same global value is used during
training and at inference.
The scalar enters each architecture in the form its symmetry structure
dictates.
For the MLP it is a tenth input alongside the nine components of
$\bar{A} / |\bar{A}|$, and for the TBNN a sixth invariant alongside
$\lambda^*_1, \dots, \lambda^*_5$.
For the G-CNN, $\operatorname{Re}_\Delta$ is built from the volume-averaged
gradient magnitude and hence invariant under the octahedral group, so it
enters the first layer as a trivial-representation channel: one learnable
weight per output channel, applied identically to all $48$ group components,
exactly like the bias $b_i$. This leaves the octahedral equivariance intact.
In every case the variant adds only a handful of parameters
(\ref{sec:training}).

The feature itself is dictated by the symmetry analysis of
\cref{sec:symmetries-navier-stokes}.
Under the full two-parameter scaling transform, $\Delta \mapsto a \Delta$,
$\bar{A} \mapsto \bar{A} / b$, and $\nu \mapsto (a^2 / b) \nu$, so
$\operatorname{Re}_\Delta$ is \emph{invariant}. The $+\mathrm{Re}$ closure
$m(\bar{A}) = \Delta^2 |\bar{A}|^2 \,
m^*(\bar{A}^*, \operatorname{Re}_\Delta)$, the shared form
\eqref{eq:shared-form} with the extra core input, is
therefore covariant under the full two-parameter group rather than only its
$a^2 = b$ subgroup: supplying $\operatorname{Re}_\Delta$ restores exactly the
viscosity dependence that the Reynolds-blind construction discards.
Galilean invariance is preserved as well, since $\operatorname{Re}_\Delta$ is
built from $\bar{A}$.
The role of this input in generalizing across Reynolds numbers is examined in
\cref{sec:results-reynolds,sec:results-tgv}.

\section{Results} \label{sec:results}

\subsection{Setup and evaluation protocol} \label{sec:results-setup}

We consider forced homogeneous isotropic turbulence in a periodic box
$\Omega = [0, 2 \pi]^3$.
A pseudo-spectral discretization is used for both DNS and LES
(\ref{sec:pseudo-spectral}),
and details of the simulation setup and data generation are given in
\ref{sec:datagen}.
The data-driven closures are trained on this forced flow across a grid of
viscosities $\nu$ and filter ratios $\Delta / h$ (\cref{tab:dns}).
We evaluate them first at a fixed network size and a representative
in-distribution operating point (\cref{sec:results-forced}), then study how
their accuracy depends on network capacity (\cref{sec:results-capacity}), then
test their generalization across Reynolds number with the filter-scale Reynolds
number supplied as an input (\cref{sec:results-reynolds}), then apply
them unchanged to a decaying Taylor--Green vortex (\cref{sec:results-tgv}), and
close by interpreting the shared error floor these experiments expose
(\cref{sec:discussion}).

We compare the closure models of \cref{sec:models}: the classical baselines
No-model, dynamic Smagorinsky, and Clark (\cref{sec:classical}), and the three
data-driven closures $m^\text{MLP}$ (\cref{sec:mlp}), $m^\text{G-CNN}$
(\cref{sec:groupconv}), and $m^\text{TBNN}$ (\cref{sec:tbnn}), which share the
pointwise form \eqref{eq:shared-form} of \cref{sec:construction} with its
input normalization $\bar{A} / |\bar{A}|$ and output prefactor
$\Delta^2 |\bar{A}|^2$.
Each data-driven closure is additionally evaluated in its $+\mathrm{Re}$
variant, which receives the global filter-scale Reynolds number
$\operatorname{Re}_\Delta$ \eqref{eq:redelta} as an extra input
(\cref{sec:redelta-input}), and the MLP also in its symmetrized form
\eqref{eq:mlp-sym}.
All data-driven models are trained using an a priori \emph{structural} loss function,
minimizing the error between predicted and reference stress tensors
computed from discretization-consistent expressions that account for the
numerical scheme (dealiasing, discrete Fourier transforms, etc.).
For details about the training procedure, see \ref{sec:training}.
Each learned closure is trained with five independent initialization seeds;
scalar metrics are reported as mean $\pm$ one standard deviation over seeds,
while field-level figures show the first seed (with seed-spread bands where
indicated).
For legibility, the field-level figures of \cref{sec:results-forced} show the
data-driven closures in their $+\mathrm{Re}$ variants only; at that
in-distribution operating point the Reynolds-blind counterparts behave
near-identically (\cref{tab:errors}).

The closure models are evaluated in two settings:
\begin{itemize}
    \item \emph{A priori}:
        The models predict the SFS using the exact filtered
        velocity field as input.
        We compute $m(\bar{u})$ for all snapshots $\bar{u}$ and models $m$.
    \item \emph{A posteriori}: The models serve as closures in LES
        simulations, initialized from the filtered DNS velocity field.
        The LES solution at time $t$ is denoted $S_t(m, \bar{u}_0)$,
        where $\bar{u}_0$ is the initial condition and $m$ is the closure model.
        The same adaptive time-stepping and forcing schemes as in the DNS are used.
        We compute $S_t$ at all snapshot times $t \in T$, at which the corresponding
        filtered DNS solutions are denoted $\bar{u}_t$.
        Because turbulence is chaotic, the resulting pointwise solution
        error~\eqref{eq:tensor-error-post} saturates once the LES and reference
        trajectories decorrelate; we therefore report its time mean and
        corroborate the model ranking with statistical diagnostics that remain
        informative after decorrelation: energy spectra, the SFS dissipation
        budget, and the distribution of sub-filter dissipation.
\end{itemize}

We report three relative errors:
\begin{itemize}
    \item The a priori closure error
        \begin{equation} \label{eq:tensor-error-prior}
            \frac{1}{|U|} \sum_{u \in U} \frac{\| m(\bar{u}) - \tau(u) \|}{\| \tau(u) \|},
        \end{equation}
        where $U$ denotes the set of all test snapshots and
        $\| T \| \coloneq \bigl( \sum_x | T(x) |^2 \bigr)^{1 / 2}$ is the
        root-sum-square over the grid points of the pointwise Frobenius norm
        $| \cdot |$.
    \item The a posteriori solution error
        \begin{equation} \label{eq:tensor-error-post}
            \frac{1}{|T|} \sum_{t \in T} \frac{\| S_t(m, \bar{u}_0) - \bar{u}_t \|}{\| \bar{u}_t \|},
        \end{equation}
        where $T$ contains all the snapshot times.
    \item The a priori closure equivariance error
        \begin{equation} \label{eq:equi-error-prior}
            \frac{1}{|U|}
            \frac{1}{|G|}
            \sum_{u \in U}
            \sum_{g \in G}
            \frac{\| g m(\bar{u}) - m(g \bar{u}) \|}{\| m(g \bar{u}) \|},
        \end{equation}
        where $g m$ and $g \bar{u}$ denote the group action of $g$ on tensor and vector fields.
\end{itemize}
We additionally report two statistics of the pointwise SFS dissipation rate
\begin{equation} \label{eq:sfs-dissipation}
    \varepsilon_\Delta \coloneq - \tau_{i j} \bar{S}_{i j},
\end{equation}
the rate at which the closure term drains resolved kinetic energy (positive
values are forward transfer, negative values backscatter): the median of
$\varepsilon_\Delta$ normalized by the reference median, and the local
backscatter fraction (the fraction of grid points with
$\varepsilon_\Delta < 0$), the a priori diagnostic introduced by Piomelli et
al.~\cite{piomelliSubgridscaleBackscatterTurbulent1991}.

\subsection{Forced homogeneous isotropic turbulence} \label{sec:results-forced}

We first fix the network size at the parameter-matched tier of roughly
$3000$ parameters (the ``p3000'' tier; \ref{sec:training}) and evaluate the
closures at a representative in-distribution operating point: a held-out forced
realization at a training viscosity ($\nu = 2.5 \e{-4}$,
$\operatorname{Re}_\lambda \approx 246$) and an interpolated filter ratio
($\Delta / h = 2.5$).
The capacity and Reynolds-number dependence are taken up in
\cref{sec:results-capacity,sec:results-reynolds}.

\begin{figure*}
    \centering
    \includegraphics[width=0.8\textwidth]{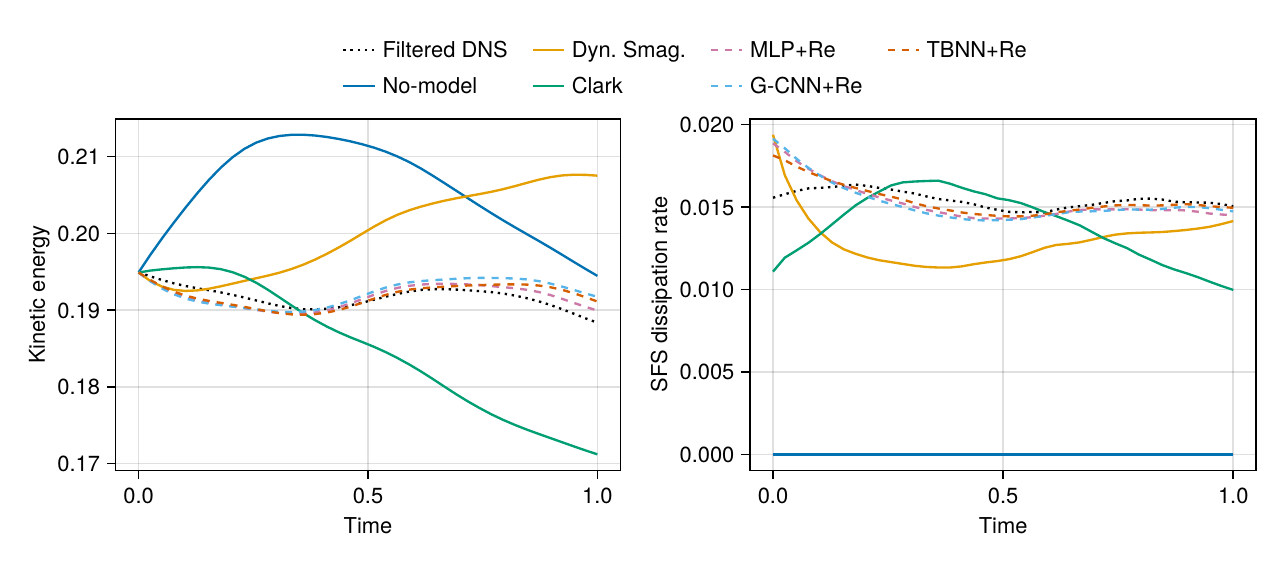}
    \caption{%
        A posteriori evolution of the total kinetic energy and the total SFS
        dissipation rate at the forced in-distribution operating point
        ($\nu = 2.5 \e{-4}$, $\Delta / h = 2.5$); the reference is the filtered
        DNS. As in the other field-level figures of this section, the
        data-driven closures are shown in their $+\mathrm{Re}$ variants
        (\cref{sec:results-setup}).
    }
    \label{fig:budget}
\end{figure*}

\begin{figure*}
    \centering
    \includegraphics[width=0.8\textwidth]{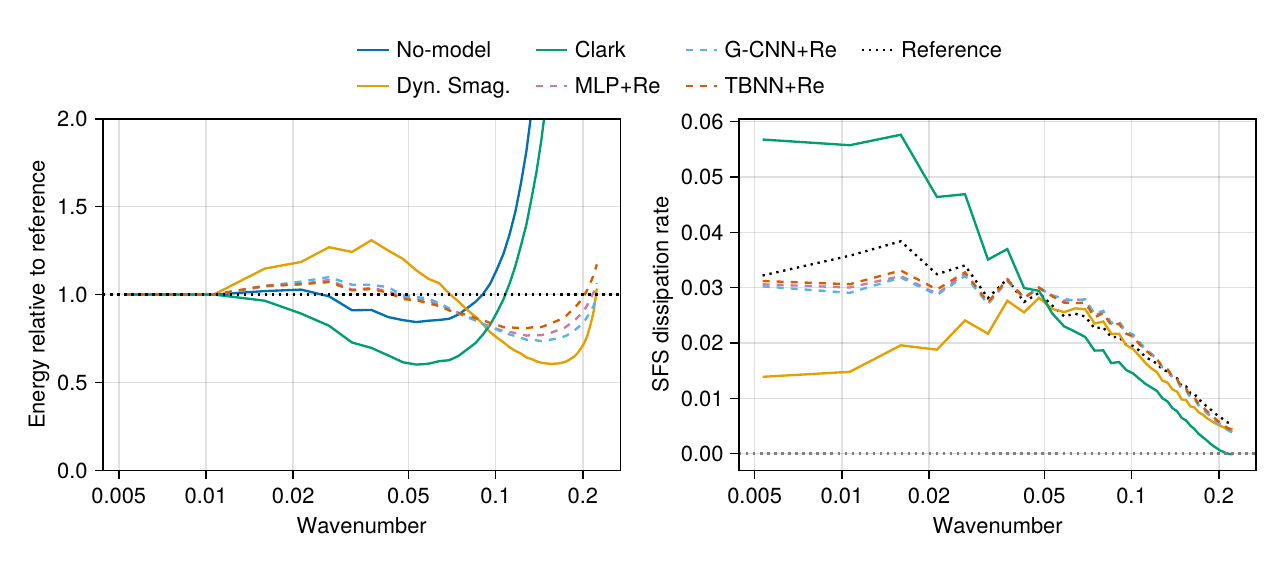}
    \caption{%
        A posteriori spectral diagnostics at the forced in-distribution
        operating point, time-averaged over the evaluation window and computed
        from each model's own LES solution; the wavenumber axis is normalized
        by the Kolmogorov length scale.
        Left: energy spectra divided by the reference (filtered-DNS) spectrum
        on a linear scale. The ratio exposes the per-model deviations; the
        under-dissipative No-model and Clark leave the frame through their
        high-wavenumber pile-up.
        Right: SFS dissipation rate spectrum, normalized by the mean total
        dissipation rate; No-model is absent since its SFS term vanishes.
    }
    \label{fig:spectra}
\end{figure*}

All models ran stably over the evaluation window ($\approx 1 \, t_\text{int}$)
in the a posteriori setting.
Clark survives this window but exhibits non-monotonic dissipation behavior
and net spectral backscatter at the highest wavenumbers
(\cref{fig:budget,fig:spectra}), consistent with the
well-known marginal stability of purely structural closures over long
rollouts.

\begin{table*}
    \centering
    \caption{
        Aggregate metrics for the forced isotropic case at the in-distribution
        operating point ($\nu = 2.5 \e{-4}$, $\Delta / h = 2.5$, ``p3000''
        size): a priori closure error
        \eqref{eq:tensor-error-prior}, a posteriori solution error
        \eqref{eq:tensor-error-post}, a priori
        equivariance error \eqref{eq:equi-error-prior}, median pointwise SFS
        dissipation $\varepsilon_\Delta$ \eqref{eq:sfs-dissipation} normalized
        by the filtered-DNS reference (values above $1$
        are over-dissipative), and the fraction of points producing local
        backscatter ($\varepsilon_\Delta < 0$; the reference fraction is
        $0.22$). The $+\mathrm{Re}$ rows receive the filter-scale Reynolds
        number as an additional input (\cref{sec:redelta-input}).
        Learned-model entries are means over five
        training seeds ($\pm$ one standard deviation). All learned models except
        the unconstrained MLP (both variants) are equivariant to machine
        precision; the classical models are exactly equivariant by construction
        and their equivariance error is not computed (--), while for No-model
        the relative error \eqref{eq:equi-error-prior} is undefined
        ($0 / 0$, N.A.).
    }
    \label{tab:errors}
    \begin{tabular}{l l l l l l}
    \toprule
    Model
                & Closure \eqref{eq:tensor-error-prior}
                & Solution \eqref{eq:tensor-error-post}
                & Equivariance \eqref{eq:equi-error-prior}
                & Median diss.
                & Backscatter \\
    \midrule
    No-model    & $1.0000$            & $0.4325$            & N.A.                       & $0.000$           & $0.0000$ \\
    Dyn. Smag.  & $0.9635$            & $0.3886$            & --                         & $1.508$           & $0.0004$ \\
    Clark       & $0.4658$            & $0.3645$            & --                         & $0.608$           & $0.2514$ \\
    MLP         & $0.4477 \pm 0.0002$ & $0.3529 \pm 0.0012$ & $0.0686 \pm 0.0015$        & $1.355 \pm 0.007$ & $0.1002 \pm 0.0011$ \\
    MLP+Re      & $0.4493 \pm 0.0005$ & $0.3529 \pm 0.0009$ & $0.0817 \pm 0.0029$        & $1.266 \pm 0.006$ & $0.1075 \pm 0.0009$ \\
    G-CNN       & $0.4490 \pm 0.0002$ & $0.3541 \pm 0.0001$ & $(1.165 \pm 0.554)\e{-14}$ & $1.379 \pm 0.004$ & $0.0975 \pm 0.0004$ \\
    G-CNN+Re    & $0.4494 \pm 0.0015$ & $0.3546 \pm 0.0004$ & $(1.386 \pm 0.627)\e{-14}$ & $1.283 \pm 0.012$ & $0.1049 \pm 0.0025$ \\
    TBNN        & $0.4414 \pm 0.0002$ & $0.3533 \pm 0.0002$ & $(1.459 \pm 0.008)\e{-15}$ & $1.340 \pm 0.003$ & $0.0979 \pm 0.0005$ \\
    TBNN+Re     & $0.4460 \pm 0.0002$ & $0.3513 \pm 0.0001$ & $(1.473 \pm 0.007)\e{-15}$ & $1.270 \pm 0.006$ & $0.0967 \pm 0.0010$ \\
    MLP (sym)   & $0.4461 \pm 0.0002$ & $0.3533 \pm 0.0004$ & $(1.437 \pm 0.004)\e{-15}$ & $1.356 \pm 0.007$ & $0.0996 \pm 0.0010$ \\
    MLP (sym)+Re & $0.4470 \pm 0.0004$ & $0.3532 \pm 0.0004$ & $(1.454 \pm 0.004)\e{-15}$ & $1.267 \pm 0.006$ & $0.1065 \pm 0.0010$ \\
    \bottomrule
\end{tabular}

\end{table*}

In \cref{tab:errors} we collect these metrics for every closure at the forced
operating point.

\paragraph{Prediction errors}
No-model is worst by construction ($\| 0 - \tau \| / \| \tau \| = 100\%$).
Dynamic Smagorinsky has an a priori error of $96\%$: as a functional model it is
designed to match energy dissipation rather than the stress tensor itself, yet it
still attains a far lower a posteriori error than No-model because it supplies
sub-filter dissipation.
The data-driven models are the most accurate in tensor prediction: at this
saturated size all learned variants cluster within one percentage point of one
another ($0.441$--$0.449$), the TBNN the most accurate ($0.441$), the
unconstrained MLP ($0.448$) essentially tied with the G-CNN ($0.449$), and all
clearly ahead of Clark ($0.466$; \cref{tab:errors}).
For a posteriori solution error the data-driven models are again essentially
tied, all within about $1\%$ of one another near $0.353$, ahead of Clark
($0.364$), dynamic Smagorinsky ($0.389$), and No-model ($0.432$).
Across the five training seeds the spread is small for every architecture at
this saturated size ($\lesssim 0.002$ in a priori error, \cref{tab:errors});
below saturation the unconstrained MLP is far less reproducible, with seed
bands an order of magnitude wider than the constrained models'
(\cref{fig:saturation}). The symmetry constraints thus do not improve
saturated mean accuracy, but they make the trained closure reproducible
across initializations already at small capacity.
That three such differently constructed networks collapse to nearly the same a
priori and a posteriori errors suggests they approximate the same underlying
pointwise map $\bar{A}(x) \mapsto \tau(x)$---the optimal structural closure in
the relative-tensor metric they were trained on.
This collapse is specific to sufficient capacity; how it emerges as the networks
grow is the subject of \cref{sec:results-capacity}, and we develop the
optimal-closure interpretation and its limits in \cref{sec:discussion}.

\paragraph{Closure equivariance errors}
All models except the unconstrained MLP have equivariance errors at machine
precision, confirming exact equivariance.
The MLP has a substantial a priori equivariance error of about $7\%$ (and the
$+\mathrm{Re}$ MLP about $8\%$), so it has not learned exact roto-reflectional
equivariance despite the isotropic training data.
Isotropy of the flow makes the joint input--output distribution invariant only
on average; it does not force a finite-data fit to be equivariant pointwise.
We did not apply octahedral augmentation to the MLP during training: on a
statistically isotropic training set the rotated and reflected copies are
already draws from the same distribution, so augmentation could at most
reduce estimation variance; the symmetrized MLP \eqref{eq:mlp-sym} provides
the stronger check, the exactly equivariant projection of the same learned
function.
Averaging the trained MLP over the octahedral group drives its equivariance
error to machine precision while leaving its accuracy essentially unchanged:
marginally improved a priori ($0.446$ versus $0.448$) and unchanged a
posteriori ($0.353$; \cref{tab:errors}), with the same parity holding across
network sizes (\cref{fig:saturation}).
The equivariance violation is therefore an underused inductive bias rather than a
fundamental limitation: the function the MLP learned already lies close to the
equivariant subspace.
This is also why the purpose-built equivariant models gain no measurable
accuracy advantage here: the unconstrained optimum is already nearly
equivariant, and the large, isotropic training set lets even the MLP approach
it directly from data (\cref{sec:conclusion}).
The MLP still preserves Galilean invariance and the $a^2 = b$ scaling subgroup
through its use of $\bar{A}$ (rather than $\bar{u}$) as input and the output
scaling by $\Delta^2 |\bar{A}|^2$.

\paragraph{Energy and dissipation evolution}
\Cref{fig:budget} shows the a posteriori evolution of the total kinetic energy
and total SFS dissipation rate (excluding the viscous dissipation rate)
over the evaluation window.
The reference kinetic energy stays nearly constant; No-model accumulates
energy rapidly due to the absence of any subgrid drain; dynamic Smagorinsky
follows the reference at first but then drifts upward, consistent with
its a posteriori dissipation falling below the reference in
the right panel of \cref{fig:spectra}.
The data-driven models track the reference kinetic energy most closely,
and their SFS dissipation rates remain close to (and slightly below)
the reference throughout the window;
they are visually indistinguishable on this plot.
Clark exhibits non-monotonic behavior, peaking at intermediate times
and falling thereafter, signaling its marginal stability in this
configuration.

\paragraph{Energy spectra}
The left panel of \cref{fig:spectra} shows the time-averaged energy spectra
for the different LES solutions.
Time averaging reduces temporal fluctuations,
providing a clearer comparison of the spectral energy distribution across models.
The wavenumber axis is normalized by the Kolmogorov length scale
(see \cref{eq:kolmogorov}), indicating how far the resolved scales are from the
dissipation range; the energy normalization cancels in the ratio.
No-model and Clark are insufficiently dissipative, exhibiting a spurious energy pile-up at the highest wavenumbers.
Dynamic Smagorinsky preserves the spectral shape but sits slightly above the reference
at intermediate wavenumbers, consistent with its excess of forward transfer there.
The data-driven models all closely match the reference spectrum, with only a slight
deficit at the highest resolved wavenumbers consistent with their mild
over-dissipation.

\paragraph{SFS dissipation rate spectra}
The right panel of \cref{fig:spectra} shows the time-averaged SFS dissipation spectra
for the different LES solutions.
These are computed \emph{a posteriori}, in the sense that the dissipation rate
is computed from the LES solution produced by the given model.
Clark is excessively dissipative at low wavenumbers and produces a \emph{net}
backscatter at the highest wavenumbers, consistent with the
pile-up at high wavenumbers seen in the left panel.
Dynamic Smagorinsky sits well below the reference at all wavenumbers,
underestimating SFS dissipation in the a posteriori setting despite
its a priori median being about $1.5$ times the reference;
this discrepancy reflects the well-known sensitivity of dissipation
to the input field, which differs between the filtered DNS and the
self-produced LES trajectory.
The data-driven models closely match the reference across the resolved range.

\begin{figure}
    \centering
    \includegraphics[width=0.9\columnwidth]{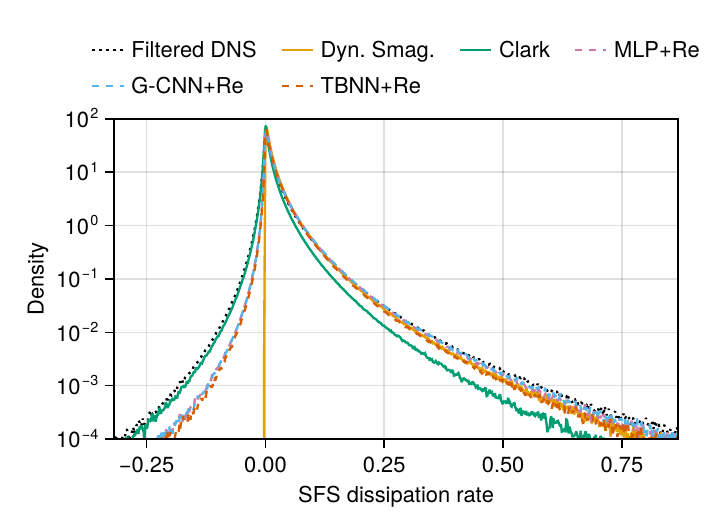}
    \caption{%
        A priori distribution of the pointwise SFS dissipation
        $\varepsilon_\Delta = - \tau_{i j} \bar{S}_{i j}$
        \eqref{eq:sfs-dissipation} (kernel density estimate over all snapshots),
        with every closure evaluated on the same filtered DNS fields. Positive
        values are forward transfer, negative values backscatter; No-model is
        omitted.
    }
    \label{fig:dissipation-density}
\end{figure}

\paragraph{Dissipation distribution and backscatter}
Beyond the aggregate errors, the distribution of the pointwise SFS dissipation
$\varepsilon_\Delta$ \eqref{eq:sfs-dissipation} separates the closures
(\cref{fig:dissipation-density}).
Dynamic Smagorinsky produces no backscatter (negative values) by construction
(the residual $0.04\%$ fraction in \cref{tab:errors} is a dealiasing artifact
of evaluating the model stress pseudo-spectrally) and concentrates its
forward transfer in a narrow band.
Clark captures the backscatter tail, but
lacks the heavy forward-transfer tail of the reference.
The data-driven models reproduce both tails closely, slightly underestimating
backscatter and slightly overestimating forward transfer.
This is summarized by the last two columns of \cref{tab:errors}.
The median pointwise SFS dissipation, normalized by the reference, is about
$1.5$ for dynamic Smagorinsky and $1.3$--$1.4$ for the data-driven models
(over-dissipative), about $0.6$ for Clark (under-dissipative), and zero for
No-model; the $+\mathrm{Re}$ input lowers the data-driven values slightly
($1.27$--$1.28$).
The backscatter fraction is $0.22$ for the reference, essentially zero for
No-model and dynamic Smagorinsky, $0.25$ for Clark (slightly too much), and about
$0.10$--$0.11$ for the data-driven models (roughly half the reference).
Both the mild over-dissipation and the backscatter deficit grow out of
distribution; \cref{sec:results-reynolds} quantifies this and shows that the
$\operatorname{Re}_\Delta$ input corrects the calibration.

\begin{table*}
    \centering
    \caption{
        Training and inference time in seconds, across the size tiers, for the
        $+\mathrm{Re}$ closures and the classical baselines.
        For perspective, generating the filtered-DNS reference trajectory for
        this evaluation point (DNS at $810^3$ over the same window, including
        filtering and diagnostics) took about $9800$ s on the same hardware:
        the LES rollouts are $60$ to $900$ times faster, depending on the
        closure.
    }
    \label{tab:timing}
    \begin{tabular}{l l r r r}
    \toprule
    Model
                & Tier
                & Parameters
                & Training $[s]$
                & Inference $[s]$ \\
    \midrule
    No-model   & --    &     0 &              $0$ &          $11.4$ \\
    Dyn. Smag. & --    &     0 &              $0$ &          $34.4$ \\
    Clark      & --    &     0 &              $0$ &          $17.8$ \\
    MLP+Re     & p120  &   125 &   $71.8 \pm 1.4$ &  $49.2 \pm 5.7$ \\
    MLP+Re     & p400  &   374 &   $75.5 \pm 0.6$ &  $45.1 \pm 0.7$ \\
    MLP+Re     & p1200 & 1,108 &   $85.8 \pm 0.4$ &  $48.3 \pm 0.5$ \\
    MLP+Re     & p3000 & 2,962 &   $95.8 \pm 0.7$ &  $53.6 \pm 0.3$ \\
    G-CNN+Re   & p120  &   118 &  $113.0 \pm 0.9$ &  $62.1 \pm 0.9$ \\
    G-CNN+Re   & p400  &   428 &  $141.9 \pm 0.9$ &  $83.6 \pm 0.9$ \\
    G-CNN+Re   & p1200 & 1,083 & $195.4 \pm 11.1$ & $113.8 \pm 1.2$ \\
    G-CNN+Re   & p3000 & 3,009 &  $275.0 \pm 1.5$ & $169.3 \pm 0.9$ \\
    TBNN+Re    & p120  &   119 &  $127.6 \pm 2.8$ &  $32.2 \pm 0.8$ \\
    TBNN+Re    & p400  &   364 &  $133.2 \pm 7.9$ &  $34.4 \pm 0.0$ \\
    TBNN+Re    & p1200 & 1,100 &  $144.3 \pm 8.7$ &  $37.0 \pm 0.4$ \\
    TBNN+Re    & p3000 & 2,918 &  $154.1 \pm 6.6$ &  $42.9 \pm 0.5$ \\
    MLP        & p3000 & 2,940 &   $91.1 \pm 0.6$ &  $44.7 \pm 0.7$ \\
    G-CNN      & p3000 & 3,005 &  $258.1 \pm 0.8$ & $153.3 \pm 2.2$ \\
    TBNN       & p3000 & 2,896 &  $145.0 \pm 0.3$ &  $39.2 \pm 0.1$ \\
    \bottomrule
\end{tabular}

\end{table*}

\paragraph{Computation time}
\Cref{tab:timing} reports training and inference times (the latter for solving
the LES equations over the evaluation window).
Only the three data-driven models require training; the dynamic-Smagorinsky
coefficient is recomputed from the resolved field at every step.
Among the cheap baselines, No-model and Clark are essentially free
($\sim 10$--$20$\,s), while dynamic Smagorinsky is about twice as expensive
because of the Germano test-filtering.
At the p3000 size the MLP is about four times slower than No-model at inference;
the G-CNN is roughly three and a half times slower than the MLP, owing to the
$48$-channel regular representation it carries despite an equal parameter
count, while the TBNN, which enforces equivariance through the tensor basis
rather than the architecture, is slightly faster than the MLP.
The $+\mathrm{Re}$ input adds a single channel and changes these costs
negligibly.
Part of the data-driven overhead is implementation-specific: unlike the
classical closures, our network evaluation reallocates each layer output at every
call, a penalty that grows with the intermediate array sizes and is therefore
largest for the G-CNN; preallocation would narrow the gap.
Both training and inference cost grow with network size (\cref{tab:timing}).
This makes the parameter efficiency of the equivariant and tensor-basis
models, quantified next (\cref{sec:results-capacity}), directly relevant to
cost: a constrained model that reaches the accuracy floor at a smaller size
is cheaper to train and to run.

\subsection{Network capacity and inductive bias} \label{sec:results-capacity}

\begin{figure*}
    \centering
    \includegraphics[width=\textwidth]{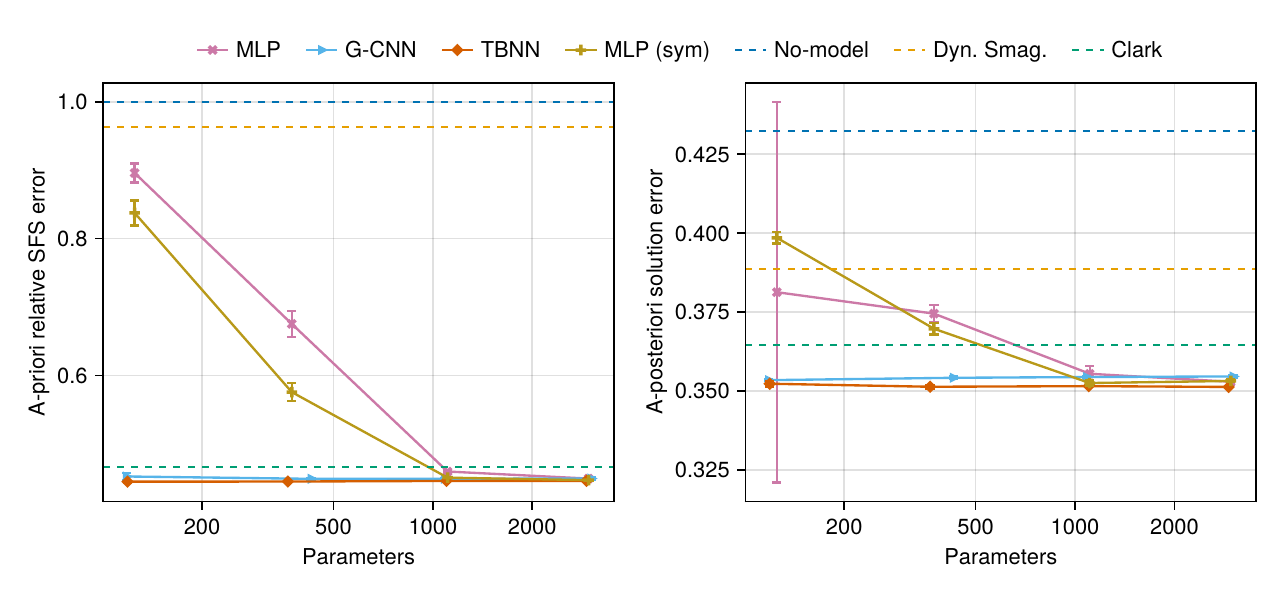}
    \caption{%
        A priori relative SFS error (left) and a posteriori solution error
        (right) versus the number of trainable parameters, at the
        in-distribution operating point. Each data-driven architecture is shown
        across its size grid; the symmetrized MLP \eqref{eq:mlp-sym} and the
        classical baselines (horizontal lines) are included for reference. Bands
        span the five training seeds.
    }
    \label{fig:saturation}
\end{figure*}

The forced comparison above fixes the network size; varying it isolates the role
of the inductive bias from that of raw capacity.
\Cref{fig:saturation} shows the a priori and a posteriori errors as the
parameter count is swept across more than an order of magnitude, with the three
architectures parameter-matched at each size.
The tensor-basis TBNN and the group-convolutional G-CNN sit at the same accuracy
floor at every size, including the smallest ($\sim 120$ parameters): their
built-in structure leaves little to learn, so they are already saturated.
The unconstrained MLP, by contrast, is far worse when small, even less
accurate than the classical Clark model at the smallest sizes, and reaches
the floor only at the largest size---a factor of $25$ more parameters than
the constrained closures need.
Even that factor is only a lower bound set by the tested grid: the direct
estimate below shows that seven \emph{constant} tensor-basis coefficients
already realize nearly all of the attainable accuracy, so the constrained
designs' true requirement is smaller still.
The symmetrized MLP tracks the raw MLP throughout, confirming that this gap is a
matter of capacity rather than of the residual equivariance error.
All three architectures converge to the same floor in both metrics, up to
percent-level differences that are resolvable across seeds but immaterial
a posteriori (\cref{tab:errors}): the inductive bias buys parameter
efficiency, not a lower attainable error.
We identify this common floor as the one-point optimal closure of Langford
and Moser~\cite{langfordOptimalFormulationsIsotropic1999}: the best predictor
of $\tau$ from the single-point gradient $\bar{A}$, which any sufficiently
flexible pointwise closure approaches and none can pass
(\cref{sec:discussion}).

\begin{figure*}
    \centering
    \includegraphics[width=\textwidth]{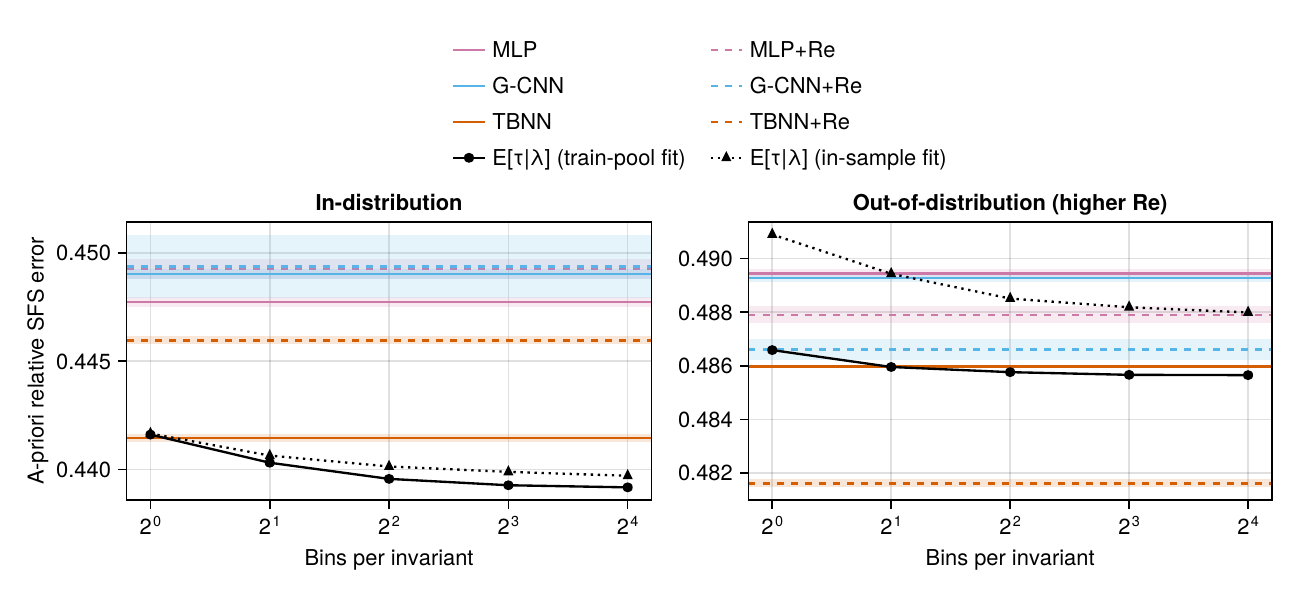}
    \caption{%
        A priori relative SFS error of the directly estimated one-point
        optimal closure versus the bin resolution, at the in-distribution
        (left) and higher-Reynolds out-of-distribution (right) test points.
        The estimator is the conditional mean of the normalized stress given
        the normalized gradient, obtained by per-bin least-squares fits of the
        seven tensor-basis coefficients over quantile bins of the gradient
        invariants; it is fitted either on the training data (circles) or
        in-sample on the test series itself (triangles). Horizontal lines mark the
        seed-mean a priori errors of the trained closures at the matched
        capacity tier (bands: $\pm$ one standard deviation over the five
        seeds); solid lines are the Reynolds-blind closures, dashed lines their
        $+\mathrm{Re}$ counterparts. The Reynolds-blind TBNN sits on the
        estimate at both points; out of distribution, the $+\mathrm{Re}$
        TBNN dips below the optimum attainable from the normalized gradient
        alone.
    }
    \label{fig:condmean}
\end{figure*}

This identification can be checked by direct estimation, because the optimal
one-point closure is a conditional average that requires no training.
Over \emph{all} functions of the normalized gradient, not just those a given
architecture can represent, the minimizer of the mean-square training
objective (\ref{sec:training}) is the conditional mean
$\mathbb{E}[\tau / (\Delta^2 |\bar{A}|^2) \mid \bar{A} / |\bar{A}|]$:
the average normalized stress over all flow configurations sharing that value
of the normalized gradient.
This conditional mean is the \emph{optimal estimator} of a priori analysis,
introduced by Moreau et
al.~\cite{moreauOptimalEstimationLargeeddy2006} as the most accurate model
attainable from a chosen set of input variables and conventionally computed,
as below, by binning those inputs; its residual is the \emph{irreducible
error} of the input set, which Vollant et
al.~\cite{vollantSubgridscaleScalarFlux2017} used both to select inputs for
neural subgrid closures and as the accuracy a trained network should attain.
Since the flow ensemble is statistically isotropic, this conditional mean is an
equivariant function of its argument and therefore takes the tensor-basis form
\eqref{eq:tbnn}: seven coefficients depending only on the invariants
$\lambda^*$, of which four are independent
(normalization fixes $\lambda^*_1 - \lambda^*_2 = |A^*|^2 = 1$).
This makes it directly measurable by counting.
We partition the four-dimensional invariant space into equal-population
(quantile) bins, $n_q$ bins per invariant, and in each bin fit the seven
coefficients of \cref{eq:tbnn} by least squares to the samples that fall in
it.
This is a histogram-style regression with no trained parameters; its only
assumption is that the coefficients vary little across a bin.
(Bins too sparse to determine seven coefficients inherit those of the
enclosing bin of a coarser partition, which affects well under $0.1\%$ of the
samples.)
Refining $n_q$ trades this resolution bias against sampling noise, so scanning
$n_q$ is a built-in robustness check.
\Cref{fig:condmean} shows the resulting a priori error, with the estimator fit
on the same pooled training data as the networks and evaluated on the test
snapshots with the same error metric.
The estimate is remarkably insensitive to the bin resolution: from a single
global bin (a \emph{constant-coefficient} tensor-basis model, seven fitted
numbers) to $16^4$ bins, the error only decreases from $0.442$ to $0.439$,
and refitting the estimator in-sample on the test data itself moves it by
less than $0.001$.
Out of distribution the in-sample fit lands slightly \emph{above} the
training-pool fit (by up to $0.005$), which reflects the error measure rather
than the estimate: the per-bin least squares minimizes the uniformly weighted
error in the normalized stress---the same objective the networks train
on---whereas the plotted relative error weights samples by $|\bar{A}|^4$, and
refitting in-sample with metric-matched weights indeed restores the expected
ordering.
The saturated networks sit on this estimate: the TBNN floor of
$0.4414 \pm 0.0002$ lies within $0.6\%$ of the refined value, and the MLP and
G-CNN floors within $2.3\%$.
That the TBNN lands closest is expected: the conditional mean takes exactly
the tensor-basis form \eqref{eq:tbnn}, as used above, so the optimum lies in
the TBNN's hypothesis space by construction, whereas the MLP and G-CNN must
approximate the basis expansion through generic layers of finite width.
The shared floor of \cref{fig:saturation} therefore \emph{is} the one-point
optimal closure, measured rather than inferred.
The flatness of the curve is itself informative: nearly all of the attainable
accuracy is already realized by the seven constant coefficients, so the
dependence of the optimal coefficients on the invariants is worth only about
$0.002$ in this error measure.
The floor is set by the irreducible one-point uncertainty, not by the
complexity of the coefficient functions.
At the held-out higher-Reynolds-number test point
(\cref{fig:condmean}, right) the same training-pool estimate gives $0.486$,
again coinciding with the Reynolds-blind TBNN ($0.4860 \pm 0.0000$), while
its $+\mathrm{Re}$ counterpart reaches $0.482$, \emph{below} what any
function of the normalized gradient alone attains here.
This is direct evidence that the $\operatorname{Re}_\Delta$ input carries
information beyond the one-point normalized gradient
(\cref{sec:results-reynolds}).
Because training and inference cost grow with size (\cref{tab:timing}), reaching
this floor at a smaller size is a concrete advantage of the equivariant and
tensor-basis designs, most valuable when data or parameters are scarce.

\subsection{Generalization across Reynolds number} \label{sec:results-reynolds}

\begin{table*}
    \centering
    \caption{%
        Window-time-averaged forced-HIT DNS statistics for the training and test
        viscosities: Taylor-microscale and integral Reynolds numbers, integral
        time $t_\text{int} = L_\text{int} / u'$, Kolmogorov scale $\eta$,
        resolution $k_{\max}\eta$, and
        dissipation $\varepsilon$. The closures are trained on the three training
        viscosities and tested on a held-out realization at a training viscosity
        and on a held-out higher Reynolds number.
        Each training DNS is filtered at the ratios
        $\Delta / h \in \{2, 3, 4\}$ and each test DNS at the interpolated and
        extrapolated ratios $\Delta / h \in \{2.5, 3.5, 5.0\}$
        (\ref{sec:datagen}), yielding the filter-scale Reynolds numbers of
        \cref{fig:trend-redelta}.
    }
    \label{tab:dns}
    \begin{tabular}{l l l l l l l l}
    \toprule
    Dataset
                & $\nu$
                & $\mathrm{Re}_\lambda$
                & $\mathrm{Re}_L$
                & $L_\text{int}/u'$
                & $\eta$
                & $k_{\max}\eta$
                & $\varepsilon$ \\
    \midrule
    Train          & $1.50\e{-4}$ & $320$ &      $6822$ & $7.18$ & $3.61\e{-3}$ & $1.46$ & $0.0200$ \\
    Train          & $2.50\e{-4}$ & $249$ &      $4128$ & $7.31$ & $5.33\e{-3}$ & $2.16$ & $0.0194$ \\
    Train          & $4.00\e{-4}$ & $199$ &      $2631$ & $7.57$ & $7.68\e{-3}$ & $3.11$ & $0.0184$ \\
    Test (in-dist) & $2.50\e{-4}$ & $246$ &      $4028$ & $7.19$ & $5.32\e{-3}$ & $2.16$ & $0.0195$ \\
    Test (OOD)     & $1.00\e{-4}$ & $394$ & $1.04\e{4}$ & $7.33$ & $2.68\e{-3}$ & $1.09$ & $0.0193$ \\
    \bottomrule
\end{tabular}

\end{table*}

\begin{figure*}
    \centering
    \includegraphics[width=\textwidth]{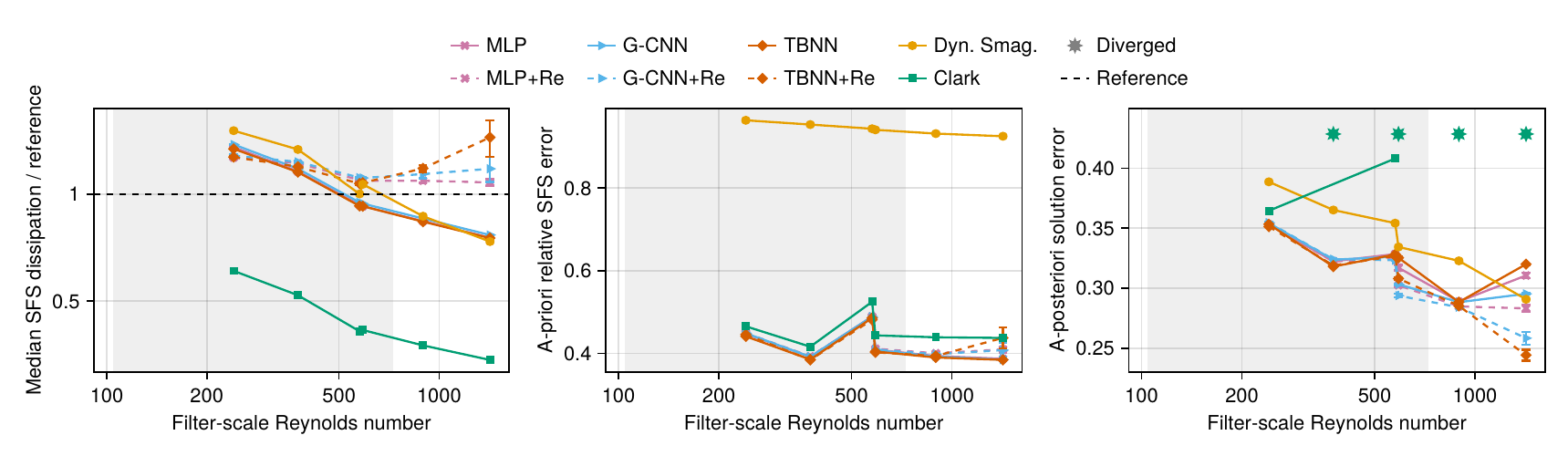}
    \caption{%
        Generalization across the filter-scale Reynolds number
        $\operatorname{Re}_\Delta$: median SFS dissipation normalized by the
        reference (left), a priori relative SFS error (middle), and a posteriori
        solution error (right) versus $\operatorname{Re}_\Delta$, over the test
        grid of viscosities and filter ratios. Each test case is placed at the
        time mean of its global $\operatorname{Re}_\Delta$ \eqref{eq:redelta}
        over the snapshot series. Solid lines are Reynolds-blind
        closures, dashed lines their $+\mathrm{Re}$ counterparts; the shaded band
        marks the training range of $\operatorname{Re}_\Delta$. Stars mark
        diverged a posteriori runs.
    }
    \label{fig:trend-redelta}
\end{figure*}

The data-driven closures so far have no explicit Reynolds-number dependence: they
regress the normalized stress $\tau / (\Delta^2 |\bar{A}|^2)$ from the normalized
gradient $\bar{A} / |\bar{A}|$ alone, so their coefficient is pinned to the
training regime.
To test how this transfers, we trained across three viscosities
(\cref{tab:dns}, $\operatorname{Re}_\lambda \approx 199$--$320$) and three
filter ratios, and evaluate on held-out data spanning the filter-scale
Reynolds number $\operatorname{Re}_\Delta$ \eqref{eq:redelta} over nearly a
decade: a new realization at a training viscosity, a held-out
higher-Reynolds-number viscosity ($\operatorname{Re}_\lambda \approx 394$),
and interpolated and extrapolated filter ratios.
\Cref{fig:trend-redelta} traces three diagnostics against
$\operatorname{Re}_\Delta$.
The a priori tensor error (middle) is essentially flat in
$\operatorname{Re}_\Delta$ and barely changed by the $+\mathrm{Re}$ input: the
\emph{structure} of the SFS is close to Reynolds-independent.
The median dissipation (left) is not: a Reynolds-blind closure slides from
over-dissipation at low $\operatorname{Re}_\Delta$ to under-dissipation at high
$\operatorname{Re}_\Delta$, drifting away from the reference as the flow leaves
the training band.
Supplying $\operatorname{Re}_\Delta$ largely corrects this: the
$+\mathrm{Re}$ TBNN and G-CNN stay far closer to the reference out into the
extrapolation range, although the $+\mathrm{Re}$ TBNN overshoots to about
$1.5$ times the reference at the highest $\operatorname{Re}_\Delta$.
The correction carries over to the a posteriori solution error (right), where
the $+\mathrm{Re}$ closures attain the lowest errors at the highest Reynolds
numbers while the parameter-free Clark model diverges at several points.
The picture is consistent: the over-dissipation is a calibration effect, set by
how much energy lies below the filter, which is a function of
$\operatorname{Re}_\Delta$; the tensor \emph{structure} is not, so feeding the
single scalar $\operatorname{Re}_\Delta$ is enough to restore the calibration
without disturbing the structure.
As with capacity, the equivariant and non-equivariant architectures do not
separate along this axis: the gain comes from the input feature, available to
all three, not from the architecture.

\subsection{Generalization to decaying Taylor--Green flow} \label{sec:results-tgv}

\begin{figure*}
    \centering
    \includegraphics[width=\textwidth]{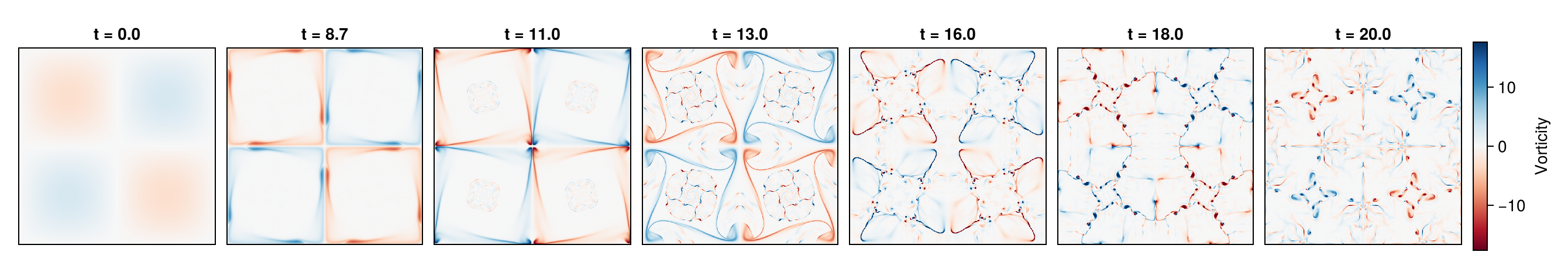}
    \caption{%
        Time evolution of the decaying Taylor--Green vortex at nominal integral
        $\operatorname{Re} = 6000$, visualized through the out-of-plane
        vorticity $\bar{\omega}_z$ of the \emph{filtered} DNS velocity on a
        two-dimensional $z$-section (only the filtered field is stored). Panels
        run from the laminar initial condition (left), through the transition to
        turbulence near peak dissipation ($t^* \approx 9$, centre-left), into the
        late viscous decay (right); time is normalized as
        $t^* = t V_0 / L_\text{TGV}$ and
        the diverging color scale is shared and zero-centered across panels, so
        the fading amplitude reflects the true decay.
    }
    \label{fig:tgv-fields}
\end{figure*}

To probe how the closures behave outside the regime they were trained on, we
apply them \emph{unchanged} to a decaying Taylor--Green vortex (TGV): an
out-of-distribution flow that is unforced, transitional, and, even at its most
intense, less turbulent than the forced training data.
Only the closure-relevant parameters are held fixed (the viscosity $\nu$,
the LES resolution $M$ and filter width $\Delta$, and the trained
network parameters); everything about the flow changes.
The initial condition is the analytic Taylor--Green field, there is no forcing,
and the flow passes through a laminar roll-up, a transition to turbulence, and
a long viscous decay---none of which appear in the forced training data
(\cref{fig:tgv-fields}).
This reuse is legitimate because the closures are amplitude- and
Reynolds-invariant \emph{in form}: the data-driven models regress the
normalized stress $\tau / (\Delta^2 |\bar{A}|^2)$ from the normalized
gradient $\bar{A} / |\bar{A}|$ (plus, for the $+\mathrm{Re}$ variants, the
filter-scale Reynolds number $\operatorname{Re}_\Delta$ \eqref{eq:redelta}),
so reusing them requires only that $\nu$, $M$, and $\Delta$ match training.

\begin{figure*}
    \centering
    \includegraphics[width=\textwidth]{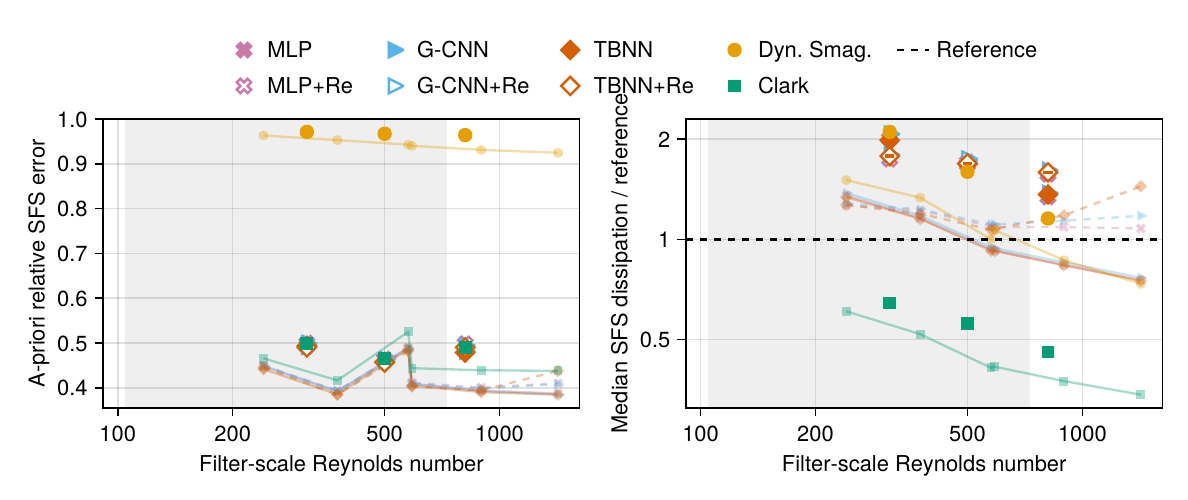}
    \caption{%
        Taylor--Green flow placed on the forced $\operatorname{Re}_\Delta$ trend
        of \cref{fig:trend-redelta}: a priori relative SFS error (left) and median
        SFS dissipation normalized by the reference (right) versus
        $\operatorname{Re}_\Delta$. Faded lines are the forced closures
        (Reynolds-blind solid, $+\mathrm{Re}$ dashed); markers are the
        Taylor--Green values, one per filter ratio
        $\Delta / h \in \{2.5, 3.5, 5.0\}$ (placed at the peak-dissipation-instant
        $\operatorname{Re}_\Delta$), Reynolds-blind versus $+\mathrm{Re}$.
        In tensor error the Taylor--Green points lie close to the forced curve,
        slightly above it; in dissipation they sit well \emph{above} it (the
        flow-type residual), and the $+\mathrm{Re}$ input pulls the
        tight-filter points back toward the curve.
    }
    \label{fig:tgv-redelta}
\end{figure*}

We run a single decaying DNS at the central training viscosity $\nu = 2.5\e{-4}$,
initialized from the analytic Taylor--Green field \eqref{eq:tgv-ic} at a nominal
integral Reynolds number
$\operatorname{Re} = V_0 L_\text{TGV} / \nu = 6000$, fixed by the initial
amplitude $V_0 = 1.5$ and the unit characteristic length
$L_\text{TGV} = 1$ of the vortex (\ref{sec:datagen}).
The \emph{realized} turbulence is far milder. Because the flow is
non-stationary, and the measured integral Reynolds number diverges in the
laminar roll-up and the late decay where the dissipation vanishes, we
characterize it at the instant of peak dissipation: the most intense and
best-developed state, and the one that dominates the sub-filter dissipation.
At that instant the flow reaches only $\operatorname{Re}_\lambda \approx 157$,
\emph{below} the forced training band ($\operatorname{Re}_\lambda \approx 199$--$320$,
\cref{tab:dns}), while the DNS stays well resolved.
So despite its larger nominal number, the Taylor--Green case is both decaying
and, even at its most turbulent, genuinely \emph{less turbulent} than the data
the closures were trained on.
We evaluate every closure on this same flow at the three test filter ratios
$\Delta / h \in \{2.5, 3.5, 5.0\}$.
Evaluated at the peak-dissipation instant, these sweep the filter-scale
Reynolds number $\operatorname{Re}_\Delta$ \eqref{eq:redelta} across
$313$--$813$, within the range covered by the forced grid of
\cref{sec:results-reynolds}, so the Taylor--Green flow can be placed on the
same $\operatorname{Re}_\Delta$ axis as that grid.
Unless noted, the figures and \cref{tab:errors-tgv} report the
interpolated filter ratio $\Delta / h = 2.5$, which matches the forced operating
point of \cref{sec:results-forced} and makes the two flows directly comparable;
the other two filter ratios enter through \cref{fig:tgv-redelta} and the text.

\begin{figure*}
    \centering
    \includegraphics[width=\textwidth]{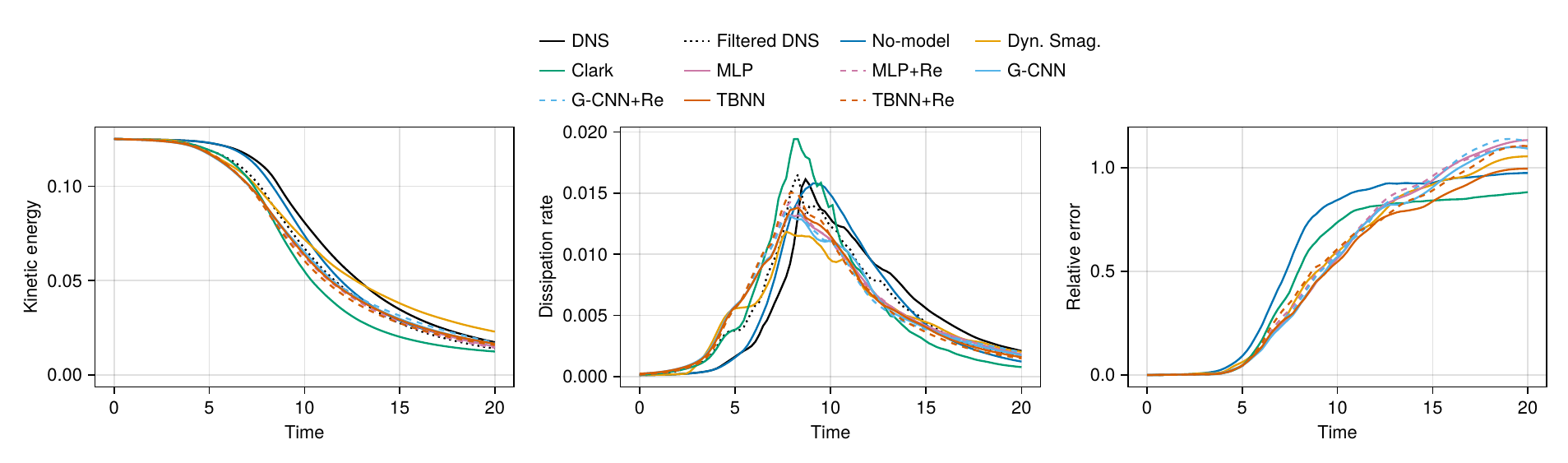}
    \caption{%
        Decaying Taylor--Green vortex at nominal integral $\operatorname{Re} = 6000$
        and filter ratio $\Delta / h = 2.5$: resolved kinetic energy (left),
        total dissipation rate (centre), and relative LES solution error
        \eqref{eq:tensor-error-post} (right) versus the normalized time
        $t^* \coloneq t V_0 / L_\text{TGV}$, for the full DNS, the filtered DNS,
        and each LES
        closure (Reynolds-blind solid, $+\mathrm{Re}$ dashed). For the LES cases,
        the dissipation includes both viscous and closure contributions; the
        solution error is measured against the filtered DNS, so it is shown for
        the closures only. The data-driven models hold the lowest error through
        the transition; Clark catches up only in the long decaying tail.
    }
    \label{fig:dissipation-tgv}
\end{figure*}

\begin{table*}
    \centering
    \caption{%
        Aggregate errors for the decaying Taylor--Green vortex at nominal integral
        $\operatorname{Re} = 6000$ and filter ratio $\Delta / h = 2.5$.
        ``Median diss.'' is the median pointwise SFS dissipation normalized by the
        filtered-DNS reference (reference $= 1$); the reference backscatter
        fraction is $0.304$.
        The $+\mathrm{Re}$ rows receive the filter-scale Reynolds number as an
        additional input. The dynamic-Smagorinsky cross-correlation is undefined
        (\,--\,): its dynamic coefficient clips to zero on the laminar
        Taylor--Green snapshots, so the predicted stress vanishes and the
        per-snapshot correlation is $0 / 0$. For No-model the cross-correlation
        is equally undefined (N.A.), its predicted stress being identically
        zero.
        Learned-model entries are means over five training seeds
        ($\pm$ one standard deviation).
        The corresponding tables for the wider filter ratios $\Delta / h = 3.5$
        and $5.0$ are given in \cref{sec:tgv-tables}.
    }
    \label{tab:errors-tgv}
    \begin{tabular}{l l l l l l l}
    \toprule
    Model
                & Tier
                & Closure \eqref{eq:tensor-error-prior}
                & Cross-corr.
                & Solution \eqref{eq:tensor-error-post}
                & Median diss.
                & Backscatter \\
    \midrule
    No-model    & --    & $1.0000$            & N.A.                & $0.5985$            & $0.000$           & $0.0000$ \\
    Dyn. Smag.  & --    & $0.9716$            & --                  & $0.5277$            & $2.105$           & $0.0056$ \\
    Clark       & --    & $0.4991$            & $0.8656$            & $0.5219$            & $0.644$           & $0.3105$ \\
    MLP         & p3000 & $0.4994 \pm 0.0005$ & $0.8586 \pm 0.0002$ & $0.5292 \pm 0.0075$ & $2.006 \pm 0.014$ & $0.1448 \pm 0.0012$ \\
    MLP+Re      & p3000 & $0.4933 \pm 0.0012$ & $0.8628 \pm 0.0007$ & $0.5531 \pm 0.0047$ & $1.744 \pm 0.007$ & $0.1601 \pm 0.0018$ \\
    G-CNN       & p3000 & $0.4991 \pm 0.0002$ & $0.8586 \pm 0.0001$ & $0.5262 \pm 0.0015$ & $2.073 \pm 0.007$ & $0.1425 \pm 0.0005$ \\
    G-CNN+Re    & p3000 & $0.4925 \pm 0.0009$ & $0.8629 \pm 0.0003$ & $0.5451 \pm 0.0016$ & $1.796 \pm 0.025$ & $0.1584 \pm 0.0014$ \\
    TBNN        & p3000 & $0.4935 \pm 0.0002$ & $0.8618 \pm 0.0001$ & $0.4933 \pm 0.0049$ & $1.987 \pm 0.009$ & $0.1375 \pm 0.0008$ \\
    TBNN+Re     & p3000 & $0.4905 \pm 0.0006$ & $0.8658 \pm 0.0005$ & $0.5224 \pm 0.0073$ & $1.780 \pm 0.006$ & $0.1495 \pm 0.0015$ \\
    \bottomrule
\end{tabular}

\end{table*}

\paragraph{Dissipation benchmark}
The left and centre panels of \cref{fig:dissipation-tgv} show the canonical
Taylor--Green benchmark---the resolved kinetic energy and total dissipation
rate versus time.
All closures track the laminar roll-up and the onset of transition closely;
differences emerge near peak dissipation ($t^* \approx 9$) and grow through the
decay.
The data-driven closures dissipate \emph{too early} and slightly under-predict
the resolved dissipation peak, whereas Clark overshoots the peak and over-decays
the late-time energy.
As in the forced case the equivariant and non-equivariant networks stay clustered
together, and the $+\mathrm{Re}$ variants depart from their Reynolds-blind
counterparts only in the decay tail---the architecture does not separate the
closures here either.

\paragraph{Structure versus magnitude}
The a priori statistics, computed on identical filtered DNS fields for every
model, locate the discrepancy in the dissipation \emph{magnitude} rather than the
tensor \emph{structure}.
The data-driven models keep a relative tensor error of $0.49$--$0.50$ and a
cross-correlation $\approx 0.86$ with the reference (\cref{tab:errors-tgv}),
about ten percent above their forced values (\cref{tab:errors}) and barely moved
by the $+\mathrm{Re}$
input: they get the shape of the stress largely right, and that shape is close
to Reynolds- and flow-type-independent.
The left panel of \cref{fig:tgv-redelta} makes this concrete: the
Taylor--Green points lie close to the forced $\operatorname{Re}_\Delta$
curve, just above it.
As in the forced case, the models also reproduce only about half the reference
backscatter fraction ($\approx 0.14$--$0.16$ versus $0.30$, \cref{tab:errors-tgv}); this
deficit is intrinsic to the models rather than a generalization artifact.

The dissipation magnitude is a different story.
At the matched operating point ($\nu = 2.5\e{-4}$, $\Delta / h = 2.5$) the
closures over-dissipate by about $2.0$--$2.1\times$ on the Taylor--Green flow,
against only about $1.3\times$ on the forced flow at the \emph{same} $\nu$ and
$\Delta$ (\cref{tab:errors,tab:errors-tgv}).
On the right panel of \cref{fig:tgv-redelta} the Taylor--Green points
therefore sit \emph{above} the forced dissipation curve at the same
$\operatorname{Re}_\Delta$.
This vertical gap is a genuine flow-type effect: less of the energy lies
below the filter in this milder, decaying, under-developed flow than at the
same $\operatorname{Re}_\Delta$ in stationary forced turbulence, and
$\operatorname{Re}_\Delta$, a single scalar built from the mean resolved
gradient magnitude, cannot see the difference.
Supplying $\operatorname{Re}_\Delta$ nonetheless removes part of the excess: at
the tightest filter, where the over-dissipation is largest, it pulls the G-CNN and
TBNN from $2.07$ and $1.99$ down toward $1.80$ and $1.78$ (\cref{tab:errors-tgv}).
The correction is only partial, and it is not uniform across filter width: at the
widest, extrapolated filter ($\Delta / h = 5.0$), where the over-dissipation is
already milder ($\approx 1.4\times$) and $\operatorname{Re}_\Delta$ sits at the top
of the swept range, the forced-calibrated mapping over-corrects and the
$+\mathrm{Re}$ dissipation rises instead ($1.38 \to 1.61$ and $1.37 \to 1.59$ for
G-CNN and TBNN).
The $\operatorname{Re}_\Delta$ input thus addresses the Reynolds-calibration
share of the bias (clearly so where that share dominates) but not the
flow-type share; on a genuinely different flow the two need not even share a
sign.

\paragraph{Solution error}
The right panel of \cref{fig:dissipation-tgv} shows the a posteriori solution
error in time.
The data-driven models hold the lowest error throughout the dynamically
active transition ($t^* \lesssim 15$), where they again stay clustered.
The best closures attain the smallest time-mean error (TBNN $0.49$, G-CNN
$0.53$), comparable to Clark ($0.52$) and dynamic Smagorinsky ($0.53$) and
well ahead of No-model ($0.60$).
Clark's error rises earlier but plateaus in the long decaying tail, narrowing its
time-mean gap; and at the two wider filter ratios its rollout diverges outright,
leaving the data-driven closures the only stable structurally-faithful option.
The $+\mathrm{Re}$ effect on this trajectory error is small and does not track the
dissipation correction: at $\Delta / h = 2.5$ the $+\mathrm{Re}$ closures dissipate
\emph{less} yet end marginally \emph{worse}, while at $\Delta / h = 5.0$ they
dissipate more yet end slightly better (\cref{tab:errors-tgv}).
We therefore read the unambiguous $\operatorname{Re}_\Delta$ generalization result
as the forced $(\nu, \Delta)$ grid of \cref{sec:results-reynolds}; on the
Taylor--Green flow the input helps the dissipation calibration where it is most
off, but that gain does not translate cleanly into a lower solution error.

\paragraph{Interpretation}
The over-dissipation follows directly from the closures' construction.
Because the Reynolds-blind models depend only on the normalized gradient
$\bar{A} / |\bar{A}|$ and not on $\operatorname{Re}_\Delta$, their learned
normalized coefficient is pinned to the forced training regime; in a less
turbulent flow, where less energy resides below the filter, the true normalized
stress is smaller and a fixed coefficient over-dissipates.
It is also consistent with the structural training objective: fit to the stress
tensor itself rather than to its dissipation $\bar{S}_{ij} \tau_{ij}$
(\cref{sec:symmetries}), the networks reproduce the SFS \emph{structure}
faithfully but are never directly optimized for the dissipation \emph{magnitude},
leaving it free to drift when the regime changes.
The filter-scale Reynolds number targets exactly the Reynolds-driven part of
this drift, and on the Taylor--Green flow it corrects that part while leaving
the flow-type residual.
This is the clearest statement of both the reach and the limit of the
$\operatorname{Re}_\Delta$ input.
In summary, the data-driven closures generalize only \emph{partially} to this
out-of-distribution flow: they remain stable, structurally faithful, and clearly
superior to the no-model and purely dissipative baselines through the dynamically
active transition, but carry a flow-type over-dissipation that
$\operatorname{Re}_\Delta$ mitigates only in part.
The parameter-free Clark model carries the opposite, under-dissipative bias and,
supplying no net dissipation of the resolved scales, becomes fragile as the filter
widens.

\subsection{Discussion: the common error floor as an optimal closure}
\label{sec:discussion}

The near-equivalence of the three data-driven closures
(\cref{sec:results-forced,sec:results-capacity}) is striking given how
differently the three networks are built: an unconstrained MLP, a G-CNN with
weights tied across the octahedral group, and a TBNN that expands the stress
in a fixed tensor basis.
We read it as evidence that, in this regime, all three approximately learn
the \emph{same} pointwise mapping from the local filtered velocity-gradient
tensor to the sub-filter stress, $\bar{A}(x) \mapsto \tau(x)$: the minimizer
of the structural training objective, i.e.\ the best predictor of $\tau$ from
$\bar{A}$ in the relative-tensor metric.
The three architectures differ in their inductive bias and in how they
parametrize this map, yet they approximate nearly the same function by every
aggregate diagnostic we measured: a priori tensor errors, a posteriori
trajectories, energy spectra, and the distribution of sub-filter dissipation.
We infer this agreement from these aggregate metrics rather than from a
direct pointwise comparison of the model outputs in function space, which we
did not perform.
The common limit itself, however, is not inferred but measured: the
conditional-mean estimate of \cref{sec:results-capacity} is the function all
three networks saturate onto.
That even the unconstrained MLP reaches this optimum, which is exactly
equivariant as an isotropic average, shows that the symmetry constraints
cost no accuracy here.

This notion of an optimal map carries two qualifications.
First, it is \emph{structural}: the loss regresses the stress tensor
(\cref{sec:symmetries}), so the learned map reproduces the SFS
\emph{structure} but is never asked to match the \emph{functional}
dissipation $\bar{S}_{ij} \tau_{ij}$.
This is why all three share the same backscatter deficit and
over-dissipation, mild in the forced regime and pronounced on the
out-of-distribution Taylor--Green flow (\cref{sec:results-tgv}); a functional
objective would define a \emph{different} optimum.
Second, the map is optimal only relative to its input, the gradient at a
\emph{single} point.
The minimizer of a mean-square loss is the conditional mean
$\mathbb{E}[\tau(u)(x) \mid \bar{A}(x)]$, the average true stress over all
flows sharing the same local filtered gradient.
This is the one-point instance of the \emph{ideal} closure of Langford and
Moser~\cite{langfordOptimalFormulationsIsotropic1999}: the conditional
average of $\tau$ over all turbulent fields whose resolved scales match the
LES field.
Because filtering is many-to-one, $\bar{A}(x)$ does not determine $\tau(x)$;
the conditional average keeps the predictable share of the stress and
averages away the fluctuating remainder, and with it part of the
backscatter, which stochastic closures reintroduce as an explicit random
term~\cite{masonStochasticBackscatterLargeeddy1992}.
The a priori tensor error therefore saturates well above zero for \emph{any}
pointwise closure on this input, not only the particular networks we use.
The direct estimate of \cref{sec:results-capacity} puts this irreducible
residual at $43.9\%$ in distribution, and the saturated networks land within
$0.6\%$ (TBNN) to $2.3\%$ (G-CNN) of it
(\cref{tab:errors}, \cref{fig:condmean}).

Langford and Moser report the analogue in isotropic turbulence: every
optimal estimate they compute leaves a residual nearly as large as the
subgrid force itself, and they argue that, if the ideal model shares this
error, it explains why subgrid closures can perform well a posteriori
despite poor a priori correlations---consistent with the low solution errors
we obtain here.
Whether the ideal model does share it remained a hypothesis for them: the
measured error of an optimal estimate is the fundamental residual of the
ideal model plus the error of the stochastic approximation to the
conditional average, and they saw no way to separate the two.
The saturation observed here separates them for the one-point input: the
direct estimate measures the fundamental residual itself, and the networks
landing within a few percent of it show that the approximation error is
small.
The magnitudes differ for a reason Langford and Moser themselves give: their
near-total residual concerns the subgrid \emph{force} under a sharp spectral
cutoff, and they note that other filters leave a larger deterministic
component in the subgrid term, consistent with the lower floor we measure
for the stress under the Gaussian filter used here.
Their ideal-closure program continued into optimal LES models estimated from
multi-point statistics~\cite{volkerOptimalLargeEddy2002,
langfordOptimalLargeeddySimulation2004,
zandonadeFinitevolumeOptimalLargeeddy2004}; Moser et
al.~\cite{moserStatisticalPropertiesSubgridScale2021} review the statistical
view of closure modeling it established, of which the saturation measured
here is a one-point, learned instance.

This saturation has precedents on both the data and the model side.
Prakash et al.~\cite{prakashInvariantDatadrivenSubgrid2022} report the data
side for an invariant one-point closure whose inputs form a complete set of
nondimensionalized gradient invariants: a single-hidden-layer network
trained on one snapshot at one filter width saturates, and additional
samples, time steps, and filter widths bring no further a priori gain.
They attribute this to the invariances embedded in their model form; the
optimal-closure reading sharpens it: once the inputs exhaust the one-point
gradient, the conditional mean is determined, and neither capacity nor data
volume can lower the floor.
On the model side, the most comprehensive estimate of Langford and Moser, a
46-term expansion subsuming the Smagorinsky, scale-similarity, and
structure-function forms, improved on their simple linear estimate by only a
few percent.
Our measurement makes the point most starkly: seven \emph{constant}
tensor-basis coefficients already come within half a percent of the optimum
(\cref{sec:results-capacity}), so on this input there is very little
function left to learn.
Conditioning on more of the resolved field, such as $\bar{A}$ in a
neighborhood of $x$, defines a different and sharper optimal estimate---a
direction we return to in the outlook of \cref{sec:conclusion}.

\section{Conclusion} \label{sec:conclusion}

Symmetries are fundamental to turbulence, and LES closures are expected to
preserve them~\cite{spalartOldFashionedFrameworkMachine2023}.
We compared three data-driven closures alongside the classical dynamic
Smagorinsky and Clark models, after first determining which continuous
Navier--Stokes symmetries the discrete, forced pseudo-spectral system
actually retains (\cref{sec:symmetries}).
The three closures share a pointwise, Galilean-invariant velocity-gradient
construction but treat rotational and reflectional symmetry differently: an
unconstrained MLP, an exactly octahedral-equivariant G-CNN, and a
tensor-basis TBNN.

The central result is a measured accuracy floor.
Swept across network sizes, the three architectures saturate to the same
a priori and a posteriori accuracy, and a direct, training-free estimate of
the conditional mean by binned least squares over the gradient invariants
(\cref{sec:results-capacity}) reproduces this floor insensitively to the bin
resolution: the $\approx 44\%$ a priori floor is the one-point optimal
closure of Langford and
Moser~\cite{langfordOptimalFormulationsIsotropic1999}, a property of the
one-point input rather than of any particular architecture
(\cref{sec:discussion}).
Despite not learning exact equivariance (an a priori equivariance error of
about $7\%$), the MLP matches the equivariant models on every accuracy
metric.
Its octahedral average, the symmetrized MLP, is exactly equivariant by
construction and matches the raw MLP throughout, confirming that the learned
function already lies close to the equivariant subspace.
We attribute this parity to the learning problem: a low-dimensional,
pointwise map trained on abundant, statistically isotropic data, which lets
even an unconstrained network find the near-equivariant optimum directly, in
line with the data-versus-built-in-symmetry trade-off discussed by McConkey
et al.~\cite{mcconkeyTurbulenceTeachesEquivariance2026a}.
The inductive bias is not idle, though: the equivariant and tensor-basis
closures sit at the floor already at the smallest tested size, which the MLP
reaches only with $25$ times more parameters and with an order-of-magnitude
wider seed-to-seed spread below saturation.
The bias thus buys parameter (and, by extension,
data~\cite{guanLearningPhysicsconstrainedSubgridscale2023}) efficiency and
reproducibility rather than a lower error floor, the regime equivariance
theory predicts when the target itself is equivariant: a variance
reduction, not a shift of the attainable
optimum~\cite{elesedyProvablyStrictGeneralisation2021}.
Of the two equivariant routes, the TBNN also enforces equivariance at
essentially no runtime cost, whereas the G-CNN's $48$-channel regular
representation makes it up to three and a half times slower than the MLP at
inference, the ratio growing with network size.

The second result concerns Reynolds numbers.
A Reynolds-blind closure regresses the normalized stress
$\tau / (\Delta^2 |\bar{A}|^2)$ with no dependence on the viscosity, so its
dissipation calibration is pinned to the training regime, while the tensor
structure it learns is essentially Reynolds-independent
(\cref{sec:results-reynolds}).
Supplying the scaling-invariant filter-scale Reynolds number
$\operatorname{Re}_\Delta$ \eqref{eq:redelta} as an input, in the spirit of
Buaria and Sreenivasan~\cite{buariaForecastingSmallscaleDynamics2023} for
velocity-gradient dynamics, and training across several viscosities removes
this pinning: the $+\mathrm{Re}$ closures hold their dissipation calibration
on held-out viscosities and filter ratios, with a residual overshoot at the
largest $\operatorname{Re}_\Delta$.
On the decaying Taylor--Green vortex, an out-of-distribution flow milder than
the training data, the closures remain stable and structurally faithful but
over-dissipate; $\operatorname{Re}_\Delta$ corrects the Reynolds-driven share
of this excess and leaves the share owed to the different sub-filter content
of a different flow type (\cref{sec:results-tgv}).
Along both axes, capacity and Reynolds number, the equivariant and
non-equivariant architectures do not separate: the generalization gain comes
from the input feature, available to all three, not from the architecture.

These conclusions carry the following limitations.
The Reynolds range spans only a factor of about two in
$\operatorname{Re}_\lambda$, its upper end bounded by the DNS resolution the
smallest viscosity demands, so the across-Reynolds trends are evidence within
a moderate range rather than an asymptotic statement.
The LES resolution is fixed at $M = 128$ throughout: the filter ratio
$\Delta / h$ is swept, but generalization across grid resolution, one
motivation for the $\Delta^2 |\bar{A}|^2$
nondimensionalization~\cite{perezhoginGeneralizableNeuralNetworkParameterization2025},
is not tested here.
All data are generated with a single (Gaussian) filter shape, and both the
measured error floor and the calibration results could shift for sharper
filters such as a spectral cutoff.
The forced a posteriori rollouts span about one integral time, so the
stability statements concern this window; the Taylor--Green rollouts, spanning
the full transition and decay, probe longer horizons.
Finally, the study considers only isotropic turbulence on a periodic box.
In anisotropic or wall-bounded flows, where anisotropy can persist down to
the LES cutoff despite Kolmogorov's hypothesis of local
isotropy~\cite{frischTurbulenceLegacyKolmogorov1995}, the sub-filter stress
need not be equivariant under arbitrary rotations, and a closure that
enforces full rotational equivariance may be over-constrained; relaxed group
convolutions that accommodate partial symmetry
breaking~\cite{wangDiscoveringSymmetryBreaking2025} provide a potential
framework for extending our approach there.

Two directions follow naturally.
Because the one-point input bounds the attainable accuracy, conditioning on
more of the resolved field, through larger convolutional stencils or
graph-based message passing, would raise the accuracy ceiling itself rather
than approach it more efficiently.
And since the structural a priori loss never constrains the dissipation
magnitude, combining the symmetry-preserving architectures with a posteriori
training~\cite{sirignanoDPMDeepLearning2020,
listLearnedTurbulenceModelling2022} could target precisely the
calibration errors that dominate out of distribution.

\section*{Software and reproducibility statement}

The code used to generate the results of this paper is archived on Zenodo at
\url{https://doi.org/10.5281/zenodo.21281963}.
It is released under the MIT license.
The archive includes the exact software environment (pinned package
versions), the random seeds, and the configuration of every experiment, so
that all simulations, trained closures, and figures can be regenerated.

\section*{Data availability}

No data sets are archived with this article.
All data used in this work are generated by the archived code (see the
software and reproducibility statement) from the provided experiment setups
and fixed random seeds, and can be reproduced with it.

\section*{CRediT author statement}

\textbf{Syver Døving Agdestein}:
Conceptualization,
Data curation,
Formal analysis,
Investigation,
Methodology,
Software,
Validation,
Visualization,
Writing -- Original Draft,
Writing -- Review \& Editing

\textbf{Benjamin Sanderse}:
Funding acquisition,
Project administration,
Supervision,
Writing -- Review \& Editing

\section*{Declaration of generative AI and AI-assisted technologies in the manuscript preparation process}

During the preparation of this work the authors used Anthropic Claude Code
to improve the language and readability of the manuscript, to review draft
versions of the text, and to assist in writing and reviewing the research
software and the data-analysis scripts. After using this tool, the authors
reviewed and edited the content as needed, verified the generated code and
analyses, and take full responsibility for the content of the publication.

\section*{Declaration of competing interest}

The authors declare that they have no known competing financial interests or
personal relationships that could have appeared to influence the work reported
in this paper.

\section*{Acknowledgements}

This work is supported by the project ``Discretize first, reduce next'' (with
project number VI.Vidi.193.105) of the research programme NWO Talent Programme
Vidi  financed by the Dutch Research Council (NWO).
This work used the Dutch national e-infrastructure with the support of the
SURF Cooperative using grant no. EINF-15798.

\appendix

\section{Parametrization of the 3D orthogonal group on Cartesian grids}
\label{sec:group-parametrization}

For vectors $x$ and tensors $\sigma$, we employ two indexing notations:
\begin{equation}
    x_i = x[i], \quad \sigma_{i j} = \sigma[i, j].
\end{equation}
Square bracket notation is used when indices involve expressions more complex
than a single symbol.
This convention also applies to the Kronecker delta $\delta_{i j} = \delta[i, j]$, which
equals $1$ if $i = j$ and $0$ otherwise.

The orthogonal group $O(3)$ comprises all rotations and reflections in 3D space.
On uniform Cartesian grids, we restrict attention to a
grid-compatible subgroup of $O(3)$, where rotations are limited to multiples of $\pi / 2$.
This restriction ensures that transformed grids align with the original grid.
We denote this subgroup as $G \subset O(3)$.
This is the \emph{octahedral group} (or ``cube'' group),
representing the symmetries of an octahedron (or equivalently, a cube).

$G$ contains $| G | = 48$ elements:
$24$ are orientation-preserving rotations
(elements of the special orthogonal group $SO(3)$),
while the remaining $24$ include a reflection.

In its physical-space representation, every element $g \in G$ maps the set of
signed coordinate axes $\{ \pm e_1, \pm e_2, \pm e_3 \}$ onto itself, so its
matrix is a \emph{signed permutation matrix}: every row and column has exactly
one non-zero entry, equal to $+1$ or $-1$.
This means that each element $g \in G$ corresponds to a unique triple
$(R, p, s)$: a rotation matrix $R \in \mathbb{R}^{3 \times 3}$, a permutation
vector $p \in \mathbb{R}^3$, and a sign vector $s \in \mathbb{R}^3$.
Using the spatial indices $(i, j) \in \{ 1, 2, 3 \}^2$, these are related through
(with no sum over $i$)
\begin{equation} \label{eq:signperm}
    R_{i j} \coloneq s_i \, \delta[p_i, j].
\end{equation}
The subscripts pick out single components of the $3$-vectors $p$ and $s$:
for each row $i$, the component $p_i \in \{ 1, 2, 3 \}$ points to the column of the
single non-zero entry, and $s_i \in \{ +1, -1 \}$ is its sign.

There are $6$ possible permutation vectors $p$, given by
\begin{equation} \label{eq:permutations}
    p \in \left\{
        \begin{pmatrix} 1 \\ 2 \\ 3\end{pmatrix},
        \begin{pmatrix} 2 \\ 3 \\ 1\end{pmatrix},
        \begin{pmatrix} 3 \\ 1 \\ 2\end{pmatrix},
        \begin{pmatrix} 3 \\ 2 \\ 1\end{pmatrix},
        \begin{pmatrix} 2 \\ 1 \\ 3\end{pmatrix},
        \begin{pmatrix} 1 \\ 3 \\ 2\end{pmatrix} \right\},
\end{equation}
and $8$ possible sign vectors $s$, given by
\begin{equation} \label{eq:signflips}
    s \in \left\{
        \begin{pmatrix} + \\ + \\ + \end{pmatrix},
        \begin{pmatrix} - \\ + \\ + \end{pmatrix},
        \begin{pmatrix} + \\ - \\ + \end{pmatrix},
        \begin{pmatrix} + \\ + \\ - \end{pmatrix},
        \begin{pmatrix} - \\ - \\ + \end{pmatrix},
        \begin{pmatrix} + \\ - \\ - \end{pmatrix},
        \begin{pmatrix} - \\ + \\ - \end{pmatrix},
        \begin{pmatrix} - \\ - \\ - \end{pmatrix}
    \right\},
\end{equation}
where $+$ and $-$ mean $+1$ and $-1$, respectively.

For a vector $x = (x_1, x_2, x_3) \in \mathbb{R}^3$, we get the transformed vector
(with no sum over $i$)
\begin{equation}
    \left(R x\right)_i = s_i x[p_i].
\end{equation}
For a tensor $\sigma \in \mathbb{R}^{3 \times 3}$, we get the transformed tensor
(with no sum over $i$ and $j$)
\begin{equation}
    \left(R \sigma R^T\right)_{i j} = s_i s_j \sigma[p_i, p_j].
\end{equation}

We encode the elements of $G$ as $g_{a b}$,
where $a \in \{ 1, \dots, 6 \}$ denotes the permutation number
in \eqref{eq:permutations} and
$b \in \{ 1, \dots, 8 \}$ denotes the sign number in \eqref{eq:signflips}.
With a ``column-major'' ordering,
the group elements are uniquely identified by a linear index
$c \coloneq a + 6 (b - 1)$, and we define the ``flat ordering''
$g_c \coloneq g_{a b}$.
Thus $G = \{ g_c \ | \ c \in \{ 1, \dots, 48 \} \}$.

\begin{definition}[Regular representation]
    The regular representation of $G$ assigns to each element $g \in G$ a matrix
    $P_g \in \{0, 1\}^{|G| \times |G|}$ with entries
    \begin{equation}
        (P_g)_{i j} = \delta[g_i, g g_j],
    \end{equation}
    where $(g_i)_{i = 1}^{|G|}$ are the group elements of $G$ ordered with the
    flat ordering and $g g_j$ is a composition of the elements $g$ and $g_j$.
\end{definition}

In other words, $P_g \in \mathbb{R}^{48 \times 48}$
is a permutation matrix that maps all the group elements to
the elements they become after applying $g$.
For an element $g_j \in G$, the composed roto-reflection $g g_j$ is equal to
the unique roto-reflection $g_i$ such that $(P_g)_{i j} = 1$.

\begin{definition}[Cayley table]
    The Cayley table of $G$ is a matrix $C \in \mathbb{N}^{|G| \times |G|}$ with entries defined such that
    \begin{equation}
        g[C_{i j}] = g_i g_j.
    \end{equation}
\end{definition}

In other words, $g_k = g_i g_j$ if and only if $C_{i j} = k$.

\section{Equivariance of the group-convolutional layers}
\label{sec:weight-sharing}

This appendix collects the equivariance computations for the three layer types
of the G-CNN of \cref{sec:groupconv}: the inner group-convolutional layers
acting between regular representations, the initial layer that lifts the
flattened tensor input into the regular representation, and the final layer
that maps back to a flattened tensor.

\paragraph{Inner layers}
Writing the linear part of the group-convolutional layer
\eqref{eq:gconv-layer} as $L$, so that
$(L_i \xi)(g) = \sum_j \sum_h k_{i j}(h^{-1} g) \, \xi_j(h)$, a change of variables
$h \leftarrow a h$ in the group sum gives, for every $a \in G$,
\begin{equation}
    \begin{split}
        \bigl(L_i \, P_a \xi\bigr)(g)
        & = \sum_{j} \sum_{h \in G} k_{i j}(h^{-1} g) \, \xi_j(a^{-1} h) \\
        & = \sum_{j} \sum_{h \in G} k_{i j}(h^{-1} a^{-1} g) \, \xi_j(h) \\
        & = (L_i \xi)(a^{-1} g)
        = \bigl(P_a \, L_i \xi\bigr)(g),
    \end{split}
\end{equation}
using the regular action $(P_a \xi)(h) = \xi(a^{-1} h)$ of \cref{eq:regular-action}
and the bijectivity of $h \mapsto a h$ on $G$.
The bias is constant across $g$ and the activation $\varphi$ acts element-wise, so
both commute with $P_a$: since only one entry in each row of $P_a$ is non-zero,
applying $\varphi$ before or after the permutation yields the same
result~\cite{cohenSteerableCNNs2016}.
Hence $P_a \, \ell(\xi) = \ell(P_a \xi)$ for all $a \in G$.

\paragraph{Initial layer}
The lifting layer of \cref{sec:groupconv} reads
$\zeta_i(g) = \varphi\bigl( \sum_\mu (Q_g c_i)(\mu) \, a(\mu) + b_i \bigr)$.
Using the orthogonality $Q_g^T = Q_{g^{-1}}$ and the representation property
$Q_{g^{-1}} Q_h = Q_{g^{-1} h}$,
\begin{equation}
    \begin{split}
        \ell_i(Q_g a)(h)
        & = \varphi\Bigl( \sum_\mu (Q_h c_i)(\mu) \, (Q_g a)(\mu) + b_i \Bigr) \\
        & = \varphi\Bigl( \sum_\mu (Q_{g^{-1} h} c_i)(\mu) \, a(\mu) + b_i \Bigr) \\
        & = \zeta_i(g^{-1} h) = (P_g \zeta_i)(h),
    \end{split}
\end{equation}
so the orbit of $c_i$ intertwines the input action $Q_g$ with the output
permutation $P_g$: $\ell(g a) = g \, \ell(a)$ for all $g \in G$.

\paragraph{Final layer}
In flattened form the output layer reads
$m = \sum_j \sum_{h \in G} \xi_j(h) \, Q_h d_j$.
Replacing $\xi$ by $P_g \xi$ and substituting $h \leftarrow g h$ in the group sum,
\begin{equation}
    \sum_j \sum_{h \in G} \xi_j(g^{-1} h) \, Q_h d_j
    = \sum_j \sum_{h \in G} \xi_j(h) \, Q_{g h} d_j
    = Q_g m,
\end{equation}
using $Q_{g h} = Q_g Q_h$; this is exactly the tensor action of $g$ on the
flattened output.
Composed, the three layer types give
$g \, m^\text{G-CNN}(\bar{u}) = m^\text{G-CNN}(g \bar{u})$ for all $g \in G$.

\section{Pseudo-spectral discretization of the incompressible Navier--Stokes equations}
\label{sec:pseudo-spectral}

\subsection{Continuous spectral form}

On a periodic box $\Omega \coloneq [0, L]^D$, we expand the velocity in
a Fourier series with integer wavenumber $k \in \mathbb{Z}^D$,
\begin{equation}
    u(x, t) = \sum_{k \in \mathbb{Z}^D} \hat{u}(k, t) \, \exp(\mathrm{i} \kappa_j x_j),
    \quad \kappa_j \coloneq \frac{2 \pi k_j}{L},
\end{equation}
where $\mathrm{i}$ is the imaginary unit and summation over repeated indices
is implied.
Spatial derivatives become wavenumber multiplications,
$\widehat{\partial_j u} = \mathrm{i} \kappa_j \hat{u}$, and the incompressible
Navier--Stokes equations take the spectral form
\begin{equation} \label{eq:spectral-ns}
    \kappa_j \hat{u}_j = 0, \quad
    \partial_t \hat{u}_i + \mathrm{i} \kappa_j \hat{\sigma}_{i j}(u)
    = - \mathrm{i} \kappa_i \hat{p},
\end{equation}
where
\begin{equation}
    \hat{\sigma}_{i j}(u) \coloneq \widehat{\check{u}_i \check{u}_j}
    - \nu \mathrm{i} (\kappa_j \hat{u}_i + \kappa_i \hat{u}_j)
\end{equation}
is the spectral stress tensor,
$\hat{p}$ is the spectral pressure,
$\nu > 0$ is the kinematic viscosity, and
$\widehat{(\cdot)}$ and $\widecheck{(\cdot)}$ denote the forward and inverse
Fourier transforms, respectively.

Incompressibility is enforced by the Leray projector
\begin{equation}
    \Pi_{i j}(k) \coloneq \delta_{i j} - \frac{\kappa_i \kappa_j}{\kappa_\ell \kappa_\ell},
    \qquad \Pi_{i j}(0) \coloneq \delta_{i j},
\end{equation}
which eliminates the pressure from \cref{eq:spectral-ns} and gives the
pressure-free momentum equation
\begin{equation} \label{eq:spectral-projected}
    \partial_t \hat{u}_i = - \Pi_{i j}(k) \mathrm{i} \kappa_\ell \hat{\sigma}_{j \ell}(u).
\end{equation}

\subsection{Discretization}

We represent the velocity field $u(x, t)$ on a Cartesian grid with $N^D$
collocated points and even $N$, with spacing $h \coloneq L / N$.
The Fourier series is replaced by the discrete Fourier transform (DFT)
\begin{align}
    \mathcal{F}_N[u](k) & \coloneq \frac{1}{N^D} \sum_{x} u(x) \exp(- \mathrm{i} \kappa_j x_j), \\
    \mathcal{F}_N^{-1}[\hat{u}](x) & \coloneq \sum_{k} \hat{u}(k) \exp(\mathrm{i} \kappa_j x_j),
\end{align}
where the sums run over the $N^D$ collocation points and over the
resolved integer wavenumbers $k \in \{ - N/2, \dots, N/2 - 1 \}^D$.
The factor $1 / N^D$ is placed on the forward transform so that
Parseval's identity reads
$\sum_{x} u(x)^2 / N^D = \sum_{k} | \hat{u}(k) |^2$ directly.
Since $u$ is real-valued in physical space, the spectrum satisfies the
Hermitian symmetry $\hat{u}(- k) = \operatorname{conj}(\hat{u}(k))$, which
the implementation exploits by storing only the half-spectrum
$k_1 \ge 0$ (a real-to-complex FFT).

\subsubsection{Numerical stress and antialiasing}

Forming the nonlinear product $u_i u_j$ on the grid generates frequencies
outside the resolved band and aliases them back onto the resolved modes.
We remove these errors with the $2/3$-rule of
Patterson and Orszag~\cite{pattersonSpectralCalculationsIsotropic1971}.
Let $T_N$ denote the spectral truncation operator
\begin{equation}
    (T_N \hat{u})(k) \coloneq \begin{cases}
        \hat{u}(k), & \text{if } | k |_\infty \le N / 3, \\
        0, & \text{otherwise},
    \end{cases}
\end{equation}
and define the numerical stress tensor by
\begin{equation} \label{eq:sigma-N}
    \sigma^N_{i j}(u) \coloneq T_N \mathcal{F}_N \! \left[
        \mathcal{F}_N^{-1}(T_N \hat{u}_i) \, \mathcal{F}_N^{-1}(T_N \hat{u}_j)
    \right] - \nu \mathrm{i} (\kappa_j \hat{u}_i + \kappa_i \hat{u}_j).
\end{equation}
The outer truncation discards the product harmonics above the cutoff, so the
right-hand side, and hence the solution, stays supported on the retained
modes.
The viscous part is exact at the resolved wavenumbers and needs no
dealiasing.
The discrete (DNS) momentum equation is then
\begin{equation} \label{eq:dns}
    \kappa_j \hat{v}_j = 0, \quad
    \partial_t \hat{v}_i
    = - \Pi_{i j}(k) \mathrm{i} \kappa_\ell \sigma^N_{j \ell}(v),
\end{equation}
where $v$ is the DNS velocity field.

\subsubsection{Time integration}

We advance \cref{eq:dns} in time with Wray's three-stage low-storage
Runge--Kutta scheme~\cite{wray1990minimal}, applying the Leray projector
$\Pi_{i j}$ at the end of every stage so that the velocity remains
discretely divergence-free throughout the step.
The time step $\Delta t$ is selected adaptively by combining the
convective and diffusive stability limits,
\begin{equation}
    \Delta t \coloneq C \, \min \! \left(
        \frac{h}{\max_x | v(x) |},
        \frac{h^2}{2 D \nu}
    \right),
\end{equation}
with safety factor $C < 1$.

\subsection{Large-eddy equations}

In spectral space, a convolutional filter reduces to a wavenumber-wise
scaling.
Let $H(k)$ denote the spectral filter kernel; the filtered field is
$\widehat{\bar{u}}(k) \coloneq H(k) \hat{u}(k)$.
The filtered DNS momentum equation reads
\begin{equation}
    \kappa_j \widehat{\bar{v}}_j = 0, \quad
    \partial_t \widehat{\bar{v}}_i
    = - \Pi_{i j}(k) \mathrm{i} \kappa_\ell \overline{\sigma^N_{j \ell}(v)}.
\end{equation}

In LES we solve on a coarser grid of size $M^D$ with $M < N$.
The resolved stress on the coarse grid is $\sigma^M(\bar{v})$, where
$\sigma^M$ is defined analogously to $\sigma^N$ in \cref{eq:sigma-N}
using $M$ for both DFTs and for the antialiasing
operator~\cite{agdesteinExactExpressionsUnresolved2026}.
Splitting the filtered DNS stress into a discretely resolved part and a
remainder gives
\begin{equation}
    \partial_t \widehat{\bar{v}}_i
    = - \Pi_{i j}(k) \mathrm{i} \kappa_\ell
        \bigl[ \sigma^M_{j \ell}(\bar{v}) + \tau^{N \to M}_{j \ell}(v) \bigr],
\end{equation}
where the \emph{discrete sub-filter stress (SFS)} is
\begin{equation} \label{eq:discrete-sfs}
    \tau^{N \to M}_{i j}(v) \coloneq
    \overline{\sigma^N_{i j}(v)} - \sigma^M_{i j}(\bar{v}).
\end{equation}
Because the spectral derivative commutes with filtering, the viscous
contributions to $\sigma^N$ and $\sigma^M$ cancel in
\cref{eq:discrete-sfs}, so $\tau^{N \to M}$ inherits only the nonlinear
product term, evaluated with the antialiasing operators of the two grids.
The trace of $\tau$ is isotropic and can be absorbed into the LES
pressure; we therefore retain only the deviatoric part
$\operatorname{dev} \tau$, which is also the output space of the closure
models in this paper.

The discrete SFS differs from the continuous physical-space expression
$\check{\tau}(u) \coloneq \overline{\check{u} \, \check{u}}
- \check{\bar{u}} \, \check{\bar{u}}$ because
$\tau^{N \to M}$ depends explicitly on $N$ and $M$ through the DFTs and
the antialiasing operators on both grids.
We numerically verified this discrete SFS expression using a DNS-aided
LES approach~\cite{baeNumericalModelingError2022}.

\section{Data generation} \label{sec:datagen}

We perform both DNS and LES on a single H100 GPU using double precision (64-bit) floating point arithmetic.
The experiment sweeps two physical axes: the kinematic viscosity $\nu$ and the
filter-to-grid ratio $\Delta / h$, where $h \coloneq L / M$ denotes the LES grid
spacing throughout.
The closures are trained on three forced realizations at
$\nu \in \{ 1.5, 2.5, 4.0 \} \e{-4}$ (one random seed each) and tested on two
held-out realizations---a new seed at a training viscosity
($\nu = 2.5 \e{-4}$, \emph{in-distribution}) and a higher Reynolds number at a
viscosity outside the training set ($\nu = 1.0 \e{-4}$, \emph{out-of-distribution}).
The resulting Taylor- and integral-scale Reynolds numbers span
$\operatorname{Re}_\lambda \approx 199$--$394$ and are collected in \cref{tab:dns}.
Combined with the shell-energy rescaling described below, each configuration
produces statistically stationary homogeneous isotropic turbulence suitable
for testing closure models.
The DNS resolution is $N \coloneq 810$ and
the LES resolution is $M \coloneq 128$ for every run.
The $2 / 3$-rule retains modes with $| k |_\infty \le N / 3$, so the
effective resolutions (resolved integer wavenumbers per axis) are
$\lfloor 2 N / 3 \rfloor + 1 = 541$ and
$\lfloor 2 M / 3 \rfloor + 1 = 85$ for DNS and LES, respectively.
The DNS resolution was chosen to maximize grid size while fitting arrays in GPU memory.
While powers of $2$ are more efficient for FFTs,
$N \coloneq 1024$ exceeds available memory.
This single fixed resolution resolves the whole viscosity range: the Kolmogorov
scale stays resolved even at the highest Reynolds number, where
$k_{\max} \eta \approx 1.1$ (\cref{tab:dns}).
Double precision is required for accurate estimation of equivariance errors.

\subsection{Initialization}

Each run uses an independent random seed; its DNS spectral velocity field
$\hat{u}(k)$ is initialized as follows:
\begin{enumerate}
    \item Assign $\hat{u}_i(k) \sim \mathcal{N}(0, 1) + \mathrm{i} \, \mathcal{N}(0, 1)$
        independently at every wavenumber.
    \item Enforce Hermitian symmetry $\hat{u}(- k) = \operatorname{conj}(\hat{u}(k))$,
        so that the corresponding physical-space field is real-valued.
    \item Project onto the divergence-free subspace:
        $\hat{u}(k) \leftarrow \Pi(k) \hat{u}(k)$.
    \item For every shell $s \in \{ 0, 1, \dots, \lfloor \sqrt{D} \, N / 2 \rfloor \}$,
        rescale the coefficients in $K(s)$ to match a prescribed shell energy,
        \begin{equation}
            \hat{u}(k) \leftarrow
                \sqrt{\frac{E_0 P(s)}{E(\hat{u}, s) \sum_{s'} P(s')}} \,
                \hat{u}(k),
            \quad \forall k \in K(s),
        \end{equation}
        where
        \begin{equation} \label{eq:shells}
            K(s) \coloneq \{ k \in \mathbb{Z}^3 \ | \ s \leq | k | < s + 1 \}
        \end{equation}
        is the shell of integer wavenumbers at level $s$,
        \begin{equation} \label{eq:shellenergy}
            E(\hat{u}, s) \coloneq \frac{1}{2} \sum_{k \in K(s)} | \hat{u}(k) |^2
        \end{equation}
        is the energy in shell $s$, $P(s) \coloneq s^{-5/3}$ (with $P(0) \coloneq 0$)
        is a prescribed inertial-range profile with logarithmic slope $-5/3$, and
        $E_0 \coloneq 0.2$ is the prescribed total kinetic energy.
\end{enumerate}
After step 4 the total energy is exactly $E_0$.
Hermitian symmetry is preserved by every subsequent step.

\subsection{Time integration}

We integrate the projected DNS equation~\eqref{eq:dns} with Wray's third-order
low-storage Runge--Kutta scheme~\cite{wray1990minimal} (see
\ref{sec:pseudo-spectral}).
The time step combines convective and diffusive stability limits,
\begin{equation} \label{eq:adaptive}
    \Delta t \coloneq C \min \! \left(
        \frac{h_\text{grid}}{\max_x | u(x) |},
        \frac{h_\text{grid}^2}{2 D \nu}
    \right),
\end{equation}
with safety factor $C \coloneq 0.35$, spatial dimension $D = 3$, and
$h_\text{grid}$ the spacing of the grid being integrated ($L / N$ for the DNS,
$L / M$ for the LES).

\subsection{Shell-energy rescaling}

To sustain the turbulence over time, we maintain a fixed energy in the two
lowest non-trivial shells by rescaling the spectral coefficients at the end
of every Runge--Kutta step.
For $s \in \{ 1, 2 \}$, we set
\begin{equation} \label{eq:rescale}
    \hat{u}(k) \leftarrow
    \sqrt{\frac{E_f(s)}{E(\hat{u}, s)}} \, \hat{u}(k),
    \quad \forall k \in K(s),
\end{equation}
where $E_f(s)$ is the energy in shell $s$ immediately after initialization
(i.e., $E_f(s) = E_0 P(s) / \sum_{s'} P(s')$ from the initialization above).
Lundgren applied a similar linear forcing to \emph{all} wavenumbers to
maintain the total kinetic energy~\cite{lundgrenLinearlyForcedIsotropic2003};
we restrict the rescaling to the two lowest shells to avoid interfering with
the scales requiring closure.
This approach for injecting energy into otherwise-decaying turbulence is widely
used, including in the Johns Hopkins turbulence
database~\cite{caoStatisticsStructuresPressure1999,liPublicTurbulenceDatabase2008}.

\subsection{Filter}

\begin{figure*}
    \centering
    \includegraphics[width=\textwidth]{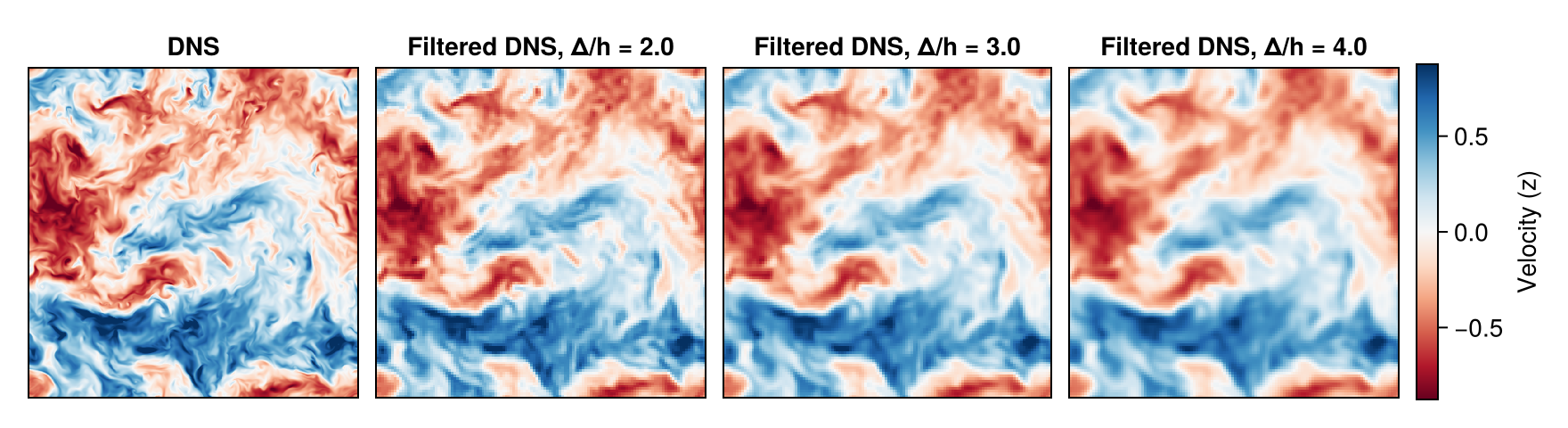}
    \caption{%
        Mid-plane slices ($z = L / 2$) of the velocity component $u_z$ for the
        central training run ($\nu = 2.5 \e{-4}$) at the end of the warm-up
        period. Left: DNS at full resolution ($N = 810$); remaining panels:
        filtered DNS on the LES grid ($M = 128$) at the three training filter
        ratios. All panels share one color range.
    }
    \label{fig:slice-filters}
\end{figure*}

For data generation we use a Gaussian filter.
To ensure the closure model accounts only for the SFS and not for any
commutator between rescaling and filtering, the filter must leave the
rescaled shells $K(1)$ and $K(2)$ unchanged.
We therefore use a modified Gaussian kernel,
\begin{equation} \label{eq:filter}
    H(k) \coloneq \begin{cases}
        1, & \text{if } | k | < 3, \\
        \displaystyle \exp \! \left( - \frac{\Delta^2 (2 \pi / L)^2 | k |^2}{24} \right),
            & \text{if } | k | \geq 3,
    \end{cases}
\end{equation}
where the filter width $\Delta$ is set by the LES grid spacing $h = L / M$
through the filter-to-grid
ratio $\Delta / h$, itself a swept axis of the experiment: training data are
filtered at $\Delta / h \in \{ 2, 3, 4 \}$ and the held-out test data at the
interpolated and extrapolated ratios $\Delta / h \in \{ 2.5, 3.5, 5.0 \}$
(\cref{sec:results-reynolds}).
The carve-out $| k | < 3$ is independent of $\Delta$, so $H$ is unity on the
rescaled shells for every ratio, and the rescaling in~\eqref{eq:rescale}
commutes with the filtering throughout the sweep.
\Cref{fig:slice-filters} shows the effect of the filter on a DNS snapshot:
the large-scale structure is retained at every ratio, while increasing
$\Delta / h$ progressively removes the small scales the closure must account
for.

\subsection{Sampling}

Each forced run is warmed up by integrating the DNS from the random initial
field for $5$ time units, by which point the flow has reached statistically
stationary turbulence (\cref{fig:evolution-data}).
We then continue the DNS and store data pairs $(\widehat{\bar{u}}, \tau(u))$
together with turbulence statistics.
Here $\tau$ denotes the discrete SFS from \cref{eq:discrete-sfs};
for brevity, we omit the superscript $N \to M$.
The save times are set relative to the \emph{measured} integral turnover
$t_\text{int} = L_\text{int} / u'$ of each run, where $L_\text{int}$ is the
integral length scale and $u'$ the root-mean-square velocity
($t_\text{int} \approx 7.2$--$7.6$,
\cref{tab:dns}), so every run is sampled over the same number of turnovers
despite its different Reynolds number.
Training and test data are drawn from \emph{separate realizations} rather than
from a temporal split of a single run:
\begin{itemize}
    \item \emph{Training.} For each of the three training runs and each training
        filter ratio, we store $8$ snapshots spread over $\approx 2\,t_\text{int}$.
        The sparse, wide spacing maximizes a priori diversity, and pooling the
        three viscosities with the three filter ratios yields
        $3 \times 3 \times 8 = 72$ gradient--stress fields spanning the training
        band of the filter-scale Reynolds number.
    \item \emph{Testing.} For each held-out run and each test filter ratio, we
        store a denser series of $40$ snapshots over $\approx 1\,t_\text{int}$,
        used as the reference for the a posteriori rollouts (initialized from the
        first snapshot and integrated across the window).
\end{itemize}
Because the test runs are independent realizations (a fresh random seed at a
training viscosity, and an entirely held-out, higher-Reynolds-number
viscosity), the reported test errors measure genuine generalization rather
than interpolation within a single trajectory.
The learned closures are in any case stateless, pointwise maps on the local
velocity gradient (\cref{sec:results}): they cannot memorize individual snapshots
and must learn the population-level mapping, and each snapshot supplies $128^3$
gradient--stress pairs.
Out-of-distribution behavior beyond the forced regime is assessed separately on
the decaying Taylor--Green vortex (\cref{sec:results-tgv}), whose data are
generated by the same pipeline, without warm-up or forcing, from the analytic
initial condition
\begin{equation} \label{eq:tgv-ic}
    u_0(x) \coloneq V_0 \begin{pmatrix}
        \phantom{-}\sin x_1 \cos x_2 \cos x_3 \\
        -\cos x_1 \sin x_2 \cos x_3 \\
        0
    \end{pmatrix}
\end{equation}
on the same $[0, 2\pi]^3$ box, so the characteristic length of the vortex is
$L_\text{TGV} = 1$ and the nominal Reynolds number is
$\operatorname{Re} = V_0 L_\text{TGV} / \nu$; we set
$V_0 = \operatorname{Re} \nu = 1.5$ to target $\operatorname{Re} = 6000$ at
$\nu = 2.5 \e{-4}$ (initial kinetic energy $V_0^2 / 8$).
The flow is sampled at $100$ instants over $20$ convective times
$t_c = L_\text{TGV} / V_0$.

\subsection{Statistics}

\begin{figure}
    \centering
    \includegraphics[width=0.9\columnwidth]{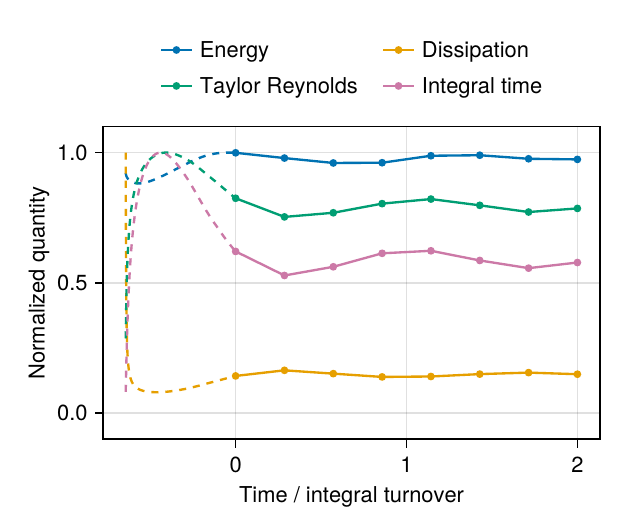}
    \caption{
        DNS time series of the
        total kinetic energy $E$,
        viscous dissipation rate $\epsilon$,
        Taylor scale Reynolds number $\mathrm{Re}_\lambda$, and
        integral time $t_\text{int}$,
        for the central training run ($\nu = 2.5 \e{-4}$); the other realizations
        are statistically equivalent.
        Time is normalized by the measured integral turnover $t_\text{int}$ of
        the run, so the sampling window ends at exactly two turnovers; the
        markers are the $8$ stored snapshots.
        Negative times (dashed lines) indicate the warm-up period.
        All time series are normalized by their respective maximum values.
    }
    \label{fig:evolution-data}
\end{figure}

\begin{figure}
    \centering
    \includegraphics[width=0.9\columnwidth]{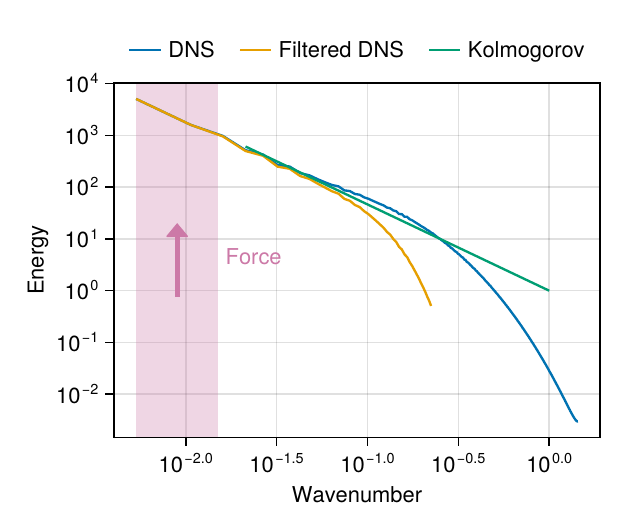}
    \caption{
        Time-averaged energy spectrum after warm-up, for the central training run
        ($\nu = 2.5 \e{-4}$, $\Delta / h = 3$).
        The banded area shows the range of wavenumbers with active
        shell-energy rescaling.
    }
    \label{fig:spectrum-data}
\end{figure}

In \cref{fig:evolution-data} we show the evolution of the total kinetic
energy
\begin{equation}
    E(\hat{u}) \coloneq \frac{1}{2} \sum_{k \in \mathbb{Z}^3} | \hat{u}(k) |^2
\end{equation}
and the viscous dissipation rate
\begin{equation}
    \epsilon \coloneq \nu \left( \frac{2 \pi}{L} \right)^2
    \sum_{k \in \mathbb{Z}^3} | k |^2 | \hat{u}(k) |^2,
\end{equation}
both normalized by their maximum values.
In \cref{fig:spectrum-data} we show the time-averaged shell spectra
$\langle E(\hat{u}, s) \rangle$ and $\langle E(\widehat{\bar{u}}, s) \rangle$
defined by~\eqref{eq:shellenergy}, where $\langle \cdot \rangle$ denotes
time-averaging.
We also show the theoretical Kolmogorov spectrum for the inertial range,
\begin{equation} \label{eq:kolmogorov}
    E_\text{Kol}(q) \coloneq C_\text{K} \langle \epsilon \rangle^{2 / 3} q^{-5 / 3},
\end{equation}
where $q$ is the (dimensional) wavenumber magnitude and
$C_\text{K} \coloneq 1.6$ is the Kolmogorov
constant~\cite{sreenivasanUniversalityKolmogorovConstant1995}.
For shell $s$, the corresponding dimensional wavenumber is
$q_s \coloneq 2 \pi s / L$.
We use the Kolmogorov normalizations
$\tilde{E} \coloneq \epsilon^{-2/3} \eta^{5/3} E$ and
$\tilde{q} \coloneq q \eta$ for plotting, where
$\eta \coloneq (\nu^3 / \epsilon)^{1 / 4}$ is the Kolmogorov length scale.
With these normalizations, the Kolmogorov spectrum equals $1$ at
$\tilde{q} = 1$.
The rescaled wavenumber shells are indicated by a banded area.
After warm-up, the three training runs reach Taylor-microscale Reynolds numbers
$\operatorname{Re}_\lambda \approx 199$--$320$ and measured integral-scale
Reynolds numbers $\operatorname{Re}_L = u' L_\text{int} / \nu$ up to
$\approx 6800$, and the held-out out-of-distribution run reaches
$\operatorname{Re}_\lambda \approx 394$ (\cref{tab:dns}); this is the training
regime that the generalization tests of
\cref{sec:results-reynolds,sec:results-tgv} refer back to.

In \cref{fig:evolution-data}, the dissipation rate peaks at the start of
the warm-up.
This occurs because the simplified initial spectrum assumes an inertial-range
slope across all wavenumbers, including near the Kolmogorov scale where viscous
effects dominate; the initial drop in dissipation reflects the depletion of
this excess high-wavenumber energy.
In \cref{fig:spectrum-data}, the Kolmogorov scales are resolved by
DNS but not by LES.
For the run shown, the rescaled shells span scales approximately $100$ times
larger than the Kolmogorov scale, while the filtered DNS spectrum decays at scales
roughly $10$ times larger than the Kolmogorov scale.

\section{Training} \label{sec:training}

\begin{figure*}
    \centering
    \includegraphics[width=0.9\textwidth]{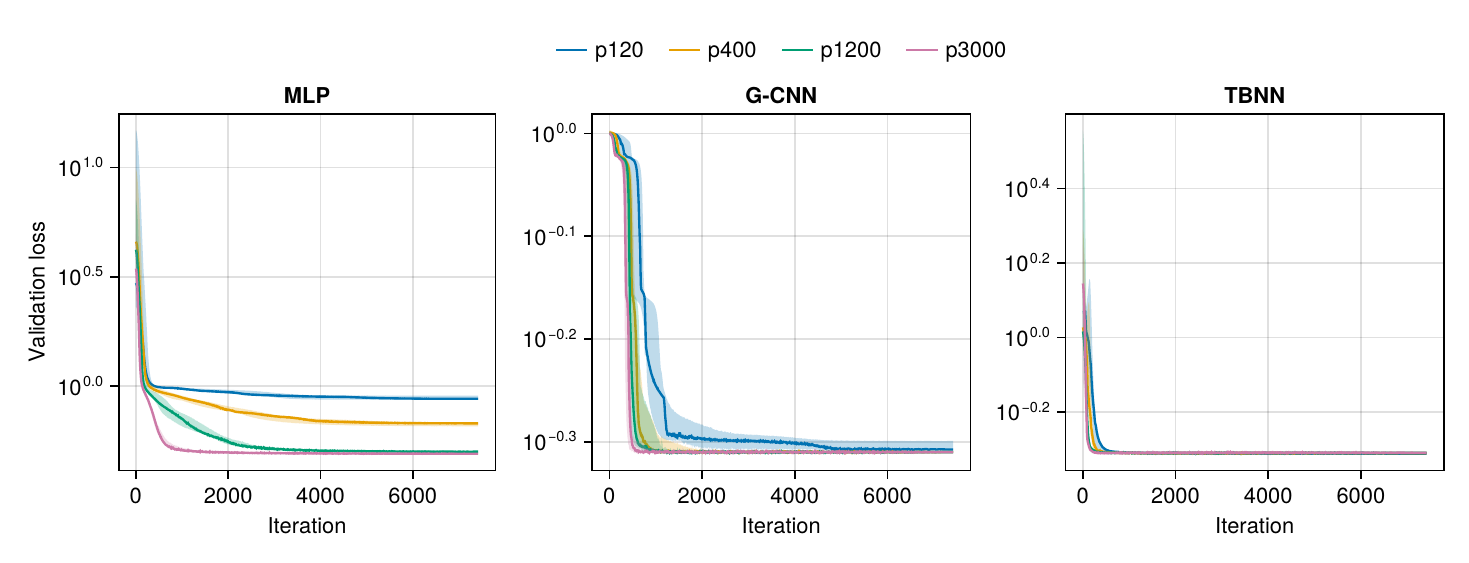}
    \caption{%
        Validation loss during training (log scale), one panel per architecture
        with a curve for each capacity tier (median over the five seeds; bands
        span the seed range). Only the $+\mathrm{Re}$ variants are shown.
    }
    \label{fig:training}
\end{figure*}

\begin{table*}
    \centering
    \caption{
        The capacity grid. Each learned closure is a pointwise
        ($1 \times 1 \times 1$) network (\cref{sec:construction}) with three hidden
        layers, a $\operatorname{GELU}$ activation after each, and a linear,
        bias-free output layer; capacity is varied by hidden width alone. Within a
        tier the three architectures are roughly parameter-matched. ``Hidden''
        lists the three hidden widths; for the G-CNN these
        count regular-representation channels, each expanding to $|G| = 48$ scalar
        feature maps after synthesis (\cref{sec:groupconv}), and its parameter
        count is that of the reduced, pre-synthesis weights updated by the
        optimizer.
        The input/output feature dimensions of the network proper are
        $9 \to 6$ (MLP), $9 \to 9$ (G-CNN), and $5 \to 7$ (TBNN), before the
        symmetrization, trace removal, or basis expansion that maps the output to
        the deviatoric stress. Parameter counts are for the Reynolds-blind
        networks; the $+\mathrm{Re}$ variant appends one input channel and adds
        only a handful of parameters (\cref{sec:redelta-input}).
    }
    \label{tab:architecture}
    \begin{tabular}{l r c r c r c r}
        \toprule
        & & \multicolumn{2}{c}{MLP} & \multicolumn{2}{c}{G-CNN} & \multicolumn{2}{c}{TBNN} \\
        \cmidrule(lr){3-4} \cmidrule(lr){5-6} \cmidrule(lr){7-8}
        Tier & Target & Hidden & Params & Hidden & Params & Hidden & Params \\
        \midrule
        p120   &     $120$    & $[3, 5, 6]$    & $122$      & $[1, 1, 1]$  & $117$      & $[3, 5, 6]$    & $116$ \\
        p400   &     $400$    & $[6, 10, 14]$  & $368$      & $[2, 2, 2]$  & $426$      & $[6, 10, 14]$  & $358$ \\
        p1200  &   $1{,}200$  & $[10, 18, 32]$ & $1{,}098$  & $[2, 3, 5]$  & $1{,}081$  & $[10, 18, 32]$ & $1{,}090$ \\
        p3000  &   $3{,}000$  & $[22, 36, 44]$ & $2{,}940$  & $[4, 5, 8]$  & $3{,}005$  & $[22, 36, 44]$ & $2{,}896$ \\
        \bottomrule
    \end{tabular}
\end{table*}

All three closures share the pointwise architecture of \cref{sec:construction}: a
$1 \times 1 \times 1$ convolution evaluated independently at each grid point,
with three hidden layers, each followed by a $\operatorname{GELU}$
activation, and a final linear layer with no activation and no bias.
Rather than a single size, each architecture is instantiated on a capacity grid of
four target parameter counts (\cref{tab:architecture}), varied by hidden width alone
with the depth fixed at three layers; within each tier the three architectures are
roughly parameter-matched, isolating the effect of the symmetry treatment from raw
network capacity.
The forced comparison of \cref{sec:results-forced} fixes the size at the ``p3000''
tier ($\approx 3000$ parameters), while the capacity study of
\cref{sec:results-capacity} sweeps the whole grid.
The MLP and G-CNN both take the nine components of the normalized VGT
$\bar{A} / |\bar{A}|$ as input.
The MLP returns six components, interpreted directly as the symmetric stress;
the G-CNN returns the full nine-component tensor, which is then symmetrized to
six components.
In both cases the trace is subsequently removed to obtain the deviatoric closure;
the output dimensions in \cref{tab:architecture} are those of the network proper,
before this symmetrization and trace removal.
For the G-CNN the hidden widths count regular-representation channels of
dimension $|G| = 48$ (\cref{sec:groupconv}); at the ``p3000'' tier, for instance,
the three hidden layers of widths $4, 5, 8$ carry $192, 240, 384$ scalar
components.
The reported parameter count is that of the reduced, pre-synthesis weights actually
updated by the optimizer, which the synthesis $\mathcal{S}$ (\cref{sec:groupconv})
expands into the full equivariant weights at each gradient step.
Counting these degrees of freedom, $\mathcal{S}$ maps the $9$ learnable weights
of a boundary layer to a full block of $9 \cdot 48 = 432$ entries, and the $48$
learnable weights of a hidden layer to a full block of $48 \cdot 48 = 2304$
entries; we scale $\mathcal{S}$ by $1 / \sqrt{|G|}$ so that its columns are
orthonormal, preserving the standard weight-initialization variance.
The TBNN instead takes the five invariants $\lambda^*$ and returns the seven
tensor-basis coefficients $\alpha_k$, which are subsequently expanded against the
tensor basis (\cref{sec:tbnn}).
Every architecture is additionally trained in a $+\mathrm{Re}$ variant that
receives the standardized logarithm of the global filter-scale Reynolds number
as an extra input
(\cref{sec:redelta-input}), whose role in generalizing across Reynolds number is
examined in \cref{sec:results-reynolds}.

Each learned closure regresses the \emph{normalized} target
$\tilde{\tau} \coloneq \tau / (\Delta^2 |\bar{A}|^2)$ from the normalized gradient
$\bar{A} / |\bar{A}|$ (\cref{sec:construction}); for the TBNN this is the stress
reconstructed from its predicted tensor-basis coefficients (\cref{sec:tbnn}).
We train it by minimizing the relative squared error aggregated over each random
mini-batch $B$,
\begin{equation}
    L_\theta \coloneq \frac{\sum_{u \in B} \| m_\theta(\bar{u}) - \tilde{\tau}(u) \|^2}
                           {\sum_{u \in B} \| \tilde{\tau}(u) \|^2},
\end{equation}
where $\theta$ denotes the model parameters.
The training pool is the union of the $72$ gradient--stress fields from the three
forced training runs crossed with the three training filter ratios
(\ref{sec:datagen}).
Because the closures are pointwise, each grid point is an independent training
sample; snapshots are folded along one spatial axis into the batch dimension,
giving a mini-batch size of $20$ in this folded layout, so each batch mixes points
drawn from across the pooled viscosities, filter ratios, and snapshots.
The last $20\%$ of the pooled snapshot sequence is held out and monitored
only to trace the convergence curve: training is fixed-budget (the networks
are tiny), with no early stopping and no checkpointing, so this holdout never
selects the returned parameters.
We optimize the loss for $20$ epochs using AdamW with default moments
$(\beta_1, \beta_2) = (0.9, 0.999)$ and weight decay $\lambda = 0$ (reducing it
to Adam), gradient-norm clipping at $1.0$, and a linear-warm-up--cosine
learning-rate schedule ($5\%$ warm-up to the peak rate $10^{-3}$, then cosine
decay to zero over the remaining steps).
Each learned family is trained with five independent initialization seeds, and
scalar metrics are reported as mean $\pm$ one standard deviation over seeds
(\cref{sec:results}).
The networks are trained in single precision while the solver and
data pipeline stay in double precision; parameters are upcast to
double precision before LES inference.

The loss evolution is shown in \cref{fig:training}, one panel per architecture and
one curve per capacity tier.
All three architectures reach a low validation loss within a few hundred
iterations, and the remaining epochs mainly confirm convergence.
For the equivariant G-CNN and TBNN the curves for the different tiers nearly
coincide, since even the smallest networks already sit close to the floor,
whereas the unconstrained MLP improves visibly with size and reaches the same
level only at the larger tiers; this is the training-time signature of the
saturation behavior of \cref{sec:results-capacity}.
The seed bands are narrowest for the constrained models, the training-time
counterpart of the seed spread visible below saturation in \cref{fig:saturation}.
That the final loss plateaus well above zero, independently of architecture and
size, is itself informative; we interpret it in \cref{sec:discussion}.

\section{Dynamic Smagorinsky model} \label{sec:dynamic-smagorinsky}

The Smagorinsky closure~\cite{smagorinskyGeneralCirculationExperiments1963}
approximates the deviatoric SFS by
\begin{equation} \label{eq:smag}
    m_{i j}(\bar{u}, \bar{\Delta}) \coloneq -2 c \bar{\Delta}^2 | \bar{S} | \bar{S}_{i j},
\end{equation}
where $\bar{S}_{i j} \coloneq (\partial_i \bar{u}_j + \partial_j \bar{u}_i) / 2$
is the resolved strain-rate tensor,
$| \bar{S} | \coloneq \sqrt{\bar{S}_{i j} \bar{S}_{i j}}$ is its Frobenius norm,
$\bar{\Delta}$ is the LES filter width, and
$c \ge 0$ is a scalar coefficient that the dynamic procedure determines from
the resolved field at every time step, rather than fixing it a priori.
The conventional factor $\sqrt{2}$ in the strain-rate magnitude
$\sqrt{2 \bar{S}_{i j} \bar{S}_{i j}}$ is absorbed into $c$ and therefore does
not affect the predicted stress.

\subsection{Germano--Lilly procedure}

Let $\widetilde{(\cdot)}$ denote a test filter of width
$\widetilde{\Delta} \coloneq 2 \bar{\Delta}$.
We use the same Gaussian kernel as the grid filter (see \cref{sec:datagen}),
so the composition $\widetilde{\overline{(\cdot)}}$ is also Gaussian with
effective width
\begin{equation} \label{eq:dyn-compwidth}
    \widetilde{\bar{\Delta}} \coloneq \sqrt{\widetilde{\Delta}^2 + \bar{\Delta}^2}
    = \sqrt{5} \, \bar{\Delta}.
\end{equation}
The Germano identity~\cite{germanoDynamicSubgridscaleEddy1991}
expresses the SFS at the combined filter level in terms
of the SFS at the grid level and a resolved tensor
$L_{i j}$ that depends only on the LES velocity,
\begin{equation} \label{eq:dyn-L}
    L_{i j} \coloneq \widetilde{\bar{u}_i \bar{u}_j} - \widetilde{\bar{u}}_i \widetilde{\bar{u}}_j.
\end{equation}
Substituting the Smagorinsky form~\eqref{eq:smag} on both sides and assuming
that $c$ varies slowly compared with the test-filter width gives
$L_{i j} \approx c M_{i j}$, where
\begin{equation} \label{eq:dyn-M}
    M_{i j} \coloneq 2 \bar{\Delta}^2 \, \widetilde{| \bar{S} | \bar{S}_{i j}}
    - 2 \widetilde{\bar{\Delta}}^2 \, | \widetilde{\bar{S}} | \widetilde{\bar{S}}_{i j}.
\end{equation}

The trace of $L$ is generally non-zero, while $M$ is trace-free by
construction (because $\operatorname{tr} \bar{S} = 0$ for incompressible
flow).
We therefore project $L$ onto its deviatoric part before contracting it with
$M$:
\begin{equation}
    L^{\text{dev}}_{i j} \coloneq L_{i j} - \frac{1}{D} L_{k k} \delta_{i j}.
\end{equation}

\subsection{Volume-averaged coefficient}

The flow generated in \cref{sec:datagen} is statistically homogeneous on the
periodic box, so we determine a single coefficient per time step by
averaging over the LES grid (Lilly's
procedure~\cite{lillyProposedModificationGermano1992}),
\begin{equation} \label{eq:dyn-c}
    c \coloneq \max \! \left(
        \frac{\langle M_{i j} L^{\text{dev}}_{i j} \rangle_\Omega}
             {\langle M_{k \ell} M_{k \ell} \rangle_\Omega},
        \, 0 \right),
\end{equation}
where $\langle \cdot \rangle_\Omega$ denotes the spatial mean over the
periodic box.
The averaging is the least-squares fit of $c M$ to $L^{\text{dev}}$ that
minimizes $\| L^{\text{dev}} - c M \|_{L^2(\Omega)}^2$.
Clipping at zero prevents the model from injecting energy back into the
resolved scales, a common stabilization for the Germano--Lilly procedure.
The closure stress is then $m_{i j}(\bar{u}, \bar{\Delta})$
from~\eqref{eq:smag}, evaluated with this coefficient $c$.

All filters, products, and strain rates in~\eqref{eq:dyn-L} and~\eqref{eq:dyn-M}
are evaluated on the LES grid:
products of physical-space fields are pseudo-spectrally dealiased with the
$2/3$-rule (\cref{sec:pseudo-spectral}), and the Gaussian filter is applied
in spectral space.

\section{Taylor--Green errors at the wider filter ratios} \label{sec:tgv-tables}

\Cref{sec:results-tgv} reports the decaying Taylor--Green vortex at the
interpolated filter ratio $\Delta / h = 2.5$ (\cref{tab:errors-tgv}), the ratio
that matches the forced operating point of \cref{sec:results-forced}.
For completeness, \cref{tab:errors-tgv-d35,tab:errors-tgv-d50} give the same
aggregate errors at the two wider ratios, $\Delta / h = 3.5$ and $5.0$, whose
filter-scale Reynolds numbers ($\operatorname{Re}_\Delta \approx 500$ and $813$)
carry the Taylor--Green flow into the upper end of the forced grid. The trends
discussed in the text are visible here: the a priori tensor error stays
essentially flat across the filter ratios, while the over-dissipation ratio falls
as the filter widens and the $+\mathrm{Re}$ input, which corrects the tight filter,
over-corrects at $\Delta / h = 5.0$.

\begin{table*}
    \centering
    \caption{%
        Aggregate errors for the decaying Taylor--Green vortex at nominal integral
        $\operatorname{Re} = 6000$ and filter ratio $\Delta / h = 3.5$; columns and
        conventions as in \cref{tab:errors-tgv} (reference backscatter fraction
        $0.279$).
    }
    \label{tab:errors-tgv-d35}
    \begin{tabular}{l l l l l l l}
    \toprule
    Model
                & Tier
                & Closure \eqref{eq:tensor-error-prior}
                & Cross-corr.
                & Solution \eqref{eq:tensor-error-post}
                & Median diss.
                & Backscatter \\
    \midrule
    No-model    & --    & $1.0000$            & N.A.                & $0.6005$            & $0.000$           & $0.0000$ \\
    Dyn. Smag.  & --    & $0.9674$            & --                  & $0.5506$            & $1.597$           & $0.0025$ \\
    Clark       & --    & $0.4647$            & $0.8821$            & --                  & $0.559$           & $0.3076$ \\
    MLP         & p3000 & $0.4634 \pm 0.0004$ & $0.8763 \pm 0.0002$ & $0.4870 \pm 0.0042$ & $1.680 \pm 0.011$ & $0.1252 \pm 0.0012$ \\
    MLP+Re      & p3000 & $0.4634 \pm 0.0007$ & $0.8765 \pm 0.0005$ & $0.5083 \pm 0.0089$ & $1.708 \pm 0.007$ & $0.1240 \pm 0.0013$ \\
    G-CNN       & p3000 & $0.4623 \pm 0.0002$ & $0.8768 \pm 0.0001$ & $0.4786 \pm 0.0055$ & $1.729 \pm 0.006$ & $0.1228 \pm 0.0005$ \\
    G-CNN+Re    & p3000 & $0.4608 \pm 0.0006$ & $0.8778 \pm 0.0002$ & $0.4786 \pm 0.0052$ & $1.744 \pm 0.007$ & $0.1225 \pm 0.0010$ \\
    TBNN        & p3000 & $0.4575 \pm 0.0002$ & $0.8792 \pm 0.0001$ & $0.4931 \pm 0.0062$ & $1.687 \pm 0.007$ & $0.1152 \pm 0.0009$ \\
    TBNN+Re     & p3000 & $0.4563 \pm 0.0003$ & $0.8797 \pm 0.0001$ & $0.4805 \pm 0.0067$ & $1.693 \pm 0.006$ & $0.1113 \pm 0.0019$ \\
    \bottomrule
\end{tabular}

\end{table*}

\begin{table*}
    \centering
    \caption{%
        Aggregate errors for the decaying Taylor--Green vortex at nominal integral
        $\operatorname{Re} = 6000$ and filter ratio $\Delta / h = 5.0$; columns and
        conventions as in \cref{tab:errors-tgv} (reference backscatter fraction
        $0.265$).
    }
    \label{tab:errors-tgv-d50}
    \begin{tabular}{l l l l l l l}
    \toprule
    Model
                & Tier
                & Closure \eqref{eq:tensor-error-prior}
                & Cross-corr.
                & Solution \eqref{eq:tensor-error-post}
                & Median diss.
                & Backscatter \\
    \midrule
    No-model    & --    & $1.0000$            & N.A.                & $0.6103$            & $0.000$           & $0.0000$ \\
    Dyn. Smag.  & --    & $0.9643$            & --                  & $0.5587$            & $1.156$           & $0.0014$ \\
    Clark       & --    & $0.4902$            & $0.8696$            & --                  & $0.459$           & $0.3112$ \\
    MLP         & p3000 & $0.4839 \pm 0.0003$ & $0.8672 \pm 0.0002$ & $0.5399 \pm 0.0052$ & $1.343 \pm 0.008$ & $0.1118 \pm 0.0013$ \\
    MLP+Re      & p3000 & $0.4981 \pm 0.0014$ & $0.8624 \pm 0.0004$ & $0.5363 \pm 0.0201$ & $1.571 \pm 0.006$ & $0.0975 \pm 0.0014$ \\
    G-CNN       & p3000 & $0.4820 \pm 0.0002$ & $0.8680 \pm 0.0001$ & $0.5156 \pm 0.0023$ & $1.377 \pm 0.005$ & $0.1091 \pm 0.0006$ \\
    G-CNN+Re    & p3000 & $0.4943 \pm 0.0011$ & $0.8648 \pm 0.0003$ & $0.4991 \pm 0.0038$ & $1.613 \pm 0.006$ & $0.0901 \pm 0.0004$ \\
    TBNN        & p3000 & $0.4789 \pm 0.0002$ & $0.8696 \pm 0.0001$ & $0.5401 \pm 0.0013$ & $1.365 \pm 0.005$ & $0.0992 \pm 0.0010$ \\
    TBNN+Re     & p3000 & $0.4902 \pm 0.0003$ & $0.8662 \pm 0.0002$ & $0.4994 \pm 0.0039$ & $1.591 \pm 0.001$ & $0.0856 \pm 0.0012$ \\
    \bottomrule
\end{tabular}

\end{table*}

\bibliographystyle{elsarticle-num}
\bibliography{data/references.bib}

\end{document}